%% file: mu.tex
\documentclass[12pt]{amsart}

\usepackage{amsmath,amsfonts,amssymb,amscd,amsthm,amsbsy,epsf}
\newtheorem{thm}{Theorem}[section]
\newtheorem{lemma}[thm]{Lemma}
\newtheorem{cor}[thm]{Corollary}

\newtheorem{claim}[thm]{Claim}
\newtheorem{subclaim}[thm]{Subclaim}

\newtheorem*{TP}{Theorem~\ref{tp}}
\newtheorem*{HandleLemma}{Handle Addition Lemma}
\newtheorem*{CORPJH}{Corollary~\ref{corPJH}}
\newtheorem*{con3.19}{Conjecture~\ref{cor:restriction}}
\newtheorem*{Tsolution}{Proposition~\ref{prop:5.1}}

\newtheorem{prop}[thm]{Proposition}
\theoremstyle{definition}
\newtheorem{assump}[thm]{Assumption}
\newtheorem*{remark}{Remark}
\newtheorem{defn}{Definition}
\numberwithin{figure}{section}
\numberwithin{equation}{section}
\def\QED#1{\newline{\rightline{Q.E.D.~(#1)\medskip}}}
\def\strut{\hbox{\vrule height 1.35em depth.65em width0pt}}
\def\myfig#1#2#3{$$\vbox{
	\centerline{\epsfysize=#1truein\epsfbox{#2.eps}}
	\centerline{\sc Figure #3}}$$}
\def\nonumfig#1#2{$$\epsfysize=#1truein\epsfbox{#2.eps}$$}
\def\doublefig#1#2#3#4#5#6{$$\vbox{
	\centerline{{\epsfysize=#1truein\epsfbox{#2.eps}}
	\hskip1.0truein{\epsfysize=#4truein\epsfbox{#5.eps}}}
	\centerline{\sc Figure #3\hskip1.8truein Figure #6}}$$}
\def\sdoublefig#1#2#3#4#5#6{$$\vbox{
	\centerline{{\epsfysize=#1truein\epsfbox{#2.eps}}
	\hskip.5truein{\epsfysize=#4truein\epsfbox{#5.eps}}}
	\centerline{\sc Figure #3\hskip2.2truein Figure #6}}$$}
\def\wdoublefig#1#2#3#4#5#6{$$\vbox{
	\centerline{{\epsfysize=#1truein\epsfbox{#2.eps}}
	\hskip1.5truein{\epsfysize=#4truein\epsfbox{#5.eps}}}
	\centerline{\sc Figure #3\hskip1.8truein Figure #6}}$$}
\def\T{{\mathcal T}}
\def\J{{\mathcal J}}

\def\zed{{\mathbb Z}}

\def\PJHsolutionUnmarked{0.1}
\def\actionofCre{1.2}
\def\enhancer{1.1}
\def\recombdirknot2seq1a{1.3}
\def\CreAm{1.4}
\def\tableA{1.5}
\def\eightfigs{1.6}
\def\CreB{1.7}
\def\111{1.8}
\def\branched{1.9}

\def\projcircle{2.1}
\def\normaleq{2.2}
\def\unique{2.3}
\def\tangcircA{2.4}
\def\tangcircBm{2.5} 
\def\tangleaddclose{2.6}
\def\tangleEqns{2.7}
\def\toruslink{2.8}
\def\convention{2.9}
\def\bandingm{2.10}
\def\intrans{2.11}
\def\unknotHopf{2.12}
\def\synapsetangle{2.13}
\def\deletionreaction{2.14}
\def\solutiontangle{2.15}

\def\standard{3.1}
\def\wagonwheel{3.2}
\def\wwgraph{3.3}
\def\nonplanar{3.4}
\def\NG{3.5}
\def\carriedbyG{3.6}
\def\solutiontangleA{3.7}
\def\tangleG{3.8}
\def\Infinitely{3.9}
\def\HALemfig{3.12}
\def\GstandardB{3.15}
\def\splitrationaltang{3.13}
\def\HatP{3.16}
\def\Betas{3.17}
\def\fouredges{3.11}
\def\counter{3.10} 
\def\tripusA{3.14}    
\def\aaa{3.18}

\def\LemfigA{4.1}
\def\LemfigB{4.2}
\def\LemfigC{4.3}
\def\LemfigD{4.4}
\def\lifesaverA{4.5}
\def\lifesaverB{4.6}
\def\hockeyballs{4.7}
\def\parallelstrandsC{4.8}

\def\Setup2{5.1}
\def\IIsetup{5.2}
\def\Reduceone{5.3}
\def\Reducetwo{5.4}
\def\IIredone{5.5}
\def\figure105{5.6}
\def\IIredtwo{5.7}
\def\Ioneset{5.8}
\def\Ionered{5.9}
\def\Itwosetup{5.10}
\def\Itworeduce{5.11}
\def\Untwist{5.12}
\def\Setupthree{5.13}
\def\threeredone{5.14}
\def\threeredtwo{5.15}
\def\threeredthree{5.16}
\def\Setupfour{5.17}
\def\Setupdl{5.18}
\def\SetupFT{5.19}
\def\SetupFell{5.20}
\def\SetupFF{5.21}
\begin{document}
\title{
Tangle analysis of difference topology experiments: applications to a Mu protein-DNA complex
%%Tangle analysis of the Mu transpososome protein complex
%%which binds three DNA segments
}
\author{Isabel K. Darcy, John Luecke, and Mariel Vazquez}


%\dedicatory{Here is a dedication}
 \begin{abstract} 
 We develop topological methods for analyzing difference topology experiments involving 
3-string tangles.  Difference topology is a novel technique used to unveil the structure of 
stable protein-DNA complexes involving two or more DNA segments. We analyze such 
experiments for the Mu protein-DNA complex. We characterize the solutions to the 
corresponding tangle equations by certain knotted graphs. By investigating planarity 
conditions on these graphs we show that there is a unique biologically relevant solution. 
That is, we show there is a unique rational tangle solution, which is also the unique 
solution with small crossing number.
 \end{abstract}

\maketitle
\baselineskip=18pt
\pagestyle{plain}              % page nos. at bottom, no headline

\input section0.tex

\input section1.tex 
\input section2.tex

\input section3.tex

\input section4.tex

\input section5.tex

%%%%%%%%%%%%%%%%%%%%%%%%%%%%%%%%%%%%%%%%%%%%%%%%%%%%%%%

\eject

~\break
Isabel K. Darcy, 
Department of Mathematics,
University of Iowa,
Iowa City, IA 52242 \hfil \break
E-mail:  idarcy@math.uiowa.edu \hfil \break
URL: http://www.math.uiowa.edu/$\sim$idarcy

~\break
John Luecke, Department of Mathematics, University of Texas at Austin, Austin, TX 78712 \hfil \break
E-mail: luecke@math.utexas.edu  \hfil \break
URL: http://www.ma.utexas.edu/$\sim$luecke/

~\break
Mariel Vazquez, San Francisco State University
Department of Mathematics, 
1600 Holloway Ave,
San Francisco, CA 94132
 \hfil \break
E-mail: mariel@math.sfsu.edu \hfil \break
URL: http://math.sfsu.edu/vazquez/\end{document}

%% file: section0.tex
In \cite{PJH}, Pathania {\it et al} determined the shape of DNA bound within the Mu 
transposase protein complex using an experimental technique called difference topology 
\cite{HJ, KBS, GBJ,PJH, PJH2, YJPH, YJH} and by making certain assumptions regarding the DNA 
shape. We show that their most restrictive assumption
(the plectonemic form described near the end of section 1) 
is not needed, and in doing so, conclude that the only biologically reasonable solution
for the shape of DNA bound by Mu transposase is the one they found \cite{PJH}
 (Figure \PJHsolutionUnmarked).  We will call this 3-string tangle the PJH
solution. The 3-dimensional ball
represents the protein complex, and the arcs represent the bound DNA.
The Mu-DNA complex modeled by this tangle is called the Mu transpososome

\myfig{1.25}{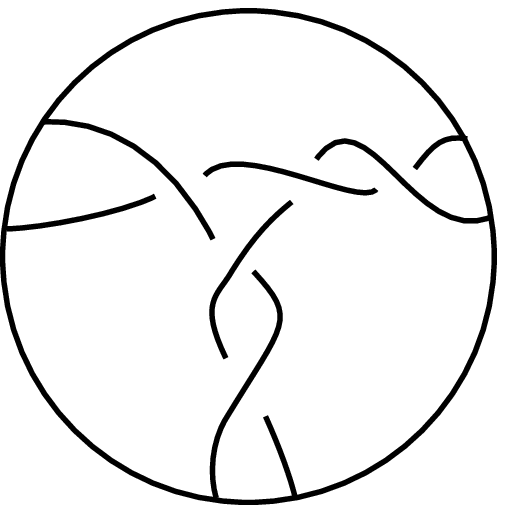}{\PJHsolutionUnmarked}

In section 1 we provide some biological background and describe eight difference topology experiments from 
\cite{PJH}. In section 2, 
we translate the biological problem of determining the shape of DNA bound by Mu into a mathematical model.  The 
mathematical model consists of a system of ten 3-string tangle equations (Figure~\normaleq). Using 2-string tangle 
analysis, we simplify this to a system of four tangle equations (Figure~\solutiontangle). In section 3 we 
characterize solutions to these tangle
		equations in terms of knotted graphs. This allows us to exhibit infinitely many different 3-string 
tangle solutions.
The 
existence of solutions different from the PJH solution raises the possibility of alternate acceptable models. In 
sections 3 - 5, we 
show that all solutions to the mathematical problem other than the PJH solution are too complex to be 
biologically reasonable, where 
the complexity is measured either by the rationality or by the minimal crossing number of the 3-string tangle 
solution.

In section 3, we show that the only {\em rational\/} solution is the PJH solution.  In particular we prove the 
following corollary.

\begin{CORPJH}
Let $\T$ be a solution tangle. 
If $\T$ is rational or split or if $\T$ has parallel strands, 
then $\T$ is the PJH solution.  
\end{CORPJH}

In section 4 we show that any 3-string tangle with fewer than 8 crossings, up to free isotopy ({\it i.e.} allowing 
the ends of the tangle to move under the
                isotopy),  must be 
either split or have parallel strands.  Thus Corollary \ref{corPJH} 
implies that any 
solution $\T$ different from the PJH solution must have at least 8 
crossings up to free isotopy. Fixing the framing of a solution tangle (the
                normal framing of section 2), and working in the category
                of tangle equivalence -- {\it i.e.} isotopy fixed on the boundary --
                we prove the following lower bound on the crossing number of
                exotic solutions:

\begin{Tsolution}    
Let $\T$ be an {in trans} solution tangle. 
If $\T$ has a projection with fewer than 10 crossings, then $\T$ 
is the PJH  tangle.
\end{Tsolution}

The framing used in \cite{PJH} is different than our normal framing.  In the context of \cite{PJH}, Proposition 
~\ref{prop:5.1} says that if Mu binds fewer than 9 crossings, then the PJH solution  is the
only solution fitting the experimental data. The PJH solution has 5 crossings. 
We interpret
Corollary~\ref{corPJH} and Proposition~\ref{prop:5.1} as
saying that the PJH solution is the only biologically reasonable model for the Mu transpososome.

Although we describe 8 experiments from \cite{PJH}, in the interest of minimizing
lab time, 
we show that only 3 
experiments ({\em in cis} deletion) are needed to prove the main result of sections 3 (Corollary 3.20) and 4.
A fourth experiment ({\em in trans} deletion) allows us to rule out some solutions (section 3) and is 
required for the analysis in 
section 5.  The remaining four experiments (inversion) were used in model design \cite{PJH}.

The results in 
sections 2-4 extend to cases such as \cite{YJPH, YJH} where
the experimental products are $(2, L)$ torus links \cite{HS} or the trefoil knot 
\cite{KMOS}.
The work in section 3 involves the analysis of knotted graphs
as in \cite{G}, \cite{ST}, \cite{T}. As a by-product, we prove the following:

\begin{TP}
Suppose $\widehat{G}$  is a tetrahedral graph with the following properties:
\begin{enumerate}
\item  There exists three edges $e_1$, $e_2$, $e_3$ such that $\widehat{G} - e_i$ is planar.
\item  The three edges $e_1$, $e_2$, $e_3$ share a common vertex.
\item  There exists two additional edges, $b_{12}$ and $b_{23}$ such that $X(\widehat{G} - b_{12})$ and $X(\widehat{G} - b_{23})$ have compressible boundary.
\item  $X(\widehat{G})$  has compressible boundary.
\end{enumerate}
Then $\widehat{G}$  is planar.
\end{TP}

 In subsection 3.4 we also give several examples of non-planar 
tetrahedral graphs to show that none of the hypothesis in Theorem \ref{tp} 
can be eliminated.

{\bf Acknowledgments}

We would like to thank D. Buck, R. Harshey, S. Pathania, and C. Verjovsky-Marcotte for helpful comments.  We would particularly like to thank M. Jayaram for many helpful discussions.
We also
thank M. Combs for numerous figures, R. Scharein and Knotplot.com for assistance with Figure \111, and A. Stasiak, University of Lausanne,   for the
electron micrograph of supercoiled DNA in Figure \branched.

J.L. would like to thank the Institute for Advanced Study
in Princeton for his support as a visiting member.  
This research was supported in part by the Institute for Mathematics and its Applications with funds provided by 
the National Science Foundation (I.D. and M.V.);
by a grant from the Joint DMS/NIGMS Initiative to Support Research in the Area of Mathematical Biology (NIH GM 
67242)
to I.D.; and by NASA NSCOR04-0014-0017, 
by MBRS SCORE S06 GM052588, and by NIH-RIMI Grant NMD000262 to M.V.

%% file: section1.tex
\section{Biology Background and Experimental Data}

Transposable elements, also called mobile elements, are fragments of DNA able to 
move along a genome by a process called transposition. Mobile elements play an 
important role in the 
shaping of a genome \cite{DMBK,S}, and they can impact the health of an organism 
by introducing genetic mutations. Of special interest is that transposition is 
mechanistically very similar to the way 
certain retroviruses, including HIV, integrate into their host genome.

Bacteriophage Mu is a system widely used in transposition studies due to the 
high efficiency of  Mu transposase (reviewed in \cite{CH}). The MuA protein 
performs the first steps required to transpose the 
Mu genome from its starting location to a new DNA location. MuA binds to 
specific DNA sequences which we refer to as attL and attR sites (named after 
{\it Left} and {\it Right} attaching regions). A third DNA sequence called the 
enhancer (E) is also required to assemble the Mu transpososome.  The Mu 
transpososome is a very stable complex consisting of 3 segments of double-
stranded DNA captured in a protein complex \cite{BM,MBM}. In this paper we are interested in 
studying the topological structure of the DNA within the Mu transpososome.

 \subsection{Experimental design}

We base our study on the difference topology experiments of \cite{PJH}. In this 
technique, circular DNA is first incubated with the 
protein(s){\footnote{Although we use the singular form of protein instead of the 
plural form, most protein-DNA complexes involve several proteins.  For example 
formation of the Mu transpososome involves four MuA proteins and the protein HU.  
Also since these are test tube reactions and not single molecule experiments, 
many copies of the protein are added to many copies of the DNA substrate to form 
many complexes.}} under study (in this case, MuA), which bind DNA.
A second protein whose mechanism is well understood is added to the reaction (in 
this case Cre).  This second protein is a protein that can cut DNA and change 
the circular DNA topology before resealing the break(s), resulting in knotted or 
linked DNA.  
DNA crossings bound by the first protein will affect the product topology.  
Hence one can gain information about the DNA conformation bound by the first 
protein by determining the knot/link type of the DNA knots/links produced by the 
second protein.

  In the experiments, first circular unknotted DNA is created containing the 
three binding sites for the Mu transpososome (attL, attR, E) and two binding sites 
for Cre (two $loxP$ sites).  We will refer to this unknotted DNA as {\em 
substrate}.
The circular DNA is first incubated with the proteins required for Mu 
transposition, thus forming the transpososome complex. This complex leaves three 
DNA loops free outside the transpososome (Figure \enhancer). The two $loxP$ sites are strategically 
placed in two of the three outside loops. The complex is incubated with Cre 
enzymes, which bind the $loxP$ sites, introduce two double-stranded breaks, 
recombine the loose ends and reseal them. 
A possible 2-string tangle model for the local action of Cre at these sites is shown in Figure \actionofCre~\cite{GGD}.  This cut-and-paste reaction may change 
the topology (knot/link type) of the DNA circle. Changes in the substrate's topology resulting 
from Cre action can reveal the structure within the Mu transpososome. 

\myfig{2.25}{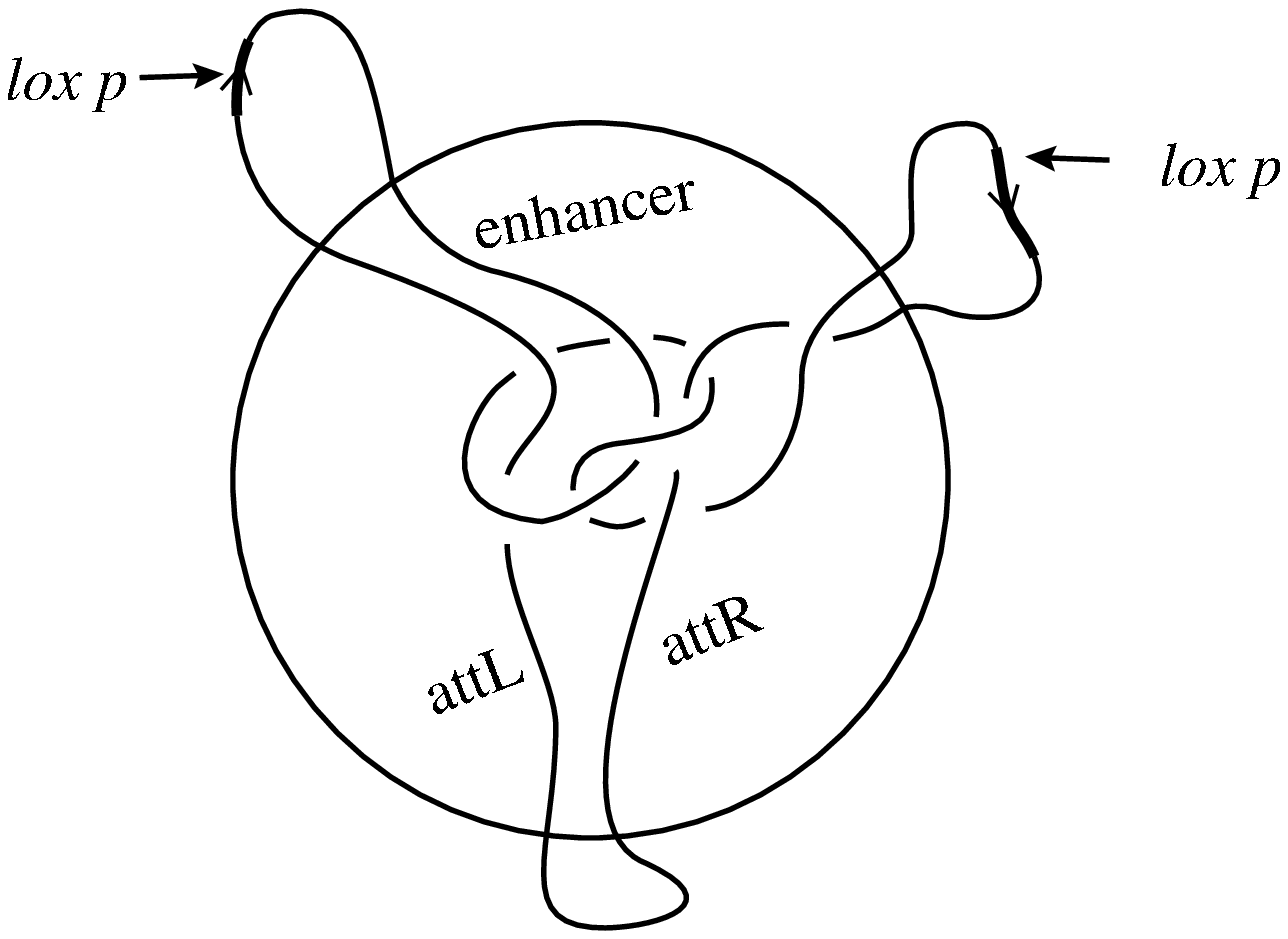}{\enhancer} 
\myfig{0.8}{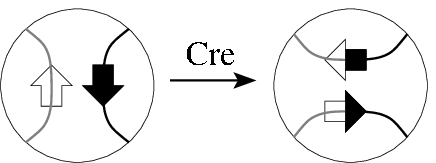}{\actionofCre} 

By looking at such topological changes, Pathania {\em et al.} \cite{PJH} deduced the 
structure of the transpososome to be the 
that of Figure ~\PJHsolutionUnmarked ~(the {\em PJH solution}).
In this paper we give a knot theoretic analysis 
that supports this deduction. We show that although there are other 
configurations that would lead to the same product topologies seen in the 
experiments, they are necessarily too complicated to be biologically reasonable.

If the orientation of both $loxP$ sites induces the same orientation on the 
circular substrate (in biological terms, the sites are {\em directly repeated}), 
then recombination by Cre 
results in a link of two components and is referred to as a {\em deletion\/} 
(Figure \recombdirknot2seq1a, left). Otherwise the sites are {\em inversely repeated}, 
the product is a knot, and the recombination is called an {\em 
inversion\/} (Figure~\recombdirknot2seq1a, right).

\myfig{1}{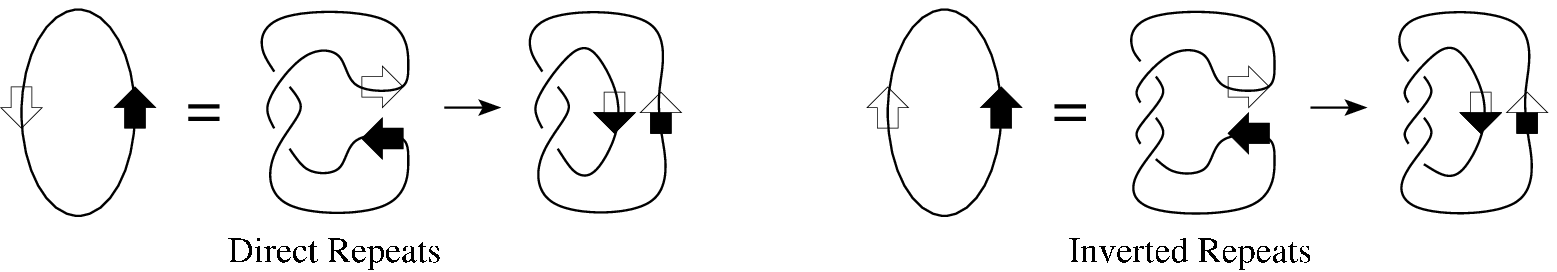}{\recombdirknot2seq1a}

In \cite{PJH} six out of eight experiments were designed by varying the relative 
positions of the $loxP$ sites and their relative orientations.  The last pair of 
experiments involved omitting one of the Mu binding sites on the circular
            substrate and placing that site 
on a linear piece of DNA to be provided ``{\em in trans\/}" as described below. 

In the first pair of experiments from \cite{PJH}, $loxP$ sites were introduced 
in the substrate on both sides of the enhancer sequence (E) (Figure~\enhancer). The 
sites were 
inserted with 
orientations to give, in separate runs, deletion and inversion products. The 
transpososome was disassembled and the knotted or linked products analyzed using 
gel electrophoresis and electron microscopy. The primary inversion products were 
$(+)$ trefoils, and the primary deletion products were 4-crossing right-hand 
torus links ($(2,4)$ torus links).

The assay was repeated, but now the $loxP$ sites were placed on both sides of 
the attL sequence. The primary products were $(2,4)$ torus links for deletion, 
and trefoils for inversion. 
In a third set of experiments, the assay was repeated again with the $loxP$ 
sites on both sides of the attR sequence. The primary products were $(2,4)$ 
torus links for deletion, and 5-crossing knots for inversion. 

Recall that in these first six experiments, the three Mu binding sites, attL, attR, and the enhancer, are all placed on the same circular DNA molecule. We will refer to these six experiments as the {\em in cis} experiments to differentiate them from the
 final set of experiments, the {\em in trans\/} assay. In the {\em in trans\/} assay, circular DNA 
substrates were created that contained attL and attR sites but no enhancer site. Each $loxP$ site was inserted between the attL and attR sites as shown in 
Figure~\CreAm. The enhancer sequence was placed on a linear DNA molecule.  
The circular substrate was incubated in solution with linear DNA 
molecules containing the enhancer sequence and with the proteins required for 
transpososome assembly.  
In this case, we say that the enhancer is provided {\em in trans}.
The $loxP$ sites in the resulting transpososome complex 
underwent Cre recombination. After the action of Cre and the disassembly of the 
transpososome (including the removal of the loose enhancer strand), the primary 
inversion products were trefoil knots, and the primary 
deletion products were (2,2) torus links (Hopf links). 

\myfig{2.0}{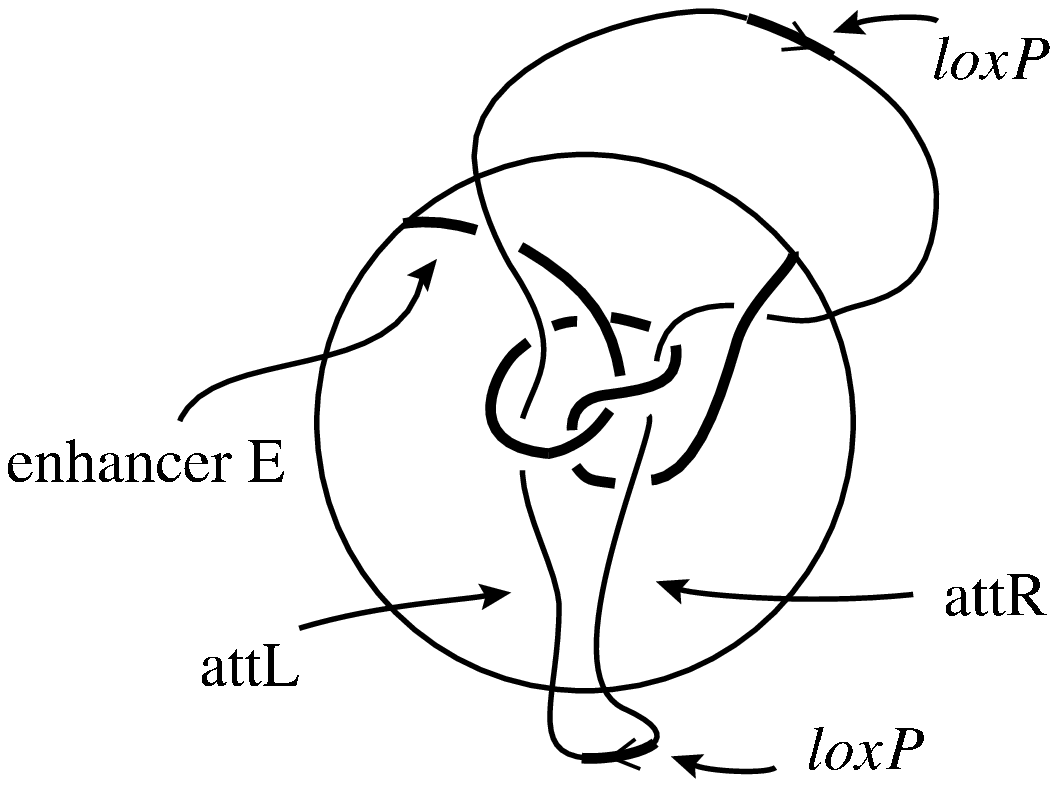}{\CreAm}

The results of these experiments are summarized in Table~\tableA\footnote{The 
chirality of the products was only determined when the $loxP$ sites were placed 
on both sides of the enhancer sequence. We here assume the chirality in 
Table~\tableA, where (2,4) torus link denotes the 4-crossing right-hand torus 
link. If any of the products are left-hand (2, 4) torus links, the results of 
sections 2 - 5 applied to these products leads to biologically unlikely 
solutions.

}. 
Vertical columns correspond to the placement of $loxP$ sites, {\it e.g.} the 
attL column shows inversion and deletion products when the $loxP$ sites were 
placed on both sides of the attL sequence.
\vskip 10pt

\centerline{$\vbox{\offinterlineskip\halign{
\enspace\hfil\hbox{#}\hfil\enspace\strut&\vrule#&
\enspace\hfil\hbox{#}\hfil\enspace&\vrule#&
\enspace\hfil\hbox{#}\hfil\enspace&\vrule#&
\enspace\hfil\hbox{#}\hfil\enspace&\vrule#&
\enspace\hfil\hbox{#}\hfil\enspace\cr
&&enhancer && attL&&attR&&{\em in trans\/}\cr
\noalign{\hrule}
Inversion&&$(+)$-trefoil&&trefoil &&5-crossing knot &&trefoil\cr
\noalign{\hrule}
Deletion&&$(2,4)$-torus 
&&$(2,4)$-torus 
&&$(2,4)$-torus 
&&$(2,2)$-torus\cr
\noalign{\vskip-6pt}
&&link&&link&&link&&link\cr
\noalign{\vskip9pt}
\multispan9\hfill\hbox{Table~\tableA
}\hfill\cr
}}$}

\vskip 10pt

\subsection{Tangle Model}
Tangle analysis is a mathematical method that models an enzymatic reaction as a 
system of tangle equations \cite{ES1, SECS}. 2-string tangle analysis has been 
successfully used to solve the topological mechanism of several site-specific 
recombination enzymes \cite{ES1, ES2, SECS, GBJ, D, VS, VCS, BV}. The Mu 
transpososome is better explained in terms of 3-string tangles. Some efforts to 
classifying rational 3-string tangles 
and solving 3-string tangle equations are underway \cite{C1, C2, EE, D1}. In 
this paper we find tangle solutions for the relevant 3-string 
tangle equations;  we characterize solutions in terms of certain knotted graphs called solution 
graphs and 
show that 
the PJH solution (Figure~\PJHsolutionUnmarked) \/ is the unique rational 
solution.

The unknotted substrate captured by the transpososome is modeled as the union of 
the two 3-string tangles $\T_0 \cup \T$, where $\T$ is the transpososome tangle 
and $\T_0$ is the tangle outside the transpososome complex. 
$\T_0 \cup \T$ is represented in Figure~\enhancer.  Notice that in this figure the 
$loxP$ sites are placed on both sides of the enhancer sequence, but the 
placement of these sites varies throughout the experiments.

Figure~\eightfigs\ shows the action of Cre on the transpososome proposed in 
\cite{PJH}. The experimentally observed products are indicated in this figure.
%% to produce the observed products.
 For 
example, $E$-inversion refers 
to the product corresponding to inversely repeated $loxP$ sites introduced on both sides of the 
enhancer sequence. However, there are 
other 3-string tangles, assuming the same action of Cre, that give rise to the 
same products. Figures~\CreB\ and ~\111~ show {two} such examples. If one replaces the tangle 
of Figure~\PJHsolutionUnmarked\  with either  that of
Figure ~\CreB~ or  ~\111~ 
in Figure~\eightfigs, the captions remain valid. 

\myfig{5.0}{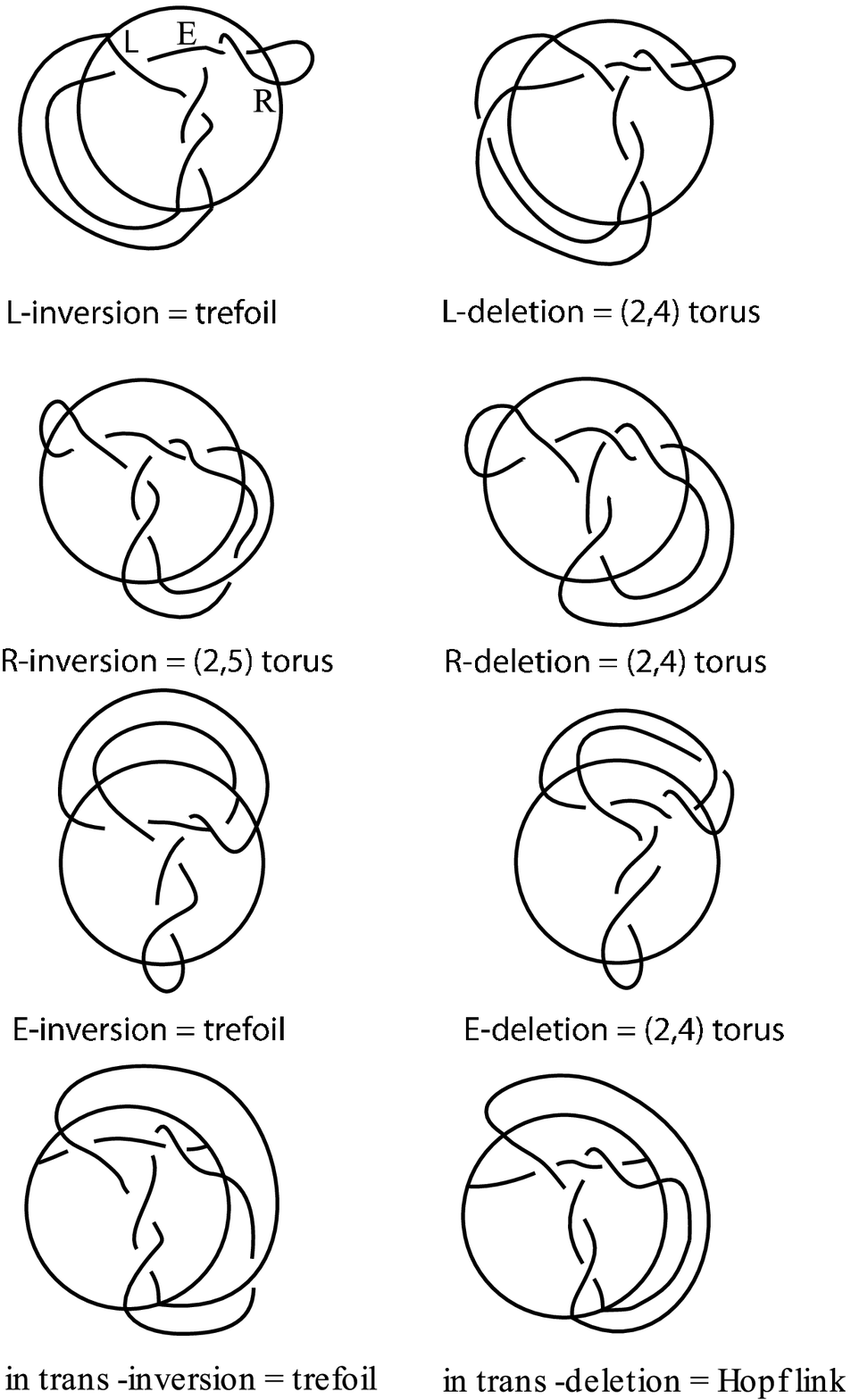}{\eightfigs}

%%\doublefig{1.2}{CreB}{\CreB}{1.08}{111}{\111}
\doublefig{1.4}{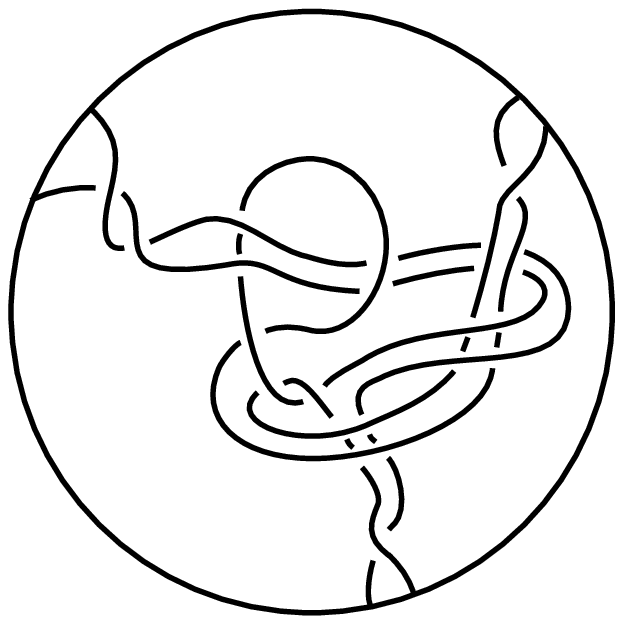}{\CreB}{1.26}{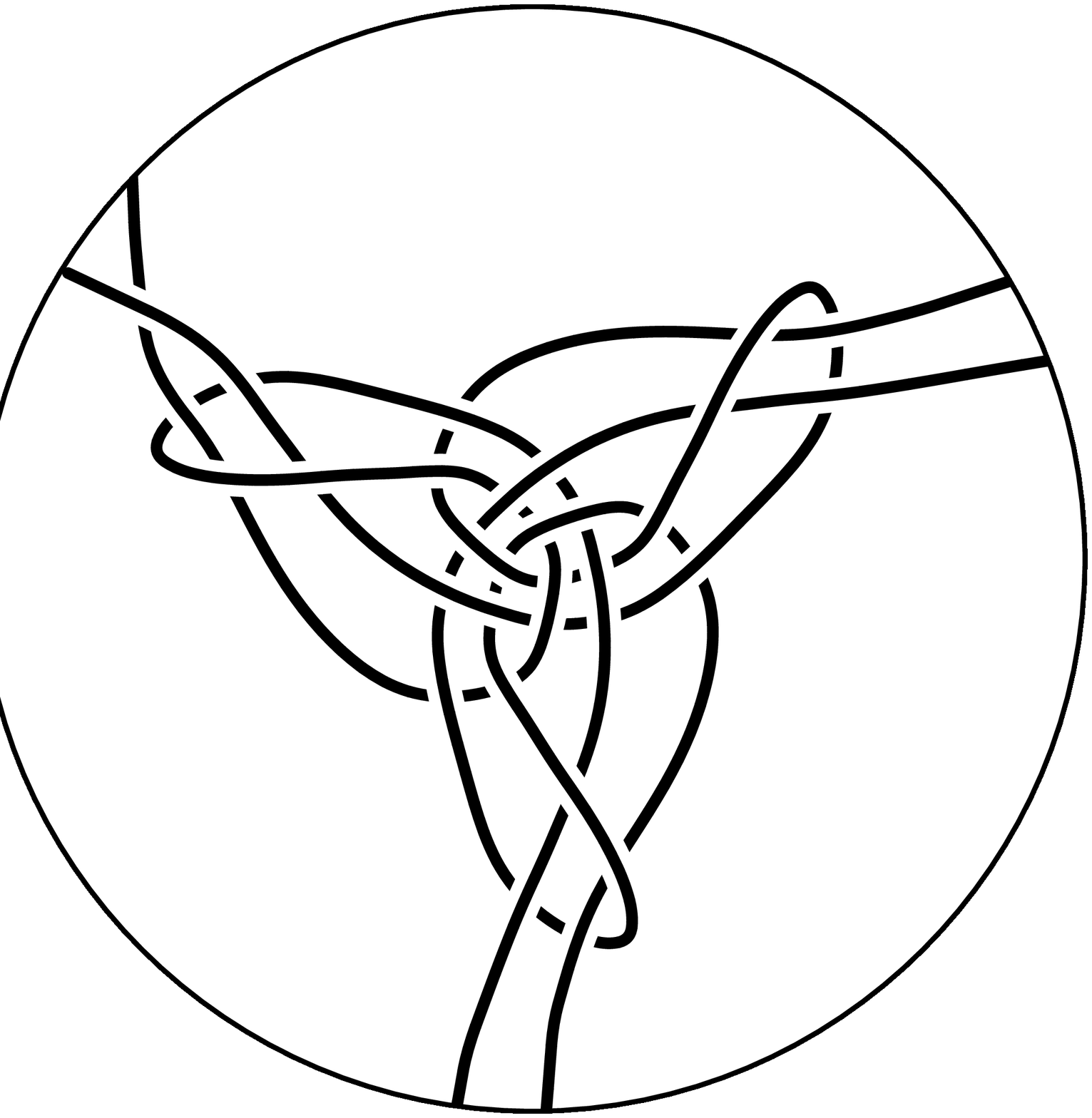}{\111}

\cite{PJH} determined the shape of DNA within the Mu transpososome to be the PJH 
solution (Figure~\PJHsolutionUnmarked) by making a restrictive assumption 
regarding this DNA conformation.  They looked at only the most biologically 
likely shape: a 3-branched supercoiled structure like that shown in 
Figure~\branched {\footnote{ Electron microscopy of supercoiled DNA courtesy of 
Andrzej Stasiak.}}.  
The $loxP$ sites were strategically{\footnote{Although we described only 8 
experiments from \cite{PJH}, they performed a number of experiments to determine 
and check effect of site placement.}} placed close to the Mu transpososome 
binding sites in order to prevent Cre from trapping random crossings not bound 
within the Mu transpososome.  
In half the experiments it was assumed that Cre trapped one extra crossing 
outside of the Mu transpososome in order to obtain the {\it loxP} sequence 
orientation of Figure \actionofCre ~(indicated by the arrows).
It was also assumed that this occurred with the higher crossing product when 
comparing inversion versus deletion products.  Hence a crossing outside of the 
Mu transpososome can be seen in Figure~\eightfigs in the case of E-deletion, L-deletion, 
R-inversion, and {\em in trans\/} inversion.  In all other cases, it was assumed that 
Cre did not trap any extra crossings outside of the Mu transpososome.  By 
assuming a branched supercoiled structure, \cite{PJH} used their experimental 
results to determine the number of crossings trapped by Mu in each of the three 
branches.  
%%In our analysis, we do not assume a branched supercoiled structure 
%%within the Mu transpososome.  
In sections 2--5 
we show that we are able 
to reach the same conclusion as \cite{PJH} without assuming a branched supercoiled structure
within the Mu transpososome.

\vskip 10pt
 
%
%\myfig{1.5}{branched}{\branched}
%
\vskip -20pt

\myfig{1.5}{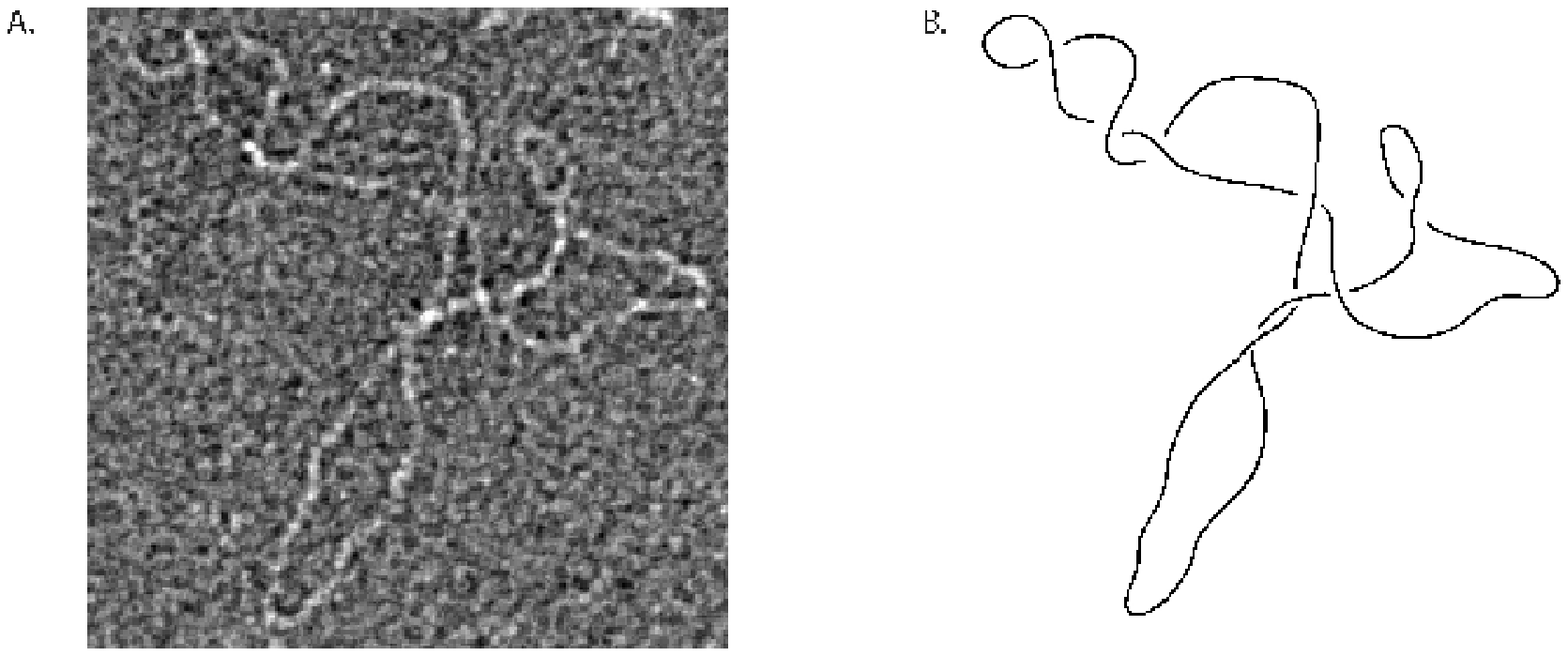}{\branched}

%% file: section2.tex
%%%%%%%%%%%%%%%%%%%%%%%%%%%%
\section{Normal Form}	%% 2

\subsection{Normal form}

The substrate  for Cre recombination in \cite{PJH} is modeled as the 3-string 
tangle union $\T\cup \T_0$. 
We here introduce a framing for $\T$ called the normal form,
which is different from that in the PJH solution (section 1).
The choice of framing affects only the arithmetic in section 2 and does not affect any of the results in sections 3 or 
4.  
The results of
section 5 on the crossing number of  $\T$ are
made with respect to this framing.

In Figure~\projcircle, let $c_1,c_2,c_3$ be the strings of $\T_0$ and 
$s_{12},s_{23},s_{31}$ be the strings of $\T$. The substrate is the union of the 
$c_i$'s and the $s_{ij}$'s. We assume there is a projection of $\T\cup \T_0$ so that 
$c_1,c_2,c_3$ are isotopic (relative endpoints) onto the tangle circle and so 
that the endpoints of $s_{ij}$ are contiguous on the tangle circle 
(Figure~\projcircle).  Note that this projection is different from  that  
 in \cite{PJH}. It is a simple matter to convert between projections, 
as described below. 
\myfig{2.0}{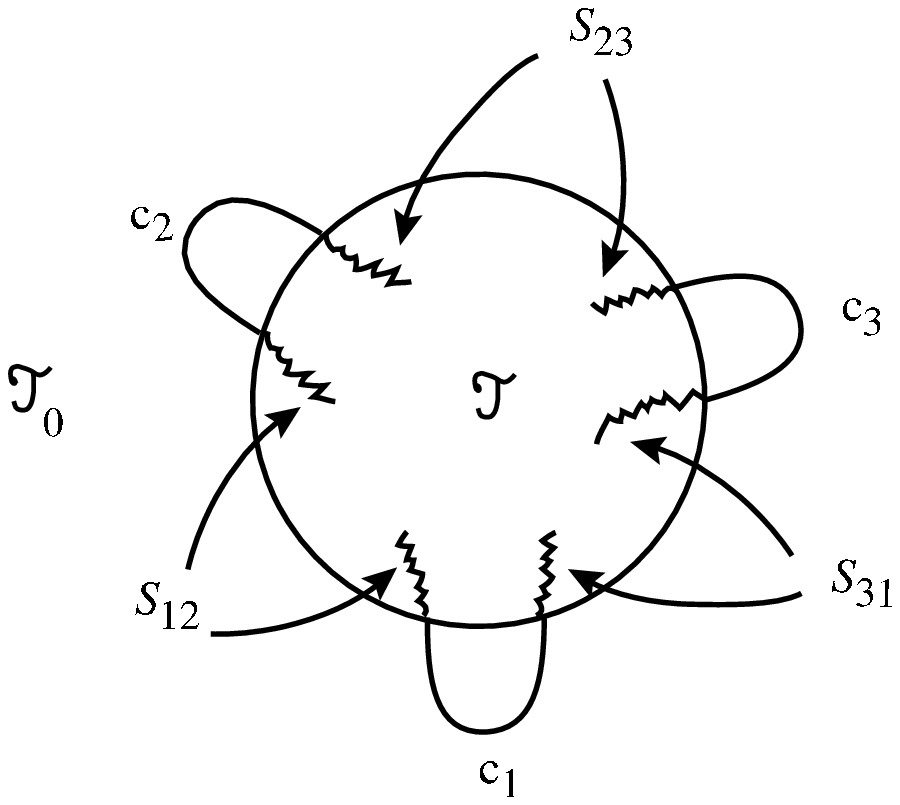}{\projcircle}

In each experiment the two recombination sites for Cre ($loxP$ sites) are located on 
two strings $c_i$ and $c_j$ ($i\neq j$). Cre bound to a pair of strings $c_i\cup 
c_j$ can be modeled as a 2-string tangle $P$ of type $\frac0{1}$. Earlier studies of 
Cre support the assumption that Cre recombination takes $P=\frac0{1}$ into 
$R=\frac1{0}$, where for both tangles, the Cre binding sites are in anti-parallel 
orientation (Figure \actionofCre) \cite{PJH, GBJ, GGD, KBS}.  
Note that from a 3-dimensional point of view, the two sites can be regarded as parallel or anti-parallel as we 
vary the 
projection \cite{SECS, VCS, VDL}.
With our choice of 
framing any 
Cre-DNA complex formed by bringing together two {\em loxP} sites ({\em 
e.g.} strings $c_i$ and $c_j$) results in $P=(0)$ with anti-parallel sites when the 
$loxP$ sites are directly repeated. Furthermore, it is possible that Cre recombination 
traps crossings outside of the Mu and Cre protein-DNA complexes.  For mathematical 
convenience we will enlarge the tangle representing Cre to include these crossings 
which are not bound by either Mu or Cre but are trapped by Cre recombination.  That 
is, the action of Cre recombination on $c_i\cup c_j$ will be modeled by taking 
$P=\frac0{1}$ into $R=\frac1{d}$ for some integer $d$.  Hence the system of tangle 
equations shown in Figure \normaleq~ can be used to model these experiments where 
the rational tangles $\frac1{d_{i}}$, $\frac1{v_{i}}$, $\frac1{d_t}$, $\frac1{v_t}$ 
represent non-trivial topology trapped inside $R$ by Cre recombination, but not 
bound by Mu. 
The tangle $\T$, representing the transpososome ({\em i.e.} the Mu-DNA complex), is assumed to remain constant throughout
the recombination event \cite{SECS, PJH}.
 Recall that the first six experiments where the three Mu binding 
sites, attL, attR, and the enhancer, are all placed on the same circular DNA 
molecule will be referred to as the {\em in cis} experiments.  The remaining two 
experiments will be referred to as the {\em in trans} experiments since the enhancer 
sequence is provided {\em in trans} on a linear DNA molecule separate from the 
circular DNA molecule containing the attL and attR sites. The tangle equation (1) in 
Figure~\normaleq~corresponds to the unknotted substrate equation from the first six 
experiments.  Equations (2)-(4) correspond to three product equations modeling the 
three {\em in cis} deletion experiments, while equations (5)-(7) correspond to the 
three product equations modeling the three {\em in cis} inversion experiments. 
Equation (8) corresponds to the unknotted substrate equation for the two {\em in 
trans} experiments, while equations (9) and (10) correspond to the product equations 
modeling {\em in trans} deletion and {\em in trans} inversion, respectively.  In addition to 
modeling experimental results in \cite{PJH}, these equations also model results in 
\cite{PJH2,YJH}.

\myfig{7.50}{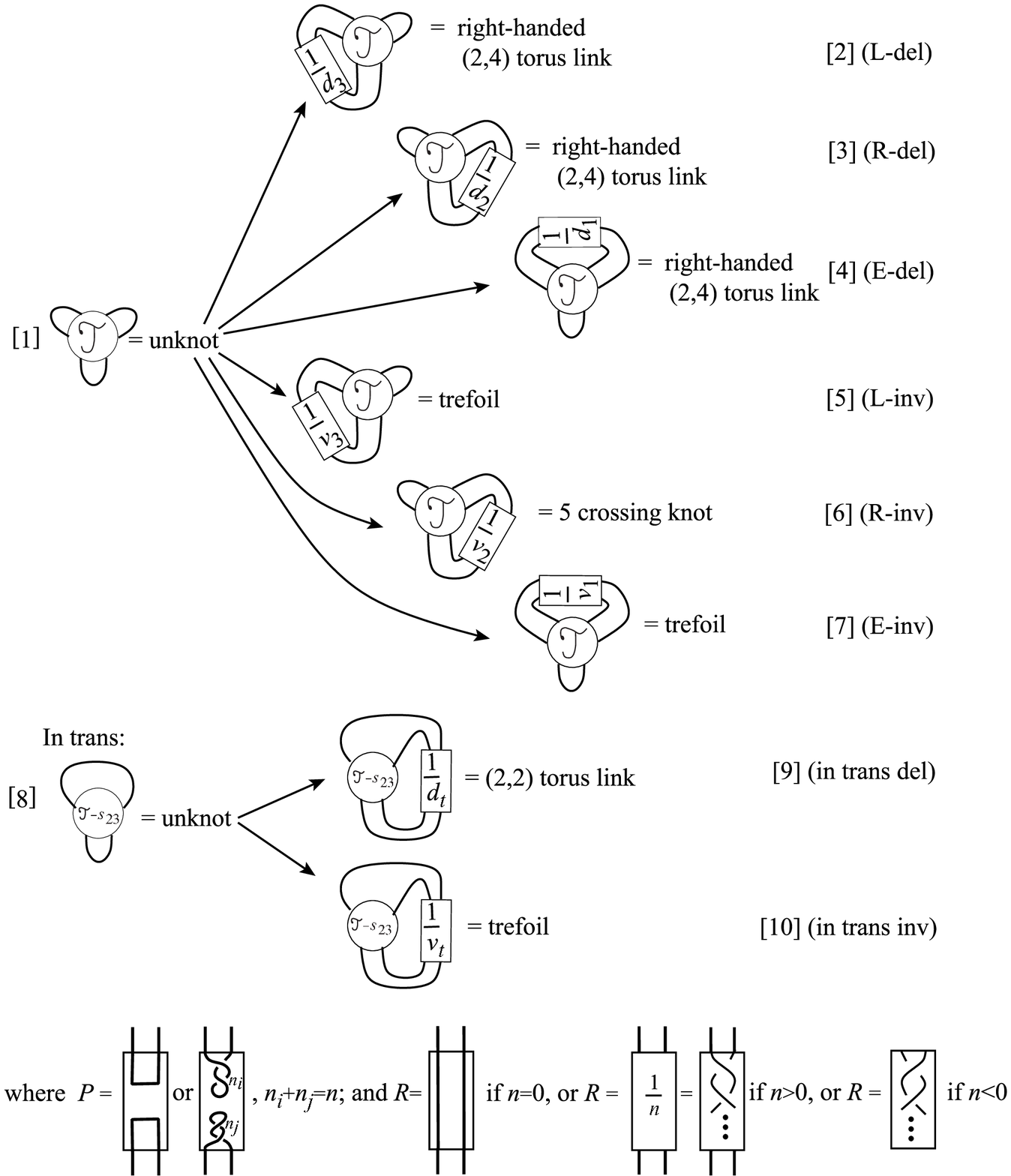}{\normaleq}

Figure~\projcircle~ partially defines a framing for the tangle $\T$. One can go further by specifying values for $d_{1}, d_{2}, d_{3}$ for the three {\em in cis} deletion experiments. 
%%A tangle \T as in Figure~\projcircle~ is in {\em normal form} if it 
%%satisfies equations (1)-(4) of Figure~\normaleq~ 
%%with $d_1 = d_2 = d_3 = 0$
 We define the {\em normal form equations} to be the system of equations 
(1) - (4) in Figure \normaleq~ (corresponding to the {\em in cis} deletion 
experiments) with the additional requirement that $d_{1} = d_{2} = d_{3} = 
0$.
 We focus first on these four equations as only the  
 {\em in cis} deletion  experiments are needed for our main results;
 but the other experiments were important experimental controls and were 
used by \cite{PJH} to determine the tangle model for the Mu transpososome.

\begin{defn}\label{defn:synapsetangle} A 3-string tangle $\T$ is called a 
{\em solution tangle} iff it satisfies the 
four normal form equations.   
\end{defn}

Note that in the normal form equations, 
the action of Cre  results in replacing $P = c_i \cup c_j = \frac0{1}$ 
with $R = \frac1{0}$ for the three {\em in cis} deletion experiments. 
If we wish to instead impose a framing where $P = c_i \cup c_j = 0/1$ is replaced by $R = 
\frac1 {d_{k}}$ for given $d_{k}$, $k =1, 2, 3$, we can easily convert between solutions.  
Suppose $\T$ is a solution to this non-normal form system of equations.  We can
move $n_i$ twists from $R$ into $\T$ at $c_i$, 
for each $i$ where $n_i + n_j = d_{k}$ (Figure \unique).  Hence $\T$ with 
$n_i$ 
crossings added inside $\T$ at $c_i$ for each $i$ is a solution tangle (for the normal form equations). Note the $n_i$ are uniquely determined. Similarly, if $\T$ 
is a solution tangle (for the normal form equations), then for given $d_k$ we can add $-n_i$ twists to $\T$ at 
$c_i$ to obtain a solution to the non-normal form equations.

\begin{remark}
For biological reasons \cite{PJH} chose  
$d_{1} = -1$, $d_{2} = 0$, $d_{3} = -1$.  Hence $n_{1} = 0$, $n_{2} = -1$, 
$n_{3} = 0$. 
This corresponds to adding a right-hand twist at $c_2$ to a solution in the 
PJH convention (such as Figure \PJHsolutionUnmarked) to obtain 
a solution tangle  (such as Figure \tangcircA). Conversely, by adding a single left-hand 
twist at $c_2$ to a solution tangle, we get the corresponding solution in 
the PJH 
convention (e.g. from the solution tangle of Figure \tangcircA~ to the PJH 
solution of
Figure \PJHsolutionUnmarked).

\end{remark}

\myfig{1.5}{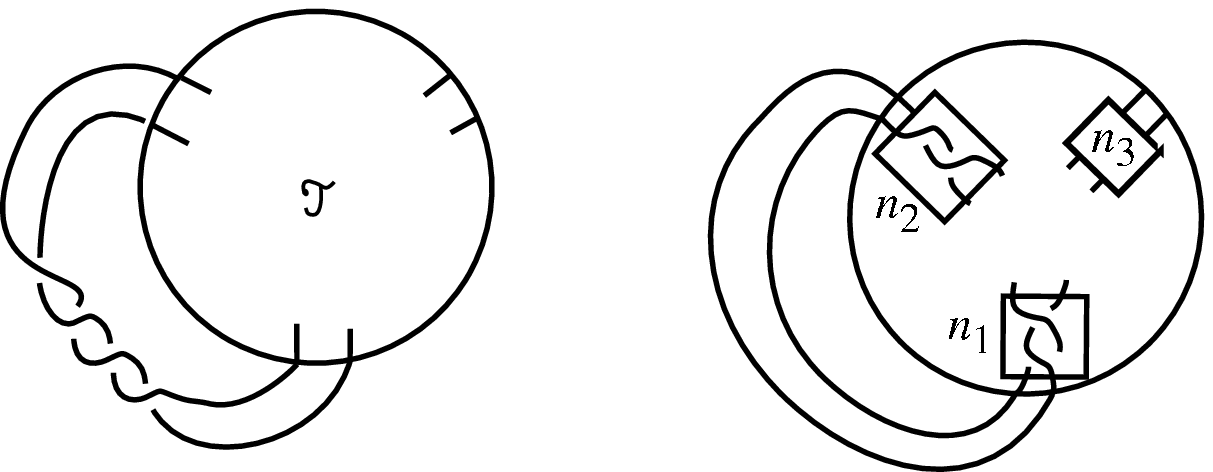}{\unique}

\begin{remark}
We will  show that knowing the number of crossings trapped by Cre outside of the 
Mu transpososome in the three {\em in cis} deletion experiments is sufficient for 
determining the number of crossings trapped outside the Mu transpososome in all 
experiments.  In other words choosing $d_{1}, d_{2}, d_{3}$ determines $v_{1}, 
v_{2}, v_{3}, v_t, $ and $d_t$.
\end{remark}

Next we apply traditional 2-string tangle calculus \cite{ES1} to the equations 
arising from the {\em in cis} deletion experiments.

\subsection{The three {\em in cis} deletion tangle equations}

With abuse of the $c_i$ notation of Figure~\projcircle, let the tangle circle be the union of arcs $c_1\cup 
x_{12} \cup c_2\cup x_{23}\cup c_3\cup x_{31}$ (Figure~\tangcircA). We think of 
these arcs as providing a framing for the tangle $\T$.
\myfig{1.5}{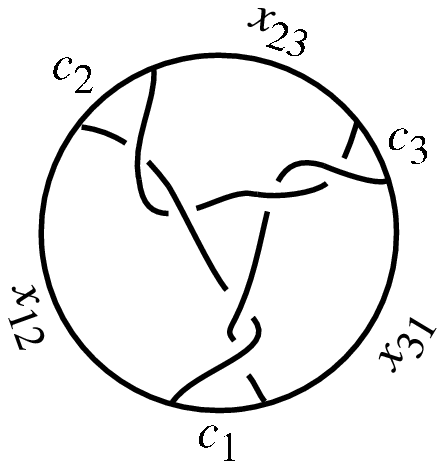}{\tangcircA}
Let $O_{i}=\T\cup c_i$ be the 2-string tangle obtained by capping off 
$s_{ji} \cup s_{ik}$ with $c_i$. Then $s_{kj}$ is one of the strings of the 2-
string tangle $O_{i}$, let $\hat s_i$ denote the other (Figure~\tangcircBm) 

\myfig{1.75}{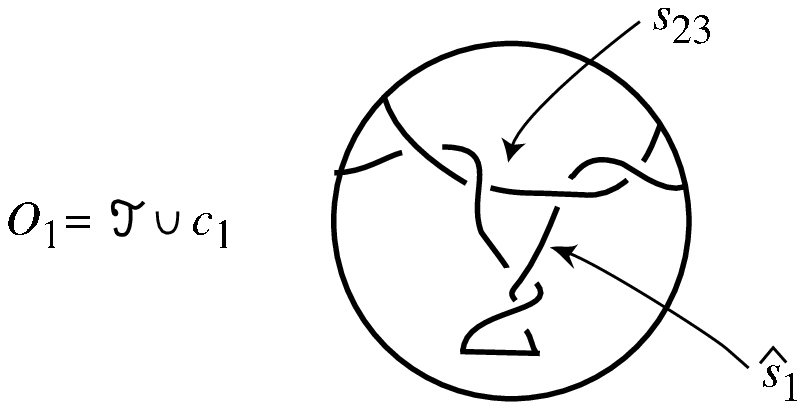}{\tangcircBm}

By capping $\T$ along each $c_i$, we approach the problem of finding all possible 
3-string tangles for the Mu transpososome by first solving three  
systems of two 2-string tangle equations (one for each $c_i$). 
Figure \tangleaddclose~illustrates the definitions of 2-string tangle 
addition (+) and numerator closure (N). Figure~\tangleEqns~shows a system 
of two 2-string tangle equations arising from Cre recombination on the Mu 
transpososome. $\T\cup c_2$ is represented by the blue and green 2-string tangle 
while the smaller pink 2-string tangle represents Cre bound to {\em loxP} sites 
before (left) and after (right) recombination.

\doublefig{.85}{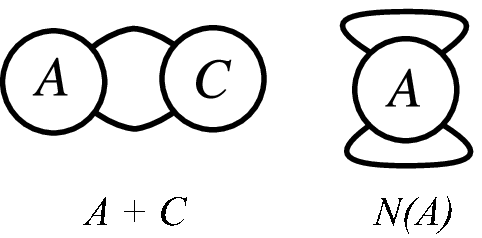}{\tangleaddclose}{.85}{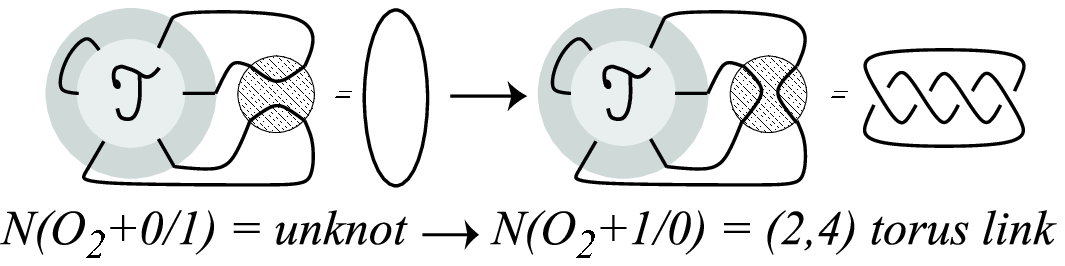}{\tangleEqns}

\begin{lemma}\label{lem:string} ({in cis} deletion)
Let $\T$ be a solution tangle. That is, consider the three  {in cis} deletion 
experiments which convert unknotted substrates to right-hand $(2,4)$ torus 
links. Let $O_{i}=\T\cup c_i$ be the 2-string tangle obtained by capping $\T$ 
along $c_i$. Then $O_i$ 
is the $-1/4$ tangle, $i=1,2,3$. 
\end{lemma}

\begin{proof}
All the {\em in cis} deletion events, modeled by replacing $P= c_j \cup 
c_k = {0 \over 1}$ by $R={1 \over 0}$ in normal form (see Figure~\toruslink), lead to identical 
systems of two 2-string tangle equations: 

\begin{equation*}  %\label{eq:deletionsystem} 
\begin{array}{rcl}
N(O_{i}+{0 \over 1})= &\text{ unknot}\\ 
\noalign{\vskip6pt}
N(O_{i}+{1 \over 0})=(2,4) &\text{ torus link}
\end{array}\qquad 
\end{equation*}

 By \cite{HS}, $O_i$ is a rational tangle. By tangle calculus, the 
system admits a unique solution $O_{i}=-1/4$, $i=1,2,3$ \cite{ES1} (see also 
\cite{D, VCS}).

\end{proof}

\myfig{1.2}{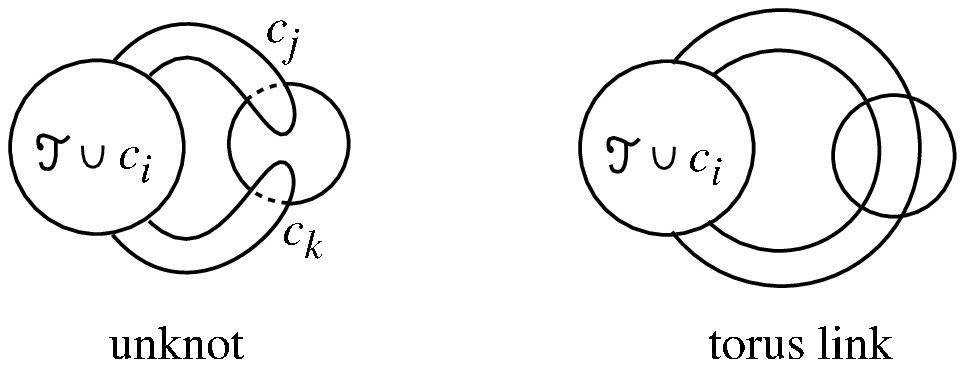}{\toruslink}

\subsection {The three {\em in cis} inversion tangle equations}

Each Cre inversion experiment is also modeled as a system of tangle equations:
\begin{equation*}  %\label{eq:deletionsystem} 
\begin{array}{rcl}
N(O_i+{0 \over 1})= &\text {unknot}\\ 
\noalign{\vskip6pt}
N(O_i+{1 \over v_{i}})= &\text{ inversion product}
\end{array}\qquad 
\end{equation*}

If we assume that the DNA-protein complex is constant for each pair of 
deletion/inversion experiments, then
by Lemma~\ref{lem:string}, $O_{1}= O_{2}= O_{3}=-1/4$. 
If the inversion 
products are known, $v_{i}$, $i = 1, 2, 3$, can be determined by tangle calculus. From \cite{PJH} we assume 
that both 
L-inversion and E-inversion produce the (+)trefoil. Also R-inversion produces a 5-crossing knot (see 
Table~\tableA\ ), which must be a 5-torus knot since $O_{2}=-1/4$.  We assume this is the (+)5-torus knot, 
for if it were the (-)5-torus knot then $v_{2} = 9$ which is biologically unlikely. Hence the only 
biologically reasonable solutions to the inversion equations, with $O_i=-1/4$, $i = 1, 2, 3$ are:

\begin{equation}\label{eq:InvSols}
\begin{split}
&\text{ L-inversion, $P={0 \over 1}$ and $R={1 \over 1}$, product 
=(+)trefoil;}\\ 
\noalign{\vskip-1pt}
&\text{ R-inversion, $P={0 \over 1}$ and R=${-1 \over 1}$, product = (+)5-torus 
knot;}\\
\noalign{\vskip-1pt}
&\text{ E-inversion, $P={0 \over 1}$ and $R={1 \over 1}$, product = 
(+)trefoil}\\
\end{split}
\end{equation}

Other possible solutions exist if $O_i$ is allowed to change with respect 
to the deletion {\em versus} inversion experiments. However, it is 
believed that the 
orientation of the Cre binding sites (inverted versus direct repeats, 
Figure~\recombdirknot2seq1a) does not affect the transpososome 
configuration.

\subsection{Linking number considerations}

%% Detecting changes in linking number introduced by an enzymatic action is 
%% important for both experimental and theoretical studies of site-specific 
%% recombination. 

In Lemma \ref{lem:synapsetangle} we compute  linking numbers related to the deletion experiments.

Assume $\T $ is a solution tangle. First, let $\hat x_i$ be the arc on the tangle circle given by $x_{ji} 
\cup c_i \cup x_{ik}$. We compute the linking number between $x_{jk} \cup s_{jk}$ and $ \hat s_i \cup \hat x_i$,  using 
the orientation induced from the tangle circle into each component and using the sign convention in 
Figure~\convention . The linking number quantifies the pairwise interlacing of arcs after capping off.   

Second, we compute the linking number of $x_{ij} \cup s_{ij}$ and $x_{ki} \cup s_{ki}$, with sign and orientation 
conventions as before. This linking number calculation quantifies the interlacing between any two arcs in $\T$.

\myfig{1.0}{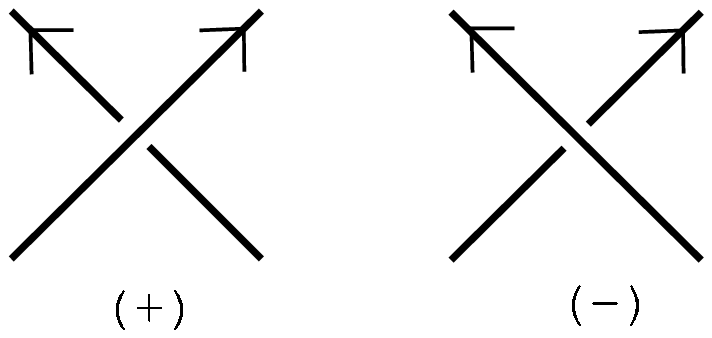}{\convention}

\begin{lemma}\label{lem:synapsetangle}
Assume $\T$ is a solution tangle. Then
%%If  {in cis} deletion on unknotted substrates produces $(2,4)$ torus %%links, then 
\begin{equation*}  %\label{eq:synapsetangle} 
\begin{array}{rcl}
\ell k(x_{jk} \cup s_{jk}, \hat s_i \cup \hat x_i) &=& -{2}\\ 
\noalign{\vskip6pt}
\ell k(x_{ij} \cup s_{ij}, x_{ki} \cup s_{ki}) &=& {-1}
\end{array}\qquad 
\{ i,j,k\} = \{1,2,3\}
\end{equation*}
\end{lemma}

\begin{proof} 
The first equation follows from Lemma~\ref{lem:string}.

To get the second equation we note that $\hat s_i\cup \hat x_i$ can be obtained 
via a banding connecting $x_{ij} \cup s_{ij}$ with $x_{ki}\cup s_{ki}$ along 
$c_i$ as indicated in Figure~\bandingm. 
\myfig{1.75}{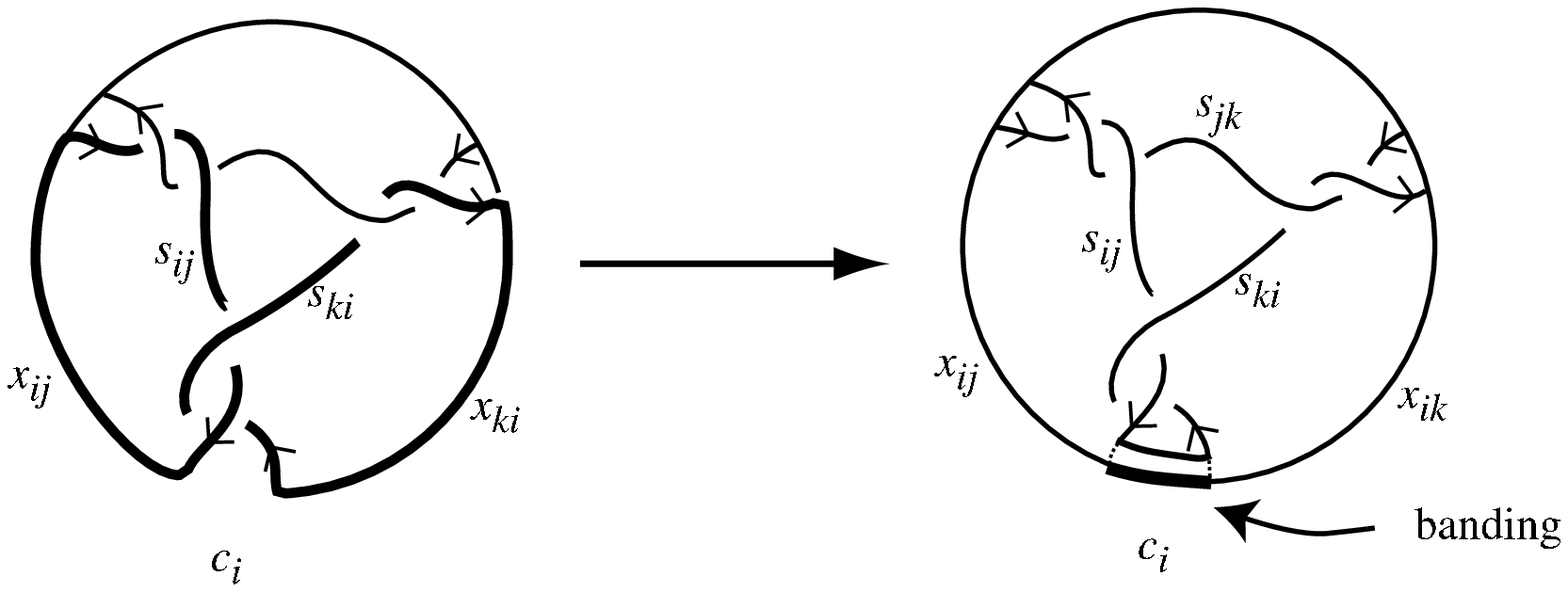}{\bandingm} 
Then for  $\{i,j,k\} = \{1,2,3\}$,
\begin{equation*}
-{2} = \ell k (x_{jk} \cup s_{jk},\hat s_i\cup \hat x_i) 
= \ell k (x_{jk} \cup s_{jk} , x_{ij}\cup s_{ij}) 
+ \ell k (x_{jk} \cup s_{jk}, x_{ki} \cup s_{ki})\ .
\end{equation*}
Solving three equations ($\{i,j,k\} = \{1,2,3\}$), we obtain  
$\ell k (x_{ij}\cup s_{ij}, x_{jk}\cup s_{jk})={-1}$. 
\end{proof}

In the next section
we will see how the {\em in cis} deletion results extend to the analysis of the 
remaining experiments ({\em in trans} deletion and inversion).

\subsection{{\em In trans} experiments}

In the {\em in trans\/} portion of \cite{PJH} (described in section 1) the 
enhancer sequence is not incorporated into the circular DNA substrate; 
it remains separate in solution on its own linear molecule when the 
transpososome is formed (Figure~\intrans). 
Assuming the transpososome complex forms the same 3-string tangle in this 
context, the transpososome writes the substrate as a union of two 2-string 
tangles: $\T'_0$, a trivial tangle outside the transpososome; and $\T' = \T- 
s_{23}$, the tangle formed by the {\em attR\/} and {\em attL\/} strands (where $s_{23}$ is the enhancer strand in $\T$).

\myfig{1.75}{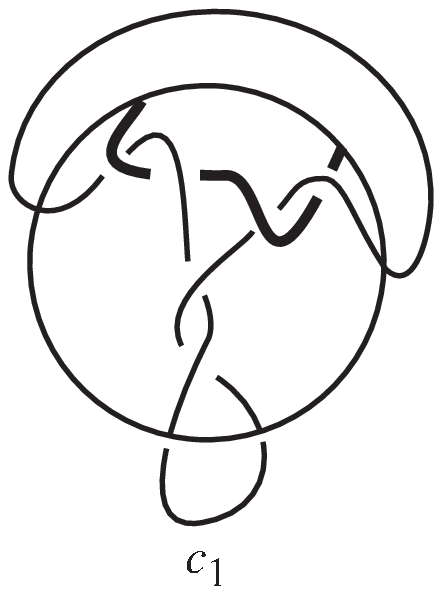}{\intrans}

\begin{lemma}\label{lem:enhancer} (in trans deletion)
Let $\T$ be a solution tangle.  Suppose $\T$ also satisfies equations (8) and (9) in Figure~\normaleq~ where $d_t$ is as shown in this figure.
%%Let $\T' =\T -s_{23}$ where $s_{23}$ is the enhancer strand in $\T$. 
Then $d_t = 0, 4$ and $\T - s_{23} = {1 \over -2}$.  
\end{lemma}

\begin{proof}
When Cre acts on the $loxP$ sites in the {\em in trans} experiment, we assume 
that it takes the ${0\over 1}$ tangle to ${1\over d_{t}}$ (
Figure~\unknotHopf ). The Cre deletion product {\em in trans\/} is a Hopf link 
({\em i.e.\/} a $(2,2)$ torus link). Hence we are solving the 2-string tangle equations, 
$N((\T - s_{23}) + {0\over 1}) =$ unknot, $N((\T - s_{23}) + {1\over d_{t}}) = 
(2, 2)$ torus link.  By \cite{BL} $\T - s_{23}$ is rational, and by tangle 
calculus \cite{ES1}, $\T - s_{23} = {1 \over {\pm 2-d_{t}}}$, which implies that 
$\ell k (x_{31} \cup s_{31}, x_{12}\cup s_{12}) = {\pm 2-d_{t} \over 2}$.  By 
Lemma~\ref{lem:synapsetangle}, $\ell k (x_{31} \cup s_{31}, x_{12}\cup s_{12}) 
=-1$.  Hence ${\pm 2-d_{t} \over 2} = -1$.  Thus $d_t = 0, 4$ and $\T - s_{23} = 
{1 \over -2}$.  
\end{proof}

\myfig{1.25}{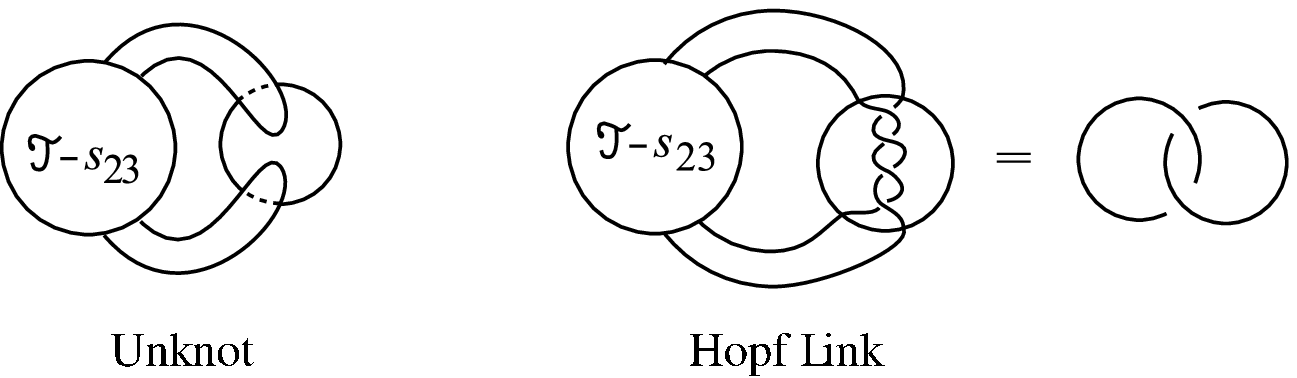}{\unknotHopf}

\vskip 5pt
\remark (\bf{{\em In trans} inversion})
\rm If $\T - s_{23} = {1 \over -2}$, then the {\em in trans} Cre inversion 
product is $N(\T - s_{23} + {1 \over v_t}) = N( {1 \over -2} + {1 \over v_t}) = 
(2, 2-v_t)$ torus knot.  A (2, 3) torus knot (i.e. the trefoil) {\em in trans} 
inversion product implies $v_t = -1$.

\subsection{Summary}

The next proposition summarizes the results in this section. 

\begin{prop}\label{prop:synapsetangle} 
Let $\T$ be a solution tangle which also satisfies the {in 
trans} deletion experiments, and let $s_{23}$ correspond to the 
enhancer strand in the \cite{PJH} experiments.  Then
%
%%\myfig{1.75}{synapsetangle}{\synapsetangle} 
%
%%Then in the Cre deletion experiments 
\begin{equation}\label{eq:enhancerPJH} 
O_{i}=\T \cup c_i = -\frac1{4} 
, ~\T - s_{23} = -\frac1{2}  
\end{equation}
The in trans deletion reaction of Cre is modeled by replacing $c_1 \cup (c_2 \cup x_{23} \cup c_3) = \frac0{1}$ by the $\frac1{0}$ or $\frac1{4}$ tangle (i.e. $d_{t}=0, 4$, see 
Figure~\deletionreaction).
\doublefig{1.5}{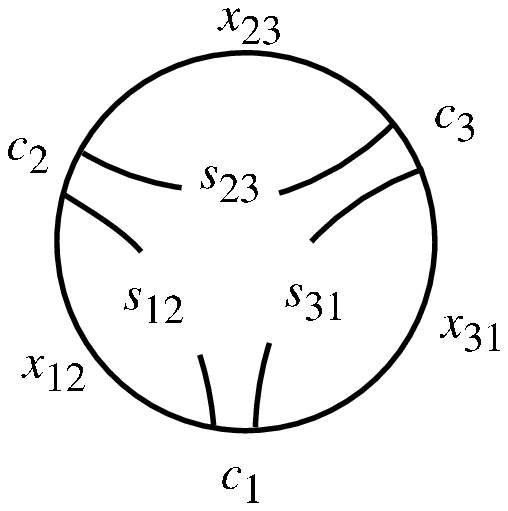}{\synapsetangle}{1.2}{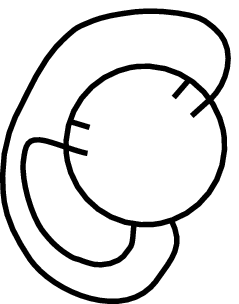}{\deletionreaction} 
%%\myfig{1.0}{deletionreaction}{\deletionreaction} 
\end{prop}

\begin{remark}
In declaring $O_i=-1/4$ we are using the framing coordinate of the tangle
circle with $c_j \cup c_k = 0/1 , x_{jk} \cup \hat x_i = 1/0$.
\end{remark}

Proposition \ref{prop:synapsetangle} can be generalized to products of deletion 
that are $(2,L_{i})$ torus links, and a $(2,L_{t})$ torus link for the {\em in 
trans} experiment:

\begin{prop}\label{prop:synapsetanglegeneral} 
Let $\T$ satisfy equations (1)-(4), (8)-(9) of Figure~\normaleq~except that the in cis
 deletion experiments produce $(2, L_i)$ torus links,
$L_i \not= \pm 2$,  for $i=1,2,3$ and the in trans deletion experiment results in a $(2, L_t)$ torus link. Assume $d_1 = d_2 = d_3 = 0$ and that $s_{23}$ corresponds to the enhancer strand. Then 
\begin{equation}\label{eq:enhancer} 
O_{i}=\T \cup c_i = -\frac1{L_{i}} 
, ~\T - s_{23} =\frac1{{L_2 + L_3 - L_1 \over 2}}.  \end{equation}
If $|L_t| = 2$, then $d_t = {-L_j + L_i + L_k \pm 4 \over 2}$.
If $|L_t| \not= 2$, then $d_t = {-L_j + L_i + L_k - 2L_t\over 2}$.
\end{prop}

Note similar results hold if $L_i = \pm 2$ for $i = 1, 2$, and/or $3$, but as 
this breaks into cases, we leave it to the reader.

\begin{remark}
Sections 3 and 4 rely on the fact the tangles in Figure \solutiontangle~ 
are all of the form $\frac1{m_i}$.  This would still be the case no matter the 
choice of framing (discussed at the beginning of this section)  or the type of the $(2, L_i)$ torus link products.   
Hence the 
results of sections 3 and 4 apply more generally.
\end{remark}

Note that a 3-string tangle $\T$ 
is a {solution tangle} iff it satisfies the condition $O_{i}=\T \cup c_i = 
-\frac1{4}$ for $i= 1, 2, 3$ (see 
Definition \ref{defn:synapsetangle} 
and Lemma \ref{lem:string}). 
This corresponds to the first three equations in Figure~\solutiontangle.
The main result in section 3 (Corollary \ref{corPJH}) only depends on the 
three {\em in cis} deletion experiments. 

\myfig{3.75}{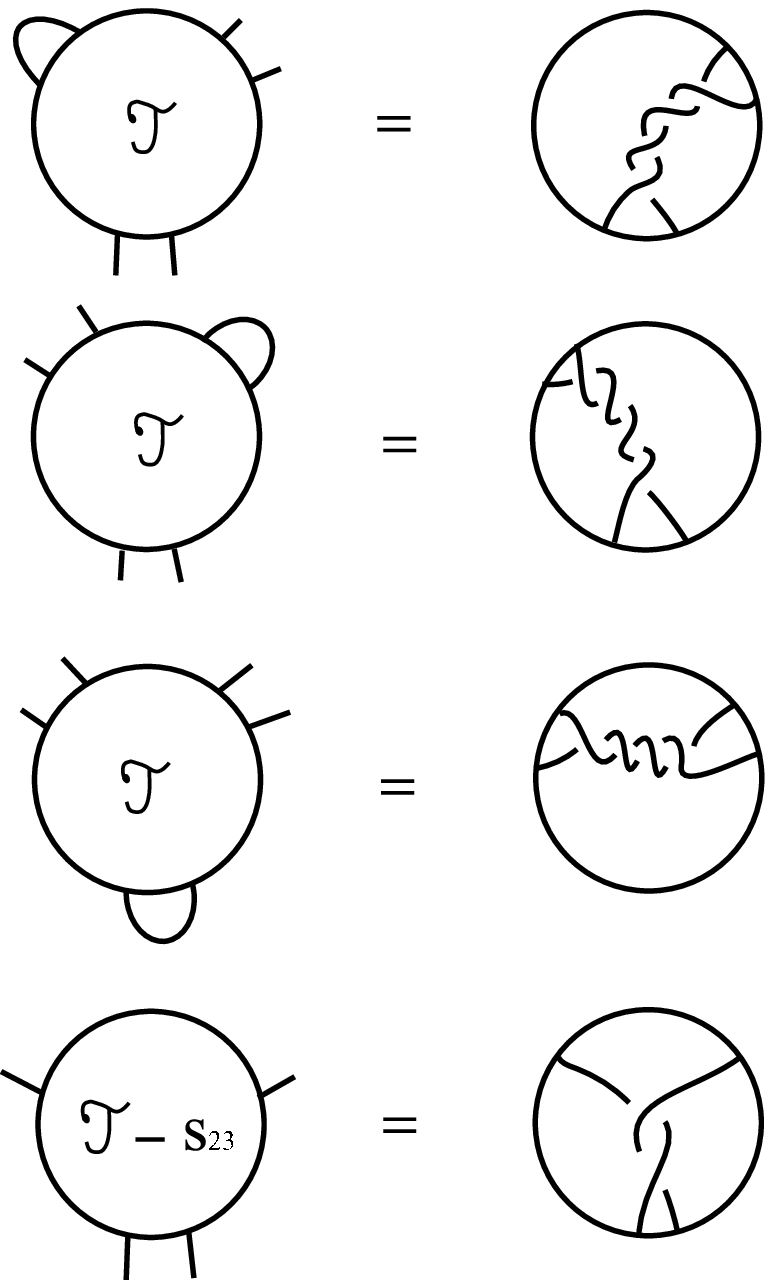}{\solutiontangle}

If in addition a solution tangle satisfies the {\em in trans} deletion equation, then we call it an {\em in 
trans 
solution tangle}:
 
\begin{defn}\label{defn:synapsetangle2}
A 3-string tangle $\T$ is called an 
{\em in trans solution tangle} 
iff it satisfies Equation \ref{eq:enhancerPJH} of Proposition ~\ref{prop:synapsetangle}.
This is equivalent to a tangle satisfying all the equations of Figure
~\solutiontangle.
%%it satisfies the three {\em in cis} deletion experiments.  
%%That is it satisfies the 
%%first three equations in Figure~\solutiontangle.
\end{defn}

Note: {\bf Proposition~\ref{prop:synapsetangle} reduces the study of 
transpososome tangles to that of {\em in trans} solution tangles.} In 
section~3 
we classify solution tangles in terms of solution graphs and use this to 
show that the only solution tangle (and hence the only {\em in trans} 
solution tangle) that is rational is the PJH tangle.

\begin{remark}
The proof of Propositions~\ref{prop:synapsetangle} and 
\ref{prop:synapsetanglegeneral} use the deletion products only. As the remaining 
analysis  depends only on Proposition~\ref{prop:synapsetangle}, we see that the 
inversion products are unnecessary for our analysis. However, the inversion 
products were used to determine the framing in \cite{PJH}.

\end{remark}

%% file: section3.tex
\section{Solution graphs}  

In sections 3.1 and 3.2 we relate solution tangles to wagon wheel graphs and 
tetrahedral graphs.  We use wagon wheel graphs to find an infinite number of solution 
tangles in section 3.1. 
Recall by Definition \ref{defn:synapsetangle} that a solution tangle is a 3-string 
tangle
satisfying the 
three {\em in cis} deletion equations ({\it i.e.}, $O_{i}=\T \cup c_i = -\frac1{4}$ for $i = 1, 2, 3$ 
in Equation 
\ref{eq:enhancerPJH} of Proposition 
%%the four equations in
%%(\ref{eq:enhancerPJH}) of Proposition
\ref{prop:synapsetangle}; first three equations of Figure \solutiontangle), 
while by Definition 
\ref{defn:synapsetangle2}, an {\em in trans}
solution tangle is a solution tangle which also satisfies the {\em in trans} deletion 
equation $\T - s_{23} = -\frac1{2}$.
 In section 3.3 we show that the only rational solution tangle is 
the PJH tangle.  
 In section 3.5, we extend this result to tangles which are split or have parallel 
strands using results in section 3.4 regarding the exterior of a solution graph. 
 In section 3.6, we show that if a solution tangle is also an {\em in trans} solution 
tangle, then its exterior is a handlebody.  We use this in section 3.7 to investigate 
the planarity of tetrahedral graphs.

\subsection{Solution tangles are carried by solution graphs}
\begin{defn}
A 3-string tangle is {\em standard} iff it can be isotoped rel endpoints 
to Figure~\standard.
A standard tangle with $n_i = {-2}$ is called the {\it PJH tangle}.

\myfig{1.5}{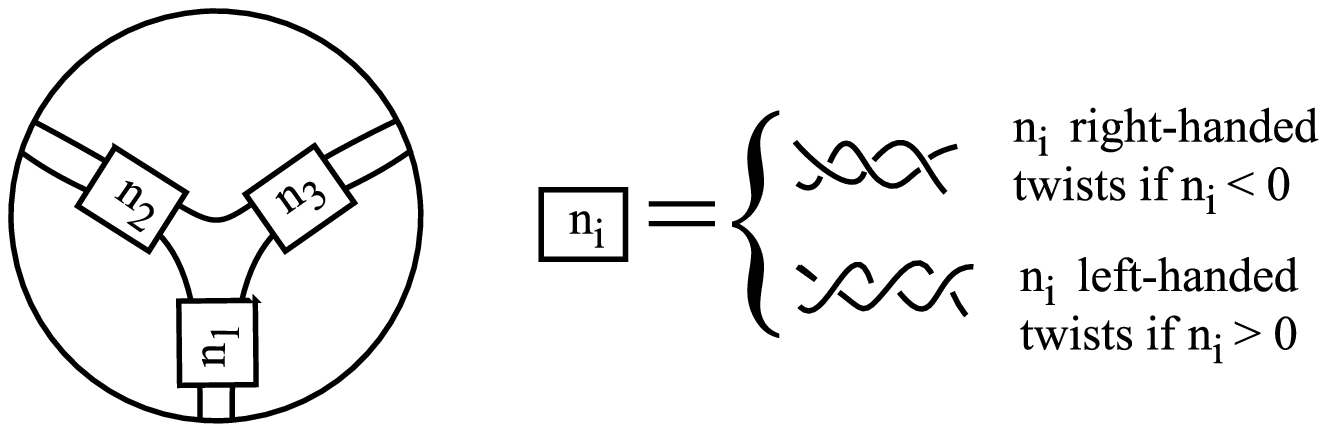}{\standard}		%% 3.1
\end{defn}

%%Recall by defintion ?? that a solution tangle is a 3-string tangle 
%%satisfying the four equations in
%%(\ref{eq:enhancerPJH}) of Proposition 
%%\ref{prop:synapsetangle} (Figure \solutiontangle), while an {\em in cis} 
%%solution tangle satisfies the first three 
%%equations of this proposition.

\begin{lemma}\label{lem:standard} 
If a solution tangle is standard, then 
it is the PJH tangle. 
\end{lemma}

\begin{proof} 
Write and solve the three equations, $n_i + n_j = -4$, describing the integral tangles 
resulting from capping along individual $c_k$.
\end{proof}

\begin{lemma}\label{lem:standard2} 
If an {in trans} solution tangle is standard, then 
it  is the PJH tangle. 
\end{lemma}

\begin{proof} 
By Lemma \ref{lem:standard}, since an {\em in trans} solution tangle is also a solution tangle.
\end{proof}

\begin{remark}
The assumption in \cite{PJH} of the branched supercoiled structure
             (see the end of section 1) is equivalent to assuming the {\em in trans} solution
             tangle is standard. In that paper, the authors deduced the 
             transpososome configuration with a method similar to that of 
             Lemma \ref{lem:standard}. The PJH solution thus 
obtained is equivalent to our PJH tangle when normal framing is imposed.
\end{remark}

\begin{defn}
The {\em abstract wagon wheel} is the graph of Figure~\wagonwheel. 
The vertices are labelled $v_1,v_2,v_3$; the edges  labelled 
$e_1,e_2,e_3,b_{12},b_{23},b_{31}$.
A {\em wagon wheel graph}, $G$, is a proper embedding of the abstract 
wagon wheel into $B^3$ with the 
endpoints of the $e_i$ in the 10, 2, and 6 o'clock positions
on the tangle circle ({\em e.g.} Figure~\wwgraph).
Two wagon wheel graphs are the {\em same} if there is an isotopy of $B^3$, 
which is fixed on $\partial B^3$, taking one graph to the other. 
\end{defn}

 \begin{defn} If a properly embedded graph lies in a properly embedded disk in the 3-ball, we call it {\em 
planar}.
 \end{defn} 

Figure~\nonplanar\ gives examples of planar and non-planar wagon wheel 
graphs.
Note that the non-planar graph in Figure~\nonplanar\ contains the knotted arc 
$ e_2 \cup b_{23} \cup e_3$.  The disk in which a planar wagonwheel lies may be taken to be the one given by the plane 
of the page, whose boundary is the tangle circle.  

\doublefig{1.3}{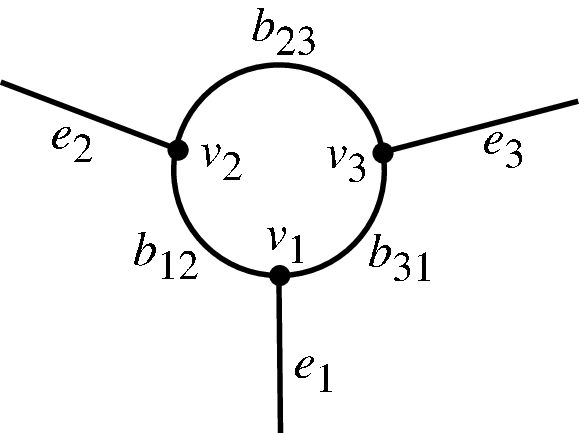}{\wagonwheel}{1.5}{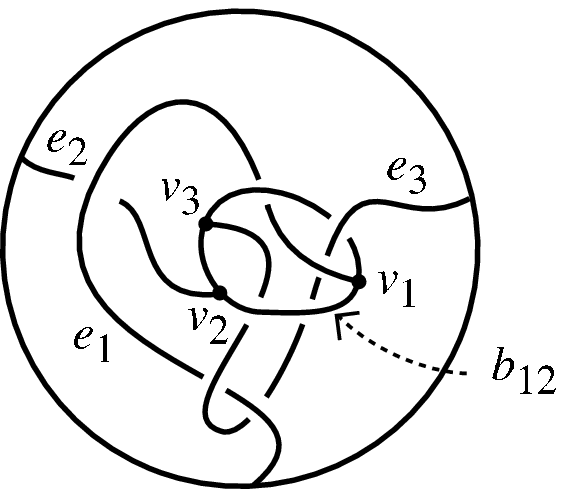}{\wwgraph} %% 3.2, 3.3

\myfig{1.5}{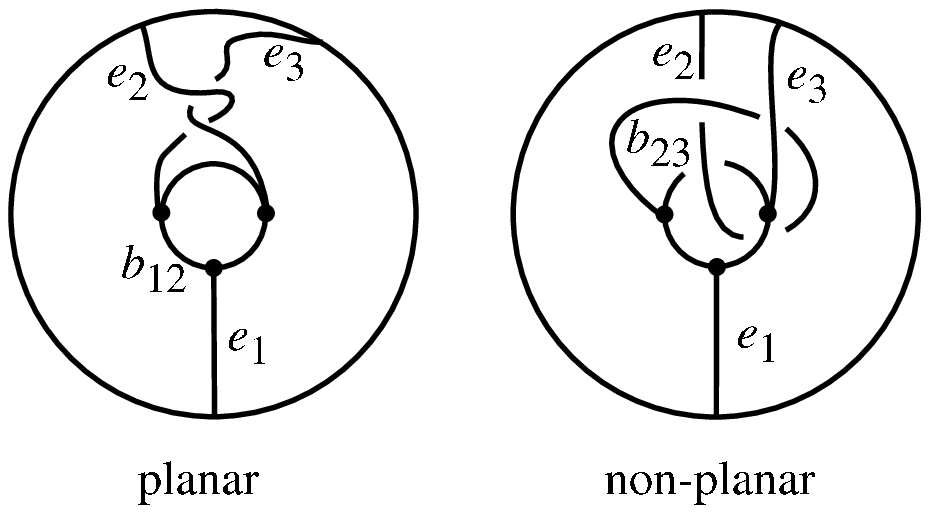}{\nonplanar}	%% 3.4

\begin{defn}
Let $G$ be a wagon wheel graph, and let 
 $N(G)$ denote a regular neighborhood of $G$ in the 3-ball. 
Let $\J_1,\J_2,\J_3$ be meridian disks of $N(G)$ corresponding to edges 
$e_1,e_2,e_3$. 
Let $J_i = \partial \J_i$. 
Let $\Gamma_{12},\Gamma_{23},\Gamma_{31}$ be meridian disks of $N(G)$ 
corresponding to edges $b_{12},b_{23},b_{31}$. 
Let $\gamma_{ij} = \partial\Gamma_{ij}$. 
See Figure~\NG.

\myfig{1.75}{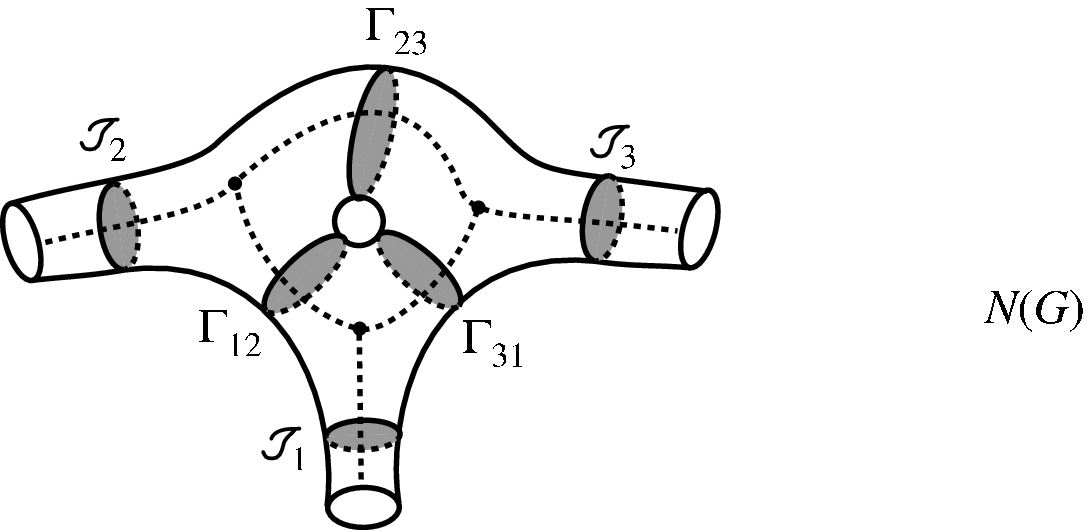}{\NG}	%% 3.5
\end{defn}

\begin{defn}
$X(G) =$ exterior of $G= B^3 - (N(G) \cup N(\partial B^3))$. 
Note that $\partial X(G)$ is a surface of genus~3. 
\end{defn}

\begin{defn}
$G$ is a {\em solution graph} iff it is a wagon wheel graph with the 
property that deleting any one of the edges $\{e_1,e_2,e_3 \}$ 
from $G$ gives a subgraph that is planar.
\end{defn}

\begin{defn}
$G$ is an {\em in trans solution graph} iff it is a wagon wheel graph with the 
property that deleting any one of the edges $\{e_1,e_2,e_3,b_{23}\}$ 
from $G$ gives a subgraph that is planar.
\end{defn}

The graphs G in Figure~\carriedbyG\ are {\em in trans} solution graphs.

\begin{defn}
Let $G$ be a wagon wheel graph. 
A 3-string tangle $\T$ is {\em carried by $G$}, written $\T(G)$, 
iff its 3-strings $s_{12},s_{23},s_{31}$ can be (simultaneously) 
isotoped to lie in $\partial N(G)$ such that 
\begin{itemize}
\item[1)] $s_{ij}$ intersects each of $J_i,J_j$ once
\item[2)] $s_{ij}$ intersects $\gamma_{ij}$ once and is disjoint from the 
remaining $\gamma_{st}$.
\end{itemize}
Figures~\carriedbyG\ gives examples of tangles carried by wagon wheel graphs.

\myfig{4.0}{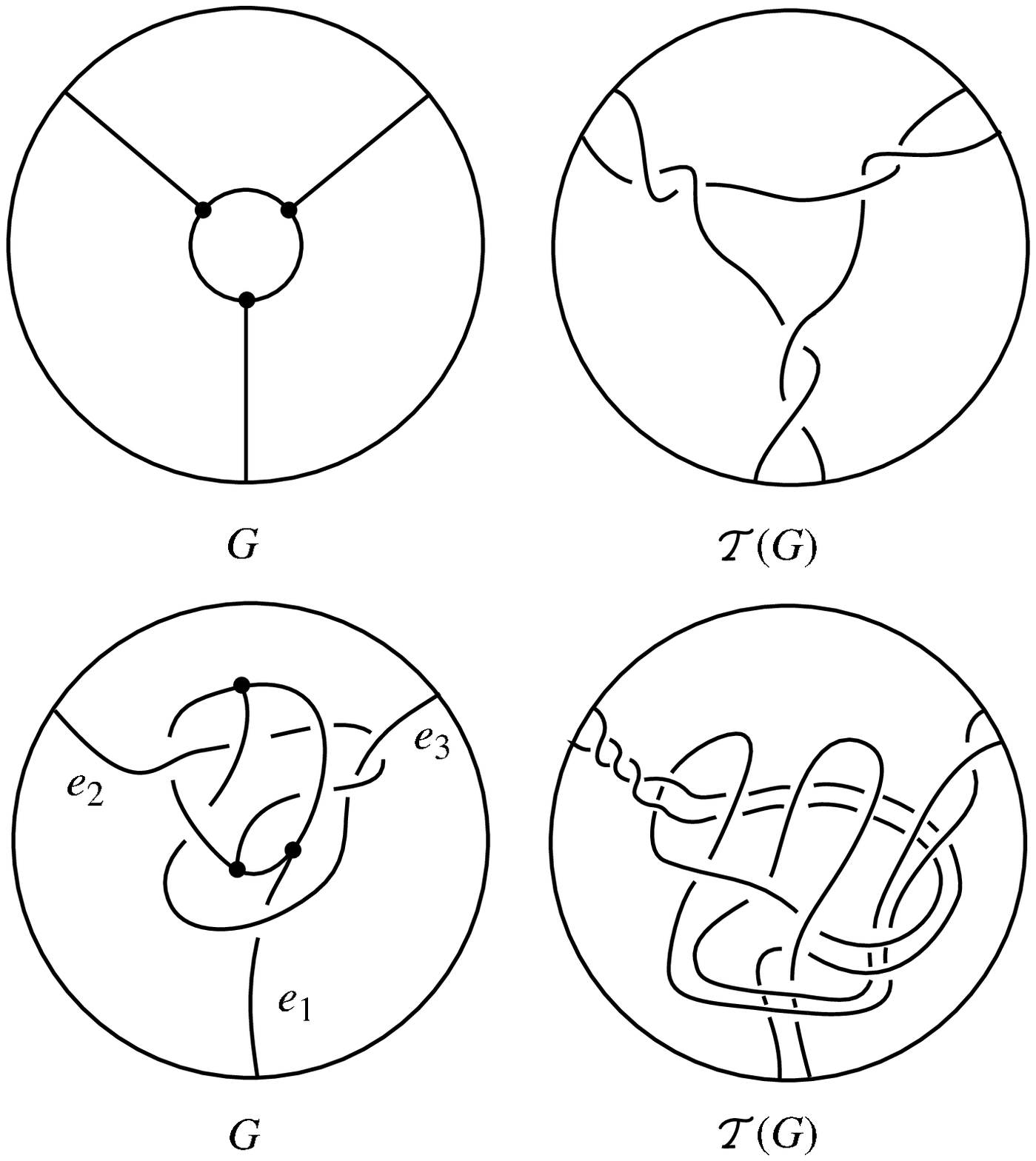}{\carriedbyG}		%% 3.6
\end{defn}

\begin{lemma}\label{lem:carriedbyG}
If $\T_1,\T_2$ are both carried by $G$, then they are isotopic in $B^3$ (rel endpoints)
up to twisting $N(G)$ along the $\J_i$ (i.e., they differ by twists at 
the ends).
\end{lemma}

\begin{proof}
Left to the reader.
\end{proof}

\begin{lemma}\label{lem:standardT}
If a 3-string tangle $\T$ is carried by a planar wagon wheel graph $G$ then $\T$ is  standard.
\end{lemma}

\begin{proof}
There is an isotopy of $B^3${, keeping the boundary of $B^3$ fixed,} taking $G$ to the 
planar wagon wheel graph of Figure \carriedbyG~(upper left). 
Apply Lemma~\ref{lem:carriedbyG} to see $\T$ as standard.
\end{proof}

The following theorem and corollary lets us work with 
%%({\em in cis}) 
%%solution 
graphs rather than 
%%({\em in cis}) 
%%solution 
tangles.
The 
%%({\em in cis}) 
solution graph economically encodes the conditions needed for a 
%%({\em in cis}) 
solution 
tangle.

\begin{thm}\label{thm:solutiontangle}
Let $\T$ be a solution tangle, then $\T$ is carried by a 
solution graph.
Conversely, a solution graph carries a unique solution tangle.
\end{thm}

\begin{proof}
Let $\T$ be a solution tangle with strands $s_{12},s_{23},s_{31}$. 
We show that $\T$ is carried by a solution graph $G$. 
Let $C = c_1 \cup s_{12} \cup c_2\cup s_{23}\cup c_3\cup s_{31}$, recalling that the $c_i$ are the capping arcs of Figure \synapsetangle.
For each $i=1,2,3$, let $P_i$ be a point in $c_i$. 
Isotope $C$ into the interior of $B^3$ by pushing each $c_i$ slightly to the 
interior. 
Under the isotopy $P_i$ 
traces out an arc $e_i$ from $P_i \in \partial B_3$ to $C$. 
\myfig{2.5}{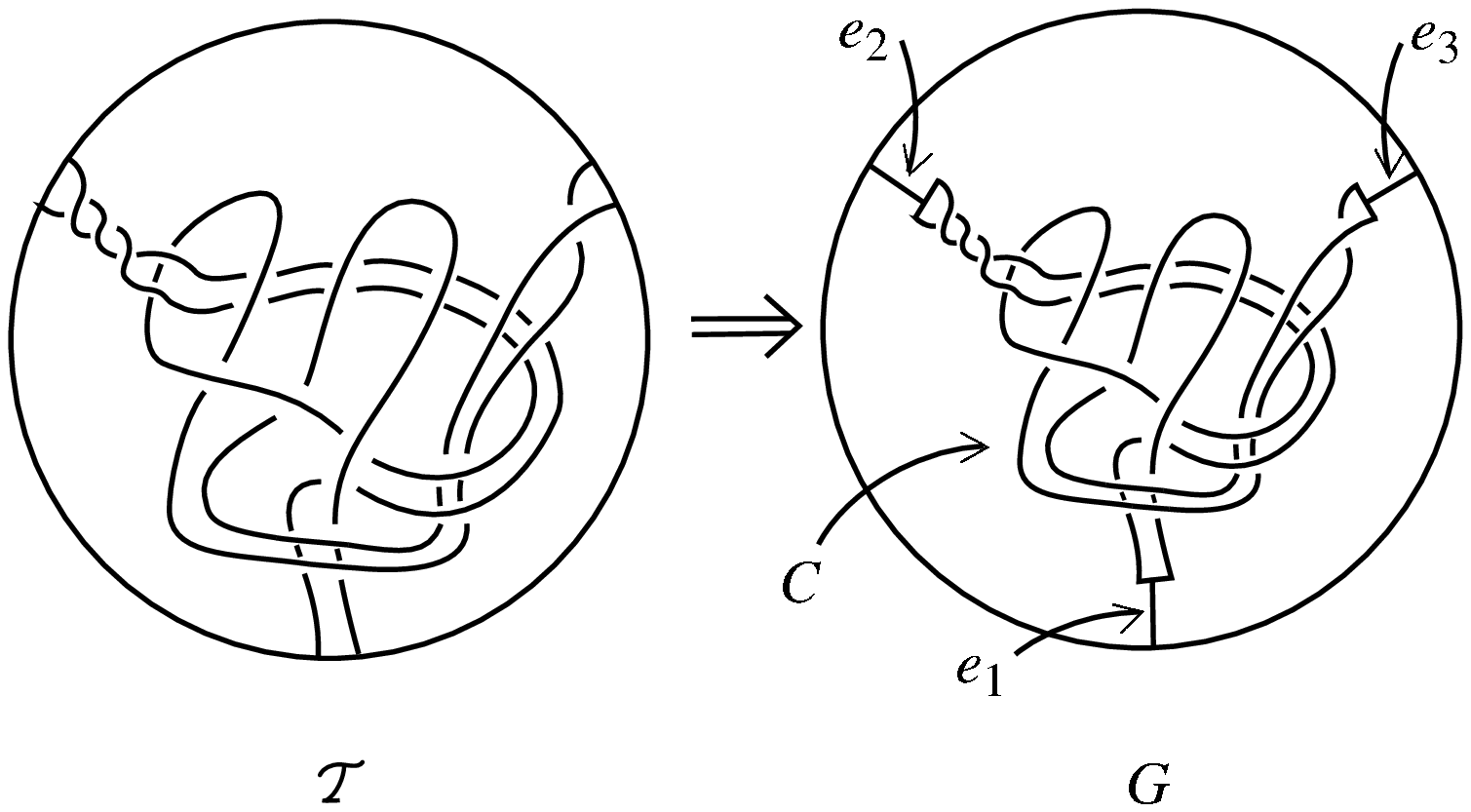}{\solutiontangleA}		%% 3.7
Define $G= C\cup e_1\cup e_2\cup e_3$.
See Figure~\solutiontangleA.  
$G$ is a wagon wheel graph that carries $\T$. 
Note that the strings $s_{ij}$ of $\T$ correspond to the edges 
$b_{ij}$ of $G$. 
The conditions for $\T$ to be a solution tangle now correspond to those
defining $G$ as a solution graph. 
For example, 
%%$G- e_1$ corresponds to the graph 
%%$b_{23}\cup (b_{12} \cup c_1 \cup b_{31}) \cup c_2 \cup c_3 \cup e_2 \cup e_3$ 
if $\T$ is a solution tangle, $s_{23}\cup (s_{12} \cup c_1 \cup s_{31}) \cup c_2 
\cup c_3$ (where $c_1$ has been pushed into the interior of $B^3$) lies in a 
properly embedded disk $D$  in $B^3$. 
Hence $G- e_1$ lies in $D$ and is planar. 

We now show that a solution graph carries a unique solution tangle. By Lemma~\ref{lem:carriedbyG} the tangles carried by $G$ are 
parameterized by integers $n_1,n_2,n_3$ as in Figure~\tangleG~(where the boxes are 
twist boxes). Denote these tangles by $\T(n_1,n_2,n_3)$. As $G$ is a 
solution graph, $\T(0,0,0)\cup c_i$ will be the 2-string
 tangle {${1 \over f_i}$}, 
where $f_i\in\zed$.
Then $\T(n_1,n_2,n_3)\cup c_i$ will be the 2-string tangle 
{${1 \over f_i+n_j+n_k}$} where $\{i,j,k\} = \{1,2,3\}$. 
\myfig{1.75}{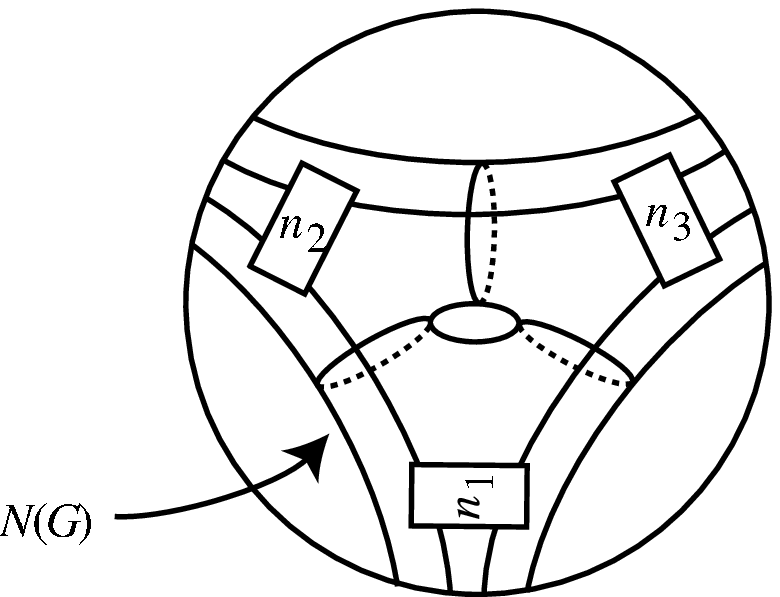}{\tangleG}		%% 3.8
Based on Lemma \ref{lem:string}, we choose $n_i,n_j,n_k$ to satisfy 
\begin{equation*}
\begin{split}
&n_1 + n_2 + f_3 = -4\\
&n_2 + n_3 + f_1 = -4\\
&n_3 + n_1 + f_2 = -4
\end{split}
\end{equation*}
Since  {$f_1 + f_2 + f_3$ is} even, this has a unique integer solution. 
Then $\T(n_1,n_2,n_3) \cup c_i = {1 \over -4}$ tangle for each $i$. Thus $\T$ is a solution tangle.

The uniqueness of $\T$ follows from the parameterization of the tangles 
by the $n_i$ and the fact that $\T \cup c_i = {1 \over -4}$ tangle.

\end{proof}

\begin{cor}\label{cor:solutiontangle}
Let $\T$ be an {in trans} solution tangle, then $\T$ is carried by an {in trans} solution graph.
Conversely, an {in trans} solution graph carries a unique {in trans} solution tangle.
\end{cor}

\begin{proof}

Let $\T$ be an {\em in trans} solution tangle. Since $\T$ is an {\em in trans} solution tangle, $\T$ is a solution tangle.  Hence $\T$ 
is carried by a solution graph, $G$. 
$G- b_{23}$ corresponds to the graph $s_{12}\cup c_1 \cup s_{31} \cup e_1$. As $\T$ 
is an {\em in trans} solution tangle, $s_{12}\cup c_1\cup s_{31}$ lies in a properly embedded disk 
$D$ in $B^3$. Hence $G- b_{23}$ lies in $D$ and is planar.  Thus $G$ is an {\em in trans} solution 
graph.

Suppose $G$ is an {\em in trans} solution graph.  Thus $G$ is a solution graph and carries a solution tangle.
Since $G$ is an {\em in trans} solution graph, $G- b_{23}$ is also planar. Then the 2-string tangle 
$\T (n_1,n_2,n_3) - s_{23} = {1 \over f}$ for some integer $f$. 
Furthermore, we have chosen $n_1, n_2,n_3$ so that $\ell k (x_{ij} \cup s_{ij}, 
\hat s_i\cup \hat x_i) =-2$.
The argument of Lemma~\ref{lem:synapsetangle} shows that $\ell k(x_{12} 
\cup s_{12}, x_{31}\cup s_{31}) =-1$.
Thus $\T(n_1,n_2,n_3) - s_{23} = {1 \over -2}$ tangle. Therefore $\T$ is an {\em in trans} solution tangle.  $\T$ is 
unique by Theorem~\ref{thm:solutiontangle}.

\end{proof}

Corollary~\ref{cor:solutiontangle} allows us to construct infinitely many 
{\em in trans}
solution tangles, via 
{\em in trans} 
solution graphs. 
Recall that a {\it Brunnian link}, $L$, is one such that once any 
component is removed the remaining components form the unlink.   
By piping a Brunnian link to a given 
{\em in trans} 
solution graph, we generate a new 
{\em in trans}
 solution graph as pictured in Figure~\Infinitely\ ~(where we 
have 
begun with the second solution graph of Figure~\carriedbyG). 
Similar results
hold for
solution tangles.

\myfig{2.25}{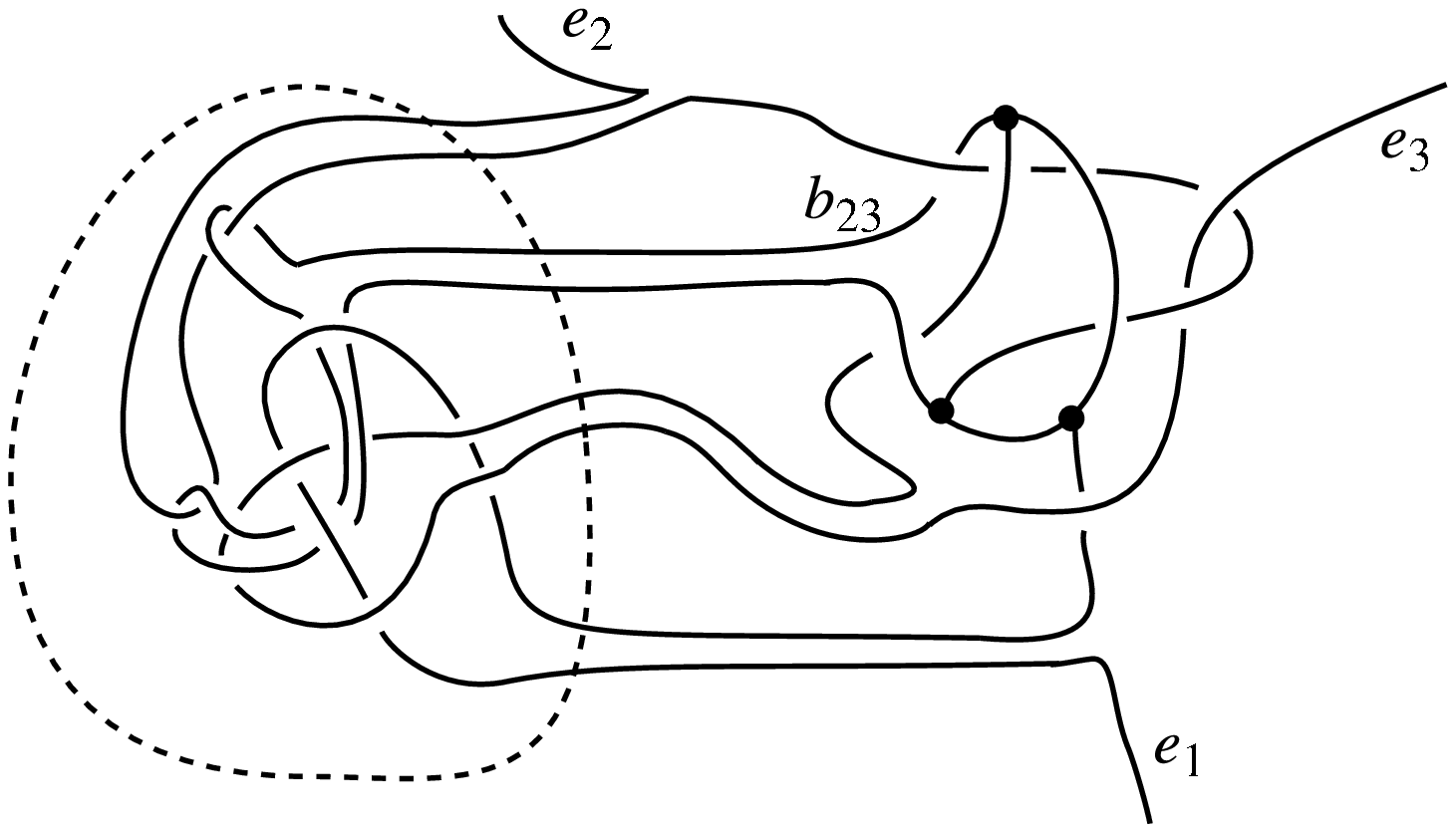}{\Infinitely}		%% 3.9
%

%%%%%%%%%%%%%%%%%%%%%%%%%%%%%%%%%%%%%%%%%%%%%%%%%%%%%%%%%%
%%%%%%%%%%%%%%%%%%%%%%%%%%%%%%%%%%%%%%%%%%%%%%%%%%%%%%%%%%

\subsection{Wagon wheel graphs and tetrahedral graphs} 

Attaching a 3-ball, $B^\prime$, to the 3-ball $B^3$ in which $G$ lies,
then collapsing $B^\prime$ to a point, gives a new graph, $\widehat{G}$,
in the 3-sphere. $\widehat{G}$ is a tetrahedral graph in $S^3$, i.e. a graph
abstractly homeomorphic to the 1-skeleton of a tetrahedron. The edges of $\widehat{G}$ are those of $G$: 
$e_1,e_2,e_3,b_{12},b_{23},b_{31}$.  The
planarity of  $G$ (i.e. lying on a disk in $B^3$) then corresponds to 
the planarity of $\widehat{G}$ in $S^3$ (i.e. lying on a 2-sphere in $S^3$). Also, the exteriors of $G$ and $\widehat{G}$,
in $B^3$ and $S^3$ (resp.), are homeomorphic.

\begin{defn} A surface $F$ 
in a 3-manifold 
$M$ is {\em compressible} if there is a disk $D$ embedded in $M$ such 
that $D\cap F= \partial D$ and $\partial D$ does not bound a disk in $F$. 
$F$ is {\em incompressible} otherwise (where $F$ is not a 2-sphere).
\end{defn}

Recall a graph is abstractly planar in $S^3$ if there exists an embedding of it 
onto a 2-sphere in $S^3$.
We are now in a position to use the following characterization of 
planar graphs in the 3-sphere.

\begin{thm}\label{Thompson} {(Theorem 2.0 of \cite{T}, see also \cite{ST})}
Let $\widehat{G}$ be an abstractly planar graph embedded in $S^3$. $\widehat{G}$ is planar
if and only if:
\begin{enumerate}
\item Every proper subgraph of $\widehat{G} $ is planar, and
\item $X(\widehat{G})=S^3 - \hbox{int}(N(\widehat{G}))$ has compressible boundary.
\end{enumerate}
\end{thm}

\begin{lemma}\label{G-bij}  Let $e_n$, $n = 1, 2, 3$, be three edges of the tetrahedral 
graph $\widehat{G}$ which share a common vertex and suppose $\widehat{G} - e_n$ is 
planar for $n = 1, 2, 3$.  Suppose there is a fourth edge, $b_{jk}$, such that 
$X(\widehat{G} - b_{jk})$ has compressible boundary.  Then   $\widehat{G} - b_{jk}$ is 
a planar graph in
$S^3$.
\end{lemma}

\begin{proof}
We can think of $\widehat{G}-b_{jk}$ as a theta-curve graph, 
$\widehat{\Theta}$, with two vertices and three edges $e_i$,
$e_j \cup b_{ij}$, and $e_k \cup b_{ki}$. 
Since $\widehat{G} - e_n$ is planar for $n = 1, 2, 3$, the subgraphs 
$(\widehat{G}-b_{jk}) - e_i$,
$(\widehat{G}-b_{jk}) - (e_j \cup b_{ij})$, and
$(\widehat{G}-b_{jk}) - (e_k \cup b_{ki})$ are all planar subgraphs.
Since $X(\widehat{G} - b_{jk})$ has compressible boundary, $\widehat{G} - b_{jk}$ is a 
planar graph in
$S^3$ by Theorem \ref{Thompson}.  
\end{proof}

A tetrahedral graph, its corresponding wagon wheel graph, and the tangles carried 
by 
the wagon wheel graph are all related by their exterior.

\begin{defn} 
Let $\T$ be a 3-string tangle with strings $s_{12},s_{23},s_{31}$. 
$X(\T) = B^3 - \text{{\em nbhd}}(s_{12} \cup s_{23}\cup s_{31} \cup \partial B^3)$.
\end{defn}

\begin{lemma}\label{lem:TcarriedbyG}
If $\T$ is carried by the wagon wheel graph $G$ and if $\widehat{G}$ is the 
corresponding tetrahedral graph, then $X(G) = X(\widehat{G})$ is isotopic 
to 
$X(\T)$ and $X(G-b_{ij}) = X(\widehat{G}-b_{ij})$ is 
isotopic to $X(\T-s_{ij})$.  
\end{lemma}

\begin{proof}
See Figure~\solutiontangleA. 
\end{proof}

\subsection{Rational tangles}
$\T$ is {\em rational} if there is an isotopy of $B^3$ taking $\T$ to a tangle containing no crossings.
Rational tangle solutions are generally believed to be
the most likely biological models. For simplification we model a protein-DNA complex
as a 3-ball (the protein) with  properly embedded strings (the DNA segments). However, it is believed 
in many cases that 
the DNA winds around the protein surface without crossing itself. Thus, if we push the
DNA inside the ball, we have a rational tangle.

\begin{lemma}\label{ratG-bij}  Suppose a rational tangle $\T $ is carried
 by $G$, then $X(G)$ has compressible 
boundary and $X(G - b_{jk})$ has compressible boundary for all $b_{jk}$.
\end{lemma}

\begin{proof}
If $\T$ is rational, then there is a disk, $D$, separating $s_{ij}$ from $s_{ik}$.  
Thus $X(G)$ has compressible boundary.
This disk $D$ also separates the two strands of $\T - s_{jk}$ and thus $X(\T - 
s_{jk})$ has
compressible boundary.
Since $X(\T - s_{jk})$ and $X(G-b_{jk})$ are homeomorphic, $X(G-b_{jk})$ has 
compressible boundary. 
\end{proof}

\begin{cor}\label{ratT-bij} Suppose ${\T}$ is rational and ${\T}$ is a solution tangle, 
then ${\T}$ is the PJH tangle. \end{cor}

\begin{proof}
By  Theorem \ref{thm:solutiontangle}, $\T$ is carried by a solution
 graph $G$.  Thus, $G - e_n$ is planar for $n 
= 1, 2, 3$.  By Lemma  \ref{ratG-bij}, $X(G)$ has compressible boundary and $X(G - 
b_{jk})$ has compressible boundary for all $b_{jk}$. Lemmas~\ref{lem:TcarriedbyG} and \ref{G-bij} imply $\widehat{G} - 
b_{jk}$ (where $\widehat{G}$ is the tetrahedral graph corresponding to $G$) is
planar for all $b_{jk}$, and hence $\widehat{G}$ is planar (Theorem \ref{Thompson}).  Thus $G$ is planar.  By Lemma 
\ref{lem:standardT}, $\T$ is standard and hence the PJH tangle 
(Lemma~\ref{lem:standard}).
\end{proof}

Proposition \ref{prop:synapsetanglegeneral} allows us to extend the results of this section to experiments in which the 
{\it in cis} deletion products are $(2,p)$ torus links (rather than specifically (2,4)). 
For example, the analog of Corollary \ref{ratT-bij} would state
that if the synaptic tangle is rational, then three $(2,p)$ torus link {\em in cis}
products would force it to be standard (plectonemic form) as in Figure~\standard.

In contrast to Corollary \ref{ratT-bij}, the examples in Figures ~\counter ~and ~\fouredges, as 
discussed 
in the proof of Theorem \ref{thm:exp}, show that no other set of three deletion experiments would be 
enough to determine that a solution ${\T}$ is standard just from the knowledge that it is rational.

 \begin{thm}\label{thm:exp}
 A system of four deletion experiments with $(2, p)$ torus link products must include 
all three {in cis} experiments to conclude that a rational tangle satisfying  this system 
is standard.
 \end{thm}

\begin{proof}
For the  tetrahedral graph in Figure~\counter,~$\widehat{G} - e_1$, $\widehat{G} - 
b_{12}$, $\widehat{G} - b_{31}$, $\widehat{G} - b_{23}$ are all planar subgraphs.  
Deleting a neighborhood of the vertex adjacent to the $e_i$'s results in a wagonwheel 
graph, $G$.  The subgraph ${G} - e_2$ carries a 2-string tangle which is not rational.  
Hence ${G} - e_2$ is not planar and thus $G$ is not planar, yet it does carry a rational tangle.  
Thus three {\em in trans} experiments and one {\em in cis} 
experiments is not sufficient to determine ${\T}$ is standard under the assumption that $\T$ is rational.

%%Recall that for the 
%%non-planar tetrahedral graph in Figure~\counter,~$\widehat{G} - e_1$, $\widehat{G} - 
%%b_{12}$, $\widehat{G} - b_{31}$, $\widehat{G} - b_{23}$ are all planar subgraphs.  
%%Deleting a neighborhood of the vertex adjacent to the $e_i$'s results in a wagonwheel 
%%graph which is not planar -- yet which carries a rational tangle that is not standard. 
%%Thus three {\em in trans} experiments or two {\em in trans} and one {\em in cis} 
%%experiments is not sufficient to determine ${\T}$ is standard.

Instead of removing a neighborhood of the vertex adjacent to the $e_i$'s, we can delete a neighborhood 
of the vertex adjacent to the edges $e_2$, $b_{12}$, and $b_{23}$ from the tetrahedral graph in Figure 
\counter.  In this case, we also obtain a wagon wheel graph carrying a rational tangle which is not 
standard.  Hence the nonplanar tetrahedral graph in Figure~\counter~can also be used to show that two 
{\em in trans} and two {\em in cis} experiments, where three of the four corresponding edges in the 
tetrahedral graph share a vertex, is also not sufficient to determine ${\T}$ is standard.

The tetrahedral graph in Figure \fouredges ~is not planar as $\widehat{G}-e_2-b_{31}$ is knotted.  
However,
 $\widehat{G} - e_1$,  $\widehat{G} - e_3$, 
$\widehat{G} - b_{23}$, and $\widehat{G} - 
b_{12}$ are all planar. Hence  this graph shows that the remaining case, two {\em in trans} and two {\em 
in cis} experiments where no three of the four corresponding edges in the tetrahedral graph share a common vertex, is also not sufficient to determine ${\T}$ is standard.

\doublefig{1.25}{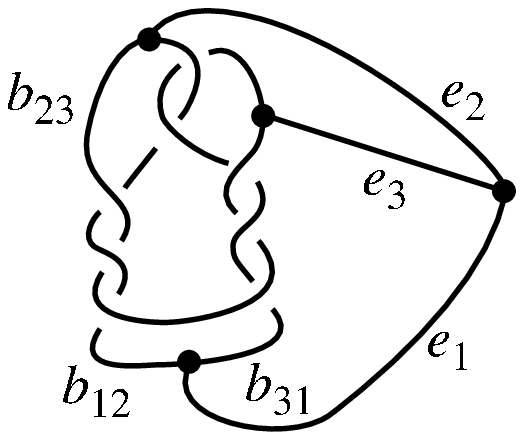}{\counter}{1.25}{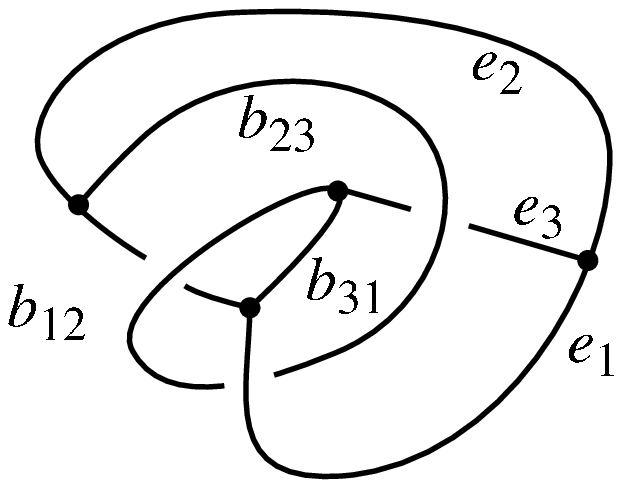}{\fouredges} 

\end{proof}
%%%%%%%%%%%%%%%%%%%%%%%%

\subsection{The exterior X(G)}

Our goal is to prove Corollary~\ref{corPJH}, which extends Corollary \ref{ratT-bij}.
%%that gives conditions under
%%which any solution tangle must be standard. 
In Section 4,
we will use Corollary ~\ref{corPJH}
to show there is a unique  small crossing solution tangle.  
%%By Lemma~\ref{lem:standardT} and Theorem 
%%\ref{thm:solutiontangle}, 
%%it suffices to show that the corresponding solution graph must be planar.
%%To this end 
We first establish some properties of the exterior of both solution
graphs and {\em in trans} solution graphs.

%%%%%%%%%%%%%%%%%%%%%%%%%%%%%%%%%%%%%%%%%%%%%%%
%%%%%%%%%%%%%%%%%%%%%%%%%%%%%%%%%%%%%%%%%%%%%%%

%%\begin{defn} A surface $F$ ($F$ is not a 2-sphere) in a 3-manifold 
%%$M$ is {\em compressible} if there is a disk $D$ embedded in $M$ such 
%%that $D\cap F= \partial D$ and $\partial D$ does not bound a disk in $F$. 
%%$F$ is {\em incompressible} otherwise.
%%\end{defn}

\begin{defn}
Let $M$ be a 3-manifold, $F$ a subsurface of $\partial M$, and $J$ a 
simple closed curve in $F$. 
Then $\tau (M,J)$ is the 3-manifold obtained by attaching a 2-handle to $M$ 
along $J$. 
That is, $\tau (M,J) = M\cup_J H$, where $H$ is a 2-handle. 
$\sigma (F,J)$ is the subsurface of $\partial\tau (M,J)$ obtained by 
surgering $F$ along $\partial H$. 
That is, $\sigma (F,J) = (F\cup \partial H) - \text{int}(F\cap \partial H)$.
\end{defn}

\begin{HandleLemma} %[Handle Addition Lemma]
{(\cite{J},\cite{CG},\cite{Sch}, \cite{Wu1})}
Let $M$ be an orientable, irreducible 3-manifold and $F$ a 
surface in $\partial M$. 
Let $J$ be a simple, closed curve in $F$. Assume $\sigma(F,J)$
is not a 2-sphere.
If $F$ is compressible in $M$ but $F-J$ is incompressible, then 
$\sigma (F,J)$ is incompressible in $\tau (M,J)$.
\end{HandleLemma}

\begin{lemma}\label{cancelb23} 
Let $G$ be a wagon wheel graph and $X(G)$ its exterior. Let $J_i$ and $\gamma_{ij}$ be the meridian curves on $\partial X(G)$ in 
Figure~\NG ~ (and defined right above the figure). Suppose $\partial X(G)$ compresses in $X(G)$.  If $G$ is a solution graph, then $\partial X(G) - (J_1 \cup J_2\cup J_3)$ is compressible in $X(G)$. If $G$ is an {in trans} solution graph, then 
$\partial X(G) - (\gamma_{23}\cup J_1 \cup J_2\cup J_3)$ is compressible in $X(G)$.
\end{lemma}

\begin{proof}
Assume $\partial X(G)$ compresses in $X(G)$. Note that X(G) is irreducible.  We will use the following four claims to prove this lemma:

\begin{claim}\label{f1} 
If $G - b_{23}$ is planar (e.g., $G$ is an in trans solution graph), then
$\partial X(G) - \gamma_{23}$ is compressible in $X(G)$.
\end{claim}

\begin{proof}
Since {$G - b_{23}$ is planar,} $\tau (X(G) ,\gamma_{23}) = X(G-b_{23})$ 
has compressible boundary.  
The claim now follows from the 
Handle Addition Lemma.
\end{proof}

If $G$ is an {\em in trans} solution graph, let $F_1 = \partial X(G) - \gamma_{23}$.  
By Claim \ref{f1}, $F_1$ is compressible 
in $X(G)$.
Else if $G$ is a solution graph, let $F_1 = 
\partial X(G)$.  In this case, $F_1$ is compressible
in $X(G)$ by hypothesis.

\begin{claim}\label{f2}
$F_2 = F_1- J_1$ is compressible in $X(G)$.
\end{claim}

\begin{proof}
$F_1$ is compressible in $X(G)$. 
On the other hand, since $\tau (X(G),J_1) = X(G-e_1)$ and $G-e_1$ is planar, 
$\sigma (F_1,J_1)$ compresses in $\tau (X(G),J_1)$, see Figure~\HALemfig.
The claim now follows from the Handle Addition Lemma.
\end{proof}
\myfig{1.5}{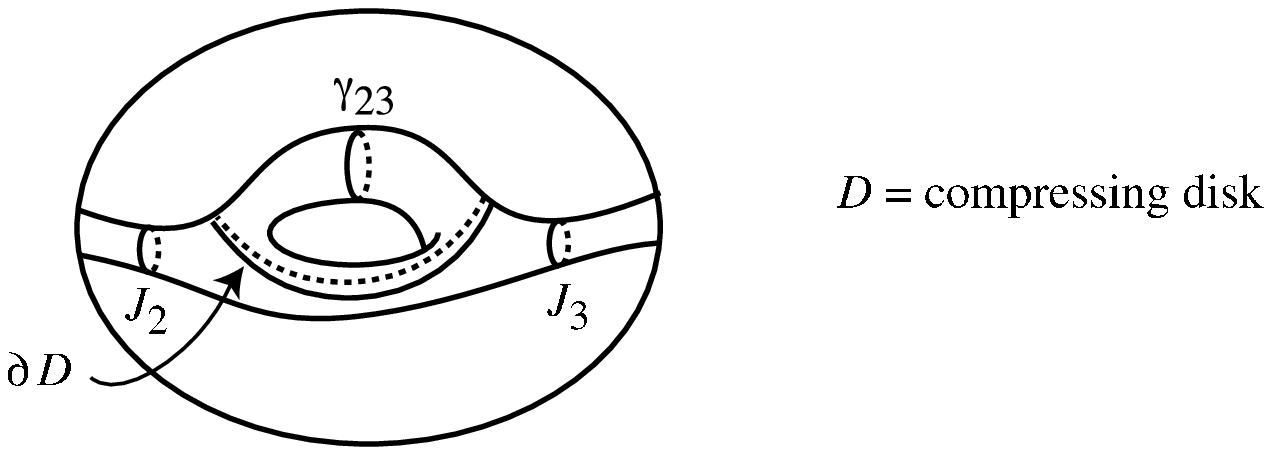}{\HALemfig}		%% 3.11

\begin{claim}\label{f3}
$F_3 = F_2 - J_2$ is compressible in $X(G)$.
\end{claim}

\begin{proof}
Since $G- e_2$ is planar, $\sigma (F_2,J_2)$ is compressible in 
$\tau (X(G),J_2)$ (as in the analog of Figure~\HALemfig). 
Thus Claim~\ref{f2} with the Handle Addition Lemma implies the claim.
\end{proof} 

\begin{claim}\label{f4}
$F = F_3  - J_3$
is compressible in $X(G)$. 
\end{claim}

\begin{proof}
$G- e_3$ is planar, hence $\sigma (F_3,J_3)$ is compressible 
in $\tau (X(G),J_3) = X(G-e_3)$ (as in the analog of Figure~\HALemfig). 
Claim~\ref{f3} with the Handle Addition Lemma says that 
$F_3  - J_3$ is compressible in $X(G)$. 
\end{proof}

Thus if $G$ is a solution
graph, then $F = \partial X(G) - (J_1 \cup J_2\cup J_3)$ is compressible in $X(G)$. 
If $G$ is an {\em in trans}
solution graph, then $F = \partial X(G) - (\gamma_{23}\cup J_1 \cup J_2\cup J_3)$ is 
compressible in
$X(G)$.
\qed (Lemma \ref{cancelb23})
\renewcommand{\qed}{}
\end{proof}

%%%%%%%%%%%%%%%%%%%%%%%%%%%

\subsection{Solution graphs}

\begin{lemma}\label{lem:cis_disk}
Let $G$ be a solution graph and $X(G)$ its exterior. 
Let $J_i$ and 
$\gamma_{ij}$ be the meridian curves on $\partial X(G)$ in Figure~\NG 
~(and defined right above the figure). If $\partial X(G)$ compresses in $X(G)$,
then there is a properly embedded disk in $X(G)$ that intersects
each $\gamma_{ij}$ algebraically once in 
$\partial X(G)$ and
is disjoint from the $J_i$.
\end{lemma}

\begin{proof}
Let $D'$ be a disk, guaranteed by Lemma \ref{cancelb23}, properly embedded in $X(G)$ 
such that $\partial D' \subset F = \partial X(G) - (J_1 \cup J_2\cup J_3)$ and 
$\partial D'$ does 
not bound a 
disk in $F$. 
Write $B^3 = X(G)\cup_{\partial X(G)} M$ where $M= \text{{\em nbhd}}(G\cup 
\partial B^3)$. 
In the notation of Figure~\NG, 
let $M- \text{{\em nbhd}}(\J_1\cup\J_2 \cup\J_3) = M_1\cup M_2$ 
where $M_1$ is $S^2\times I$ and $M_2$ is a solid torus. If  
$\partial D'$ lies in $\partial M_1$, we can form a 
2-sphere, $S$, by capping off $D'$ with a disk in $M_1$ (take the disk 
bounded by $\partial D'$ on $\partial M_1$ and push in slightly). 
That is, $S$ is a 2-sphere which intersects $\partial X(G)$ in $\partial D'$.
But then $S$ would have to be 
non-separating (there would be an arc from $\J_1$ to $\J_2$, say, on 
$\partial M_1$ which intersected $\partial D'$ once and an arc connecting
$\J_1$ to $\J_2$ on $\partial M_2$ which missed $\partial D'$, hence $S$, 
altogether). 
But 2-spheres  in a 3-ball are separating.

Thus $\partial D'$ must lie in $\partial M_2$.
%%If $\partial D'$ is essential in $\partial M_2$, $D'$ is the desired disk
%%and we are done. 
If $\partial D'$
	is essential in $\partial M_2$, then $\partial D'$ must intersect
	any meridian of the solid torus $M_2$ algebraically once (else
	$B^3$  would contain a lens space summand). Then $D'$ is the 
	desired disk, and we are done.
So assume $\partial D'$ bounds a disk, $D''$, 
in $\partial M_2$. Then the disks $\text{{\em nbhd}}(\J_1\cup\J_2 \cup\J_3) 
\cap \partial M_2$ must lie in $D''$, as $S = D' \cup D''$ is separating.
Since $G - e_1$ is planar, $M_2$ is an unknotted solid torus in 
$B^3$. Thus there is a properly embedded disk 
$D \subset {B^3 - \text{int} M_2}$ such that $\partial D$ is a longitude of $\partial M_2$ and 
hence
is essential. Furthermore, we can arrange that $D$ is disjoint from $D''$.
Then $D \cap S$ is a collection of trivial circles in $D'$ and $D$ can be
surgered off $S$ without changing $\partial D$. In particular, $D$ is properly
embedded in $X(G)$ since the $e_i$ lie on the opposite side of $S$ from $D$.
$D$ is the desired disk.

\end{proof}

%%%%%%%%%%%%%%%%%%%%%%%%%%%%%%%%%%%%%%

\begin{defn}
Let $\T$ be a 3-string tangle. 
%%$\T$ is {\em rational} if there is an isotopy of $B^3$ taking $\T$ to the 
%%standard tangle of Figure~\splitrationaltang (a). 
$\T$ is {\em split} if there is a properly embedded disk separating 
 two of its strands from the third 
(Figure~\splitrationaltang (a)). 
Strands $s_1,s_2$ of $\T$ are {\em parallel}  if there is a disk $D$ in 
$B^3$ such that $\text{int}(D)$ is disjoint from $\T$ and 
$\partial D = s_1 \cup \alpha \cup s_2 \cup\beta$ where $\alpha,\beta$ 
are arcs in $\partial B^3$ (Figure~\splitrationaltang (b)). 
\myfig{1.75}{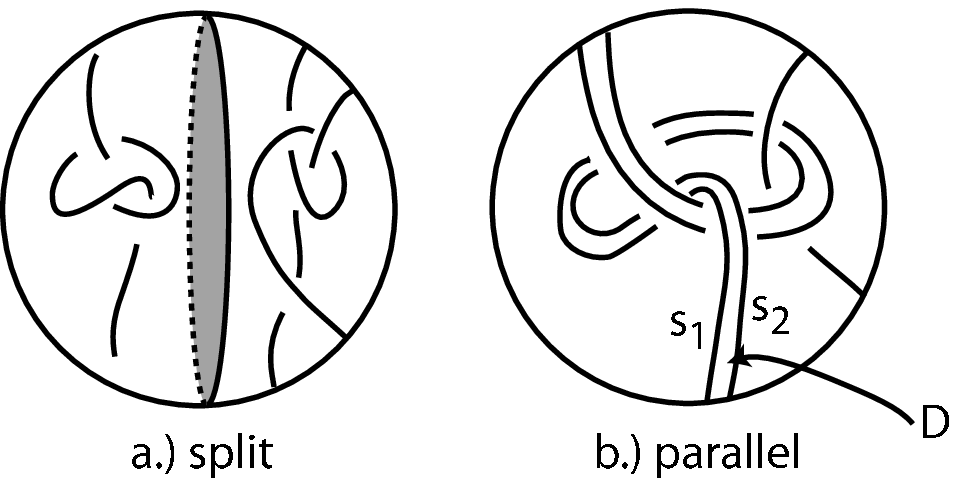}{\splitrationaltang}		%% 3.13
Note: Rational tangles are split and have parallel strands.
\end{defn}

\begin{lemma}\label{lem:lemPJH}     %% 3.7 
Let $\T$ be a 3-string tangle carried by a solution graph $G$ . 
If $\T$ is rational, split, or if $\T$ has parallel strands{,} 
then $X(G)$ and $X(G-b_{ij})$, for each $i,j$, have compressible boundary.
\end{lemma}

\begin{proof}
If $\T$ is rational, split, or has parallel strands then 
$\partial X(\T)$ is compressible in $X(\T)$. 
By Lemma~\ref{lem:TcarriedbyG}, $\partial X(G)$ is compressible in  $X(G)$.

We now argue that $\partial X(G-b_{ij})$ is compressible in $X(G-b_{ij})$.
Below $\{i,j,k\} = \{1,2,3\}$.

If $\T$ is split, then there is a disk, $D$, separating one strand
from two. If $s_{ij}$ is on the side of $D$ containing
two, then $D$ still separates the two strands of $\T - s_{ij}$ and
$X(\T-s_{ij})$ still has compressible boundary. By Lemma~\ref{lem:TcarriedbyG},
$X(G-b_{ij})$ has compressible boundary. So we assume $s_{ij}$ is on 
one side of $D$ and $s_{jk}, s_{ki}$ on the other. As a disk properly
embedded in $X(G)$, $D$ separates the curves $\gamma_{ij},\gamma_{jk}$ of
$\partial X(G)$. By Lemma~\ref{lem:cis_disk}
there is a disk $D^\prime$ in X(G) 
that intersects $\gamma_{jk}$ on $\partial X(G)$
algebraically a non-zero number of times. After
 surgering along $D$ (while preserving the non-zero algebraic intersection 
with $\gamma_{jk}$), we may assume that $D^\prime$ is disjoint from
$D$, and hence from $\gamma_{ij}$. But then $D^\prime$ is a compressing disk 
for $G-b_{ij}$.

If $\T$ has parallel strands, we argue in two cases according to
which strands are parallel:

First, assume $s_{jk},s_{ki}$ are parallel. Then they are parallel in 
$\T - s_{ij}$ and $X(\T - s_{ij})$ has compressible boundary. By 
Lemma~\ref{lem:TcarriedbyG} then $X(G-b_{ij})$ has compressible boundary.

So, assume that $s_{ij}$ and $s_{jk}$ are parallel. Then
there is a disk, $D$, properly embedded in $X(G)$ that intersects 
$\gamma_{ij}$ exactly once, $\gamma_{jk}$ exactly once, and is disjoint 
from $\gamma_{ki}$.
By Lemma~\ref{lem:cis_disk}, there is also a disk, $D^\prime$, in
$X(G)$ that intersects $\gamma_{ki}$ algebraically 
once.
%%a non-zero number of times. 
After possibly surgering
$D^\prime$ along $D$ we may assume the $D$ and $D^\prime$ are disjoint.
But $X(G - \hbox{{\em nbhd}}(D))$ is homeomorphic to $X(G-b_{ij})$, and
$D^\prime$ becomes a compressing disk in $X(G-b_{ij})$.

\end{proof} 

%%%%%%%%%%%%%%%%%%%%%%%%%%%%%%%%%

\begin{cor}\label{corPJH}       %% 3.7 
Let $\T$ be a solution tangle. 
If $\T$ is rational or split or if $\T$ has parallel strands{,} 
then $\T$ is the PJH tangle.
\end{cor}

\begin{proof} Let $G$ be a solution graph carrying $\T$. Then Lemmas~\ref{lem:lemPJH}, \ref{lem:TcarriedbyG} and \ref{G-bij} imply 
that $\widehat{G} - b_{jk}$ planar for all $b_{jk}$ where $\widehat{G}$ is the associated tetrahedral graph.  Thus 
Theorem~\ref{Thompson} implies $\widehat{G}$ is planar.  Thus $G$ is planar. By Lemma~\ref{lem:standardT}, $\T$ is standard.  Now 
Lemma~\ref{lem:standard} says $\T$ is the PJH tangle. \end{proof}

\begin{remark}
Since {\em in trans} solution tangles are also solution tangles, Corollary 
\ref{corPJH} also 
applies to {\em in trans} solution tangles.
Although the three {\em in cis} deletion equations suffice to determine ${\T}$ 
if ${\T}$ is rational, split or has parallel strands, the {\em in trans} deletion 
experiment does rule out 
some exotic solutions.  For example, the tangle in Figure \tripusA~ satisfies 
all the {\em in cis} experimental results ({{\em i.e.} it is a solution tangle), but does not satisfy 
the {\em in 
trans} results ({\em i.e.} it is not an {\em in trans} solution tangle}).

\myfig{1.25}{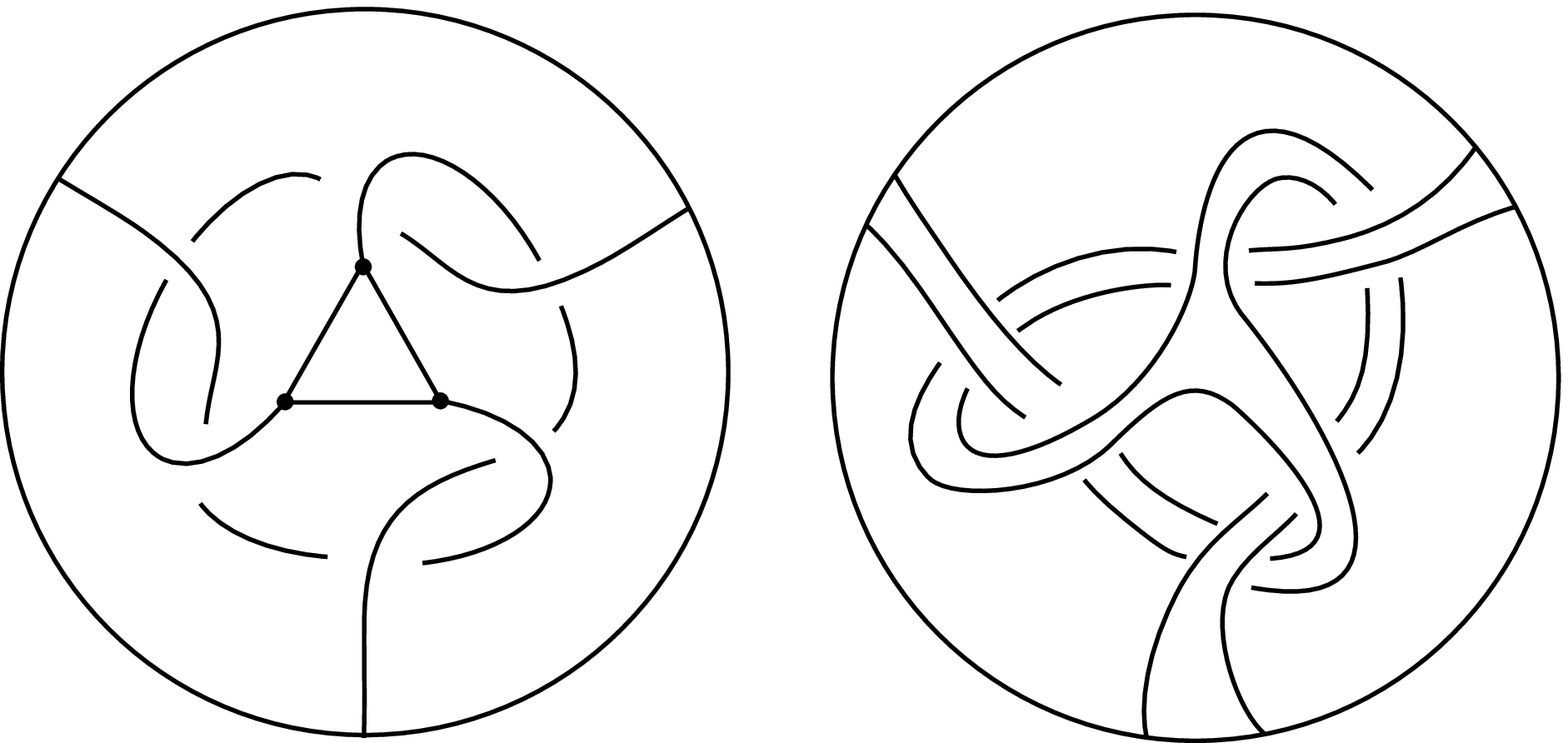}{\tripusA} 
\end{remark}

%%%%%%%%%%%%%%%%%%%%%%%%%%%%

\subsection{{\em In trans} solution graphs}

In Lemma \ref{cancelb23a} we will show that the exterior of an {\em in trans} solution graph is a 
handlebody.  This need not be true for solution graphs. The solution graph of Figure~\tripusA ~contains an incompressible genus 2 surface in its
exterior (it is based on Thurston's {\em tripos} graph). Thus its exterior is not a
handlebody. 
%% as will be concluded for solution graphs in 
%%Lemma \ref{cancelb23a}.
We will also use Lemma \ref{cancelb23a} to further explore planarity of tetrahedral
graphs in the next section.

\begin{lemma}\label{cancelb23a}
Let $G$ be an {in trans} solution graph and $X(G)$ its exterior. Let $J_i$ and 
$\gamma_{ij}$ be the meridian curves on $\partial X(G)$ in Figure~\NG ~
(and defined right above the figure). If $\partial X(G)$ compresses in $X(G)$,
then there is a properly embedded disk in $X(G)$ that intersects
$\gamma_{23}$ exactly once, intersects each of $\gamma_{12}, \gamma_{31}$ algebraically once, and is disjoint from  the 
$J_i$. Furthermore,
$X(G)$ is a handlebody.
\end{lemma}

\begin{proof}

By Lemma \ref{cancelb23}, there is a disk $D'$
properly embedded in $X(G)$ 
such that $\partial D' \subset F = \partial X(G) - (\gamma_{23}\cup J_1 \cup 
J_2\cup J_3)$ and $\partial D'$ does not bound a 
disk in $F$. 
Write $B^3 = X(G)\cup_{\partial X(G)} M$ where $M= \text{{\em nbhd}}(G\cup 
\partial B^3)$. 
In the notation of Figure~\NG, 
let $M- \text{{\em nbhd}}(\J_1\cup\J_2 \cup\J_3\cup\Gamma_{23}) = M_1\cup M_2$ 
where $M_1$ is $S^2\times I$ and $M_2$ is a 3-ball. 
$\partial D'$ lies in $\partial M_i$ for some $i$, hence we can form a 
2-sphere, $S$, by capping off $D'$ with a disk in $M_i$ (take the disk 
bounded by $\partial D'$ on $\partial M_i$ and push in slightly). 
That is, $S$ is a 2-sphere which intersects $\partial X(G)$ in $\partial D'$.
If $\partial D'$ lay on $\partial M_1$, then $S$ would have to be 
non-separating (there would be an arc from $\J_1$ to $\J_2$, say, on 
$\partial M_1$ which intersected $\partial D'$ once and an arc connecting
$\J_1$ to $\J_2$ on $\partial M_2$ which missed $\partial D'$, hence $S$, 
altogether). 
But 2-spheres  in a 3-ball are separating.
\myfig{2.0}{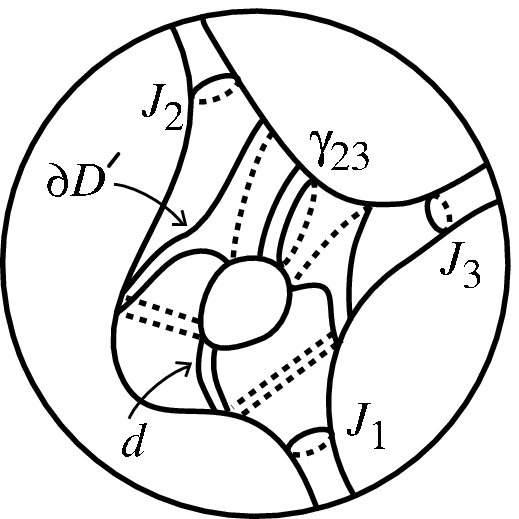}{\GstandardB}		%% 3.12
%\doublefig{2.5}{GstandardA}{\GstandardA}{2.5}{GstandardB}{\GstandardB}
%
Thus $\partial D'$ must lie in $\partial M_2$. 
But then in $\partial M_2$, the disks $\text{{\em nbhd}}(\J_1\cup\J_2 \cup\J_3) 
\cap \partial M_2$ must lie on opposite sides of $\partial D'$ from the
disks $\Gamma_{23}^{\prime} \cup \Gamma_{23}^{\prime \prime} = 
\text{{\em nbhd}}(\Gamma_{23}) \cap \partial M_2$, as $S$ is separating.
Thus $S$ can be isotoped  to intersect $N(G)$ in 
the curve gotten by connecting the two disks $\Gamma_{23}^{\prime}$, 
$\Gamma_{23}^{\prime \prime}$
by a band ${d}$ lying
in $\partial$N(G). See Figure~\GstandardB. We may view S as a 2-sphere
which intersects $G$ in only two 
points and separates off a subarc $b'_{23}$ of $b_{23}$.   
Now $b'_{23}$ must be unknotted in the 3-ball bounded by $S$, 
otherwise deleting $e_1$ would 
not give a planar graph. Hence, there exists a disk $D$ in the 3-ball bounded by $S$ such that the boundary of $D$ 
consists of the arc $b'_{23}$ and an arc in $\Gamma_{23}^{\prime} \cup
\Gamma_{23}^{\prime \prime} \cup {d}.$ Then D can be taken 
 to be a disk,
properly embedded in $X(G)$ whose boundary is disjoint from $J_1\cup J_2
\cup J_3$ and whose boundary intersects $\gamma_{23}$ exactly once.  Finally, as $\gamma_{12}, \gamma_{31}$ are 
isotopic to 
             $\gamma_{23}$ on the boundary of the solid torus component of 
$M - nbhd(\J_1 \cup \J_2 \cup \J_3)$,              
%%$\text{nbhd}(\J_1 \cup \J_2 \cup \J_3) \cap \partial M_2$,
             $\partial D$ intersects these curves algebraically once. 
Furthermore, $X(G \cup D)$ is homeomorphic to $X(G-\gamma_{23})$, which 
is a handlebody. Thus $X(G)$ is a handlebody.

\end{proof}

%%%%%%%%%%%%%%%%%%%%%%%%%%%%%%%%%%%%%%%%%%%%%%%%%
%%%%%%%%%%%%%%%%%%%%%%%%%%%%%%%%%%%%%%%%%%%%%%%%%

%%%%%%%%%%%%%%%%%%%%%%%%%%%%%%%%%%%%%%%%%%%%%%%%%
%%%%%%%%%%%%%%%%%%%%%%%%%%%%%%%%%%%%%%%%%%%%%%%%%

\subsection{Tetrahedral graph planarity}
In the above model, each mathematical assumption about the planarity of the {\em in trans}
solution graph after deleting an edge corresponds to an experimental result.
For the economy of experiment, then, it is natural to ask how much information
is necessary to conclude the planarity of the graph. 
Recall that by collapsing the 
outside 
2-sphere to a vertex, a wagon wheel graph in $B^3$ becomes a {tetrahedral} 
graph $\widehat{G}$  
in $S^3$.
We discuss this issue of economy in the context of tetrahedral graphs and in the 
spirit of Theorem~\ref{Thompson} (see also \cite{ST}, \cite{Wu2}, \cite{G}).
In particular, Theorem \ref{thm:oppositeedges} and  Corollary \ref{tp} say that, in contrast to 
Theorem~\ref{Thompson}, to check the planarity of a tetrahedral graph, 
one does not need to check the planarity of {\em all\/} subgraphs.
We end with examples showing that Corollary \ref{tp}  is sharp.

\begin{thm}\label{thm:oppositeedges}
Let $\widehat{G}$  be a tetrahedral graph.  	
Then $\widehat{G}$  is planar if and only if
\begin{enumerate}
\item
the exterior of $\widehat{G}$  has compressible boundary; and
\item
there is an edge, $\epsilon^\prime$, of $\widehat{G}$  such that 
for any edge $e \neq \epsilon^\prime$ of $\widehat{G}$ , $\widehat{G}-e$ is planar.
\end{enumerate}
\end{thm}

\begin{proof}
If $\widehat{G}$  is planar then clearly the two conditions hold. We prove the 
converse: the two listed conditions on $\widehat{G}$  guarantee that $\widehat{G}$  is planar. 
Let $\epsilon$ be the edge of $\widehat{G}$ which does not share an endpoint
with $\epsilon^\prime$. 

Let $e_1, e_2, e_3 , e_4$ be the four edges of $\widehat{G}$  other than $\epsilon,
\epsilon^\prime$. Let $d_i$ be the meridian disk 
corresponding to the edge $e_i$ in {\em nbhd}$(\widehat{G})$ and let $m_i = \partial d_i$
be the corresponding meridian curve on the 
boundary of $X(\widehat{G})$, the exterior of $\widehat{G}$ .

\begin{claim}\label{handlebody}
$X(\widehat{G})$  is a handlebody.
\end{claim}

\begin{proof}
Pick a vertex, $v$, of $\widehat{G}$  not incident to $\epsilon^\prime$. Then removing a 
neighborhood of $v$ from $S^3$, $\widehat{G}$  becomes an {\em in trans} solution graph $G$
in a 3-ball. Furthermore, $X(G)$ is homeomorphic to $X(\widehat{G})$ . The
Claim now follows from Lemma~\ref{cancelb23a}.
\end{proof}

Recall the following from \cite{G}.

\begin{thm}\label{thm:gordon}
Let $\mathcal C$ be a set of $n+1$ disjoint simple loops in the boundary
of a handlebody X of genus n, such that $\tau(X;{\mathcal C}^\prime)$ is
a handlebody for all proper subsets ${\mathcal C}^\prime$ of $\mathcal C$.
Then $\cup \, \mathcal C$ bounds a planar surface $P$ in $\partial X$ such
that $(X,P) \cong (P \times I, P \times {0})$, where $I = [0,1]$.
\end{thm}

Set $X = X(\widehat{G})$ and ${\mathcal C} = \{m_1, m_2, m_3, m_4\}$. Since $\widehat{G} - e_i$ is planar 
for $i = 1, ..., 4$,
Theorem~\ref{thm:gordon} implies that there is a planar surface $P$
in $X(\widehat{G})$  whose boundary is $\{m_1, m_2, m_3, m_4\}$ and such that
$(X(\widehat{G}),P) \cong (P \times I, P \times {0})$. Capping $P \times {1 \over 2}$ with
the meridian disks $\{ d_1, d_2, d_3, d_4 \}$ gives a 2-sphere $\widehat P$
in $S^3$ disjoint from $\widehat{G}$  except for a single point intersection with
each of the edges $e_1, e_2, e_3, e_4$. Let $\widehat{G} - {\widehat P} = G_1 \cup
G_2$. Labeling the edges as in $\widehat{G}$ , we take $G_1$ with edges
$\{\epsilon, e_1, e_2, e_3, e_4 \}$ and $G_2$ with edges
$\{\epsilon^\prime, e_1, e_2, e_3, e_4 \}$. 
Furthermore, we may label so that
$e_1, e_2$ share a vertex in $G_1$ and $e_2,e_4$ share a vertex in $G_2$.
See Figure~\HatP.

\myfig{1.75}{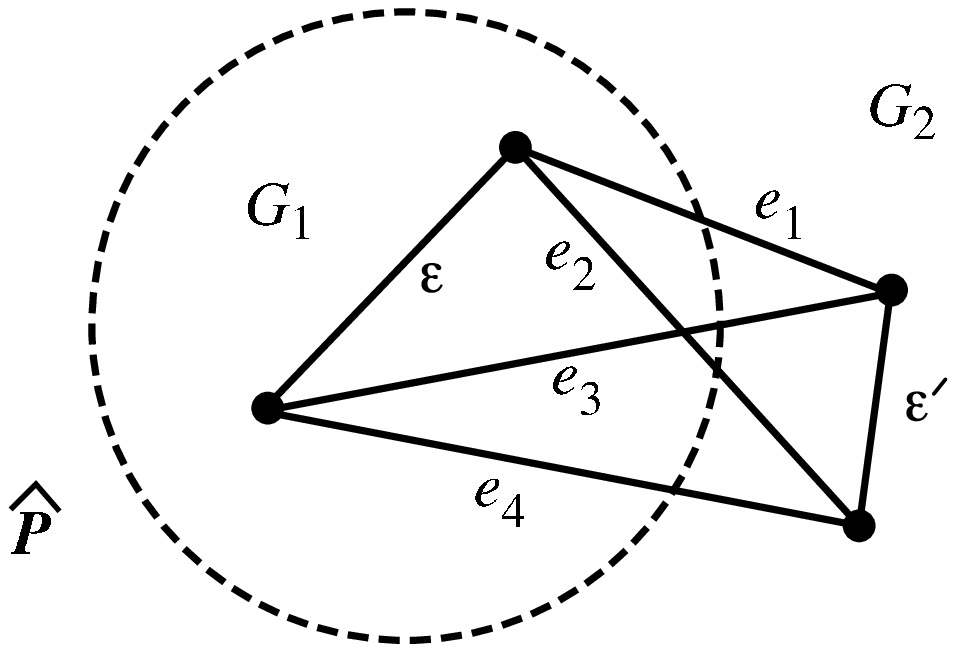}{\HatP}		%% 3.14

Let $N_1 = \hbox{{\em nbhd}}(G_1)$ and $N_2 = \hbox{{\em nbhd}}(G_2)$. Then
{\em nbhd}$(\widehat{G})$ is the union of $N_1$ and $N_2$ along the meridian disks 
$\{d_1, d_2, d_3, d_4\}$. Now the pair $(N_i,\cup d_j) = (F_i \times I,
\cup \alpha_i^j \times I)$ where $i \in \{1,2\}$, $j \in \{1,2,3,4\}$ and
where each $F_i$ is a disk, each $\alpha_i^j$ is an arc in $\partial F_i$.
Each $G_2$ can be taken to lie in $F_i \times \{1/2\}$. 
Write $S^3 = B_1 \cup_{\widehat P} 
B_2$ where the $B_i$ are 3-balls. Because   
$(X(\widehat{G}),P) \cong (P \times I, P \times {0})$, we may extend 
$F_i \times \{1/2\}$ to a 
properly embedded disk $F_i^\prime$ in $B_i$ containing $G_2$. By a proper
choice of the initial $F_i$ we may assume that $F_i^{\prime} - 
(\cup \alpha_i^j)$
is a union of four arcs $\beta_i^j$, labelled as in Figure~\Betas.

%%\smallskip

\myfig{2.0}{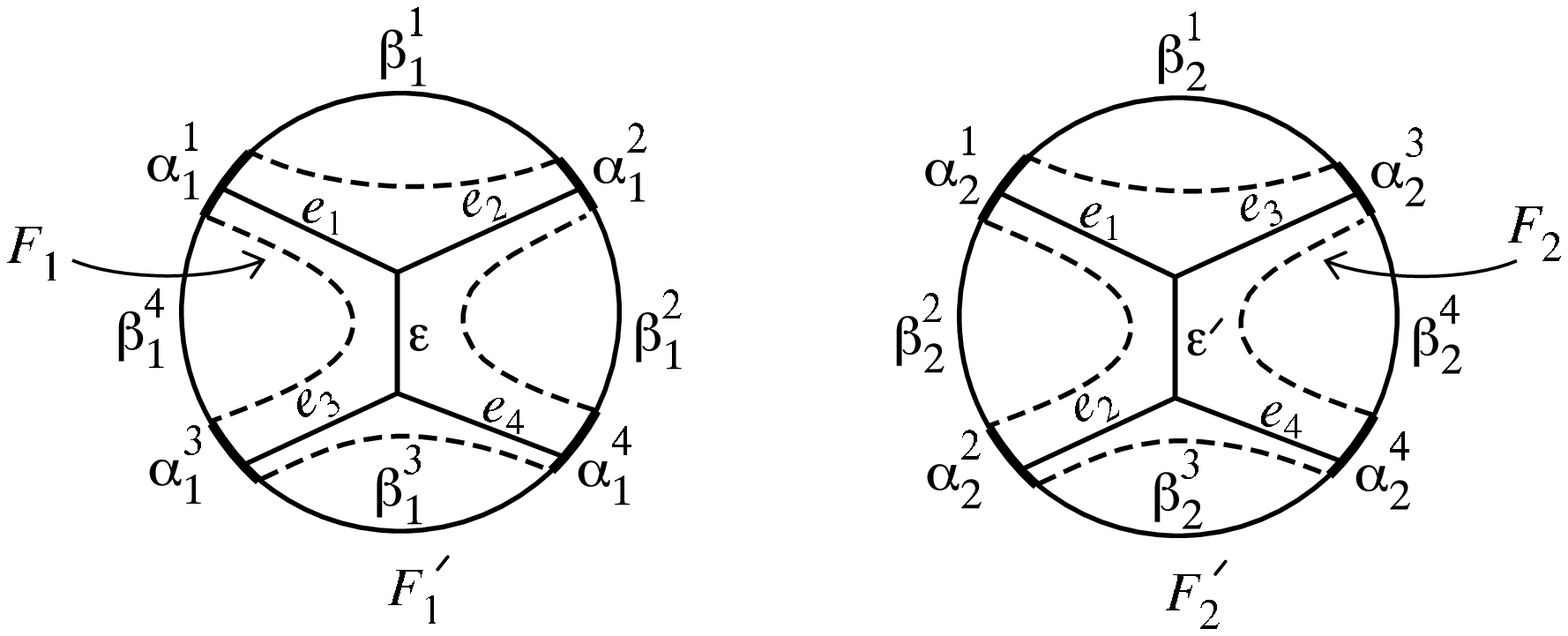}{\Betas}	%% 3.15

In particular, 
\begin{itemize}
\item
$\beta_1^1$ is isotopic to $e_1 \cup e_2$ and $\beta_1^3$ is
isotopic to $e_3 \cup e_4$ in $B_1$ (keeping the endpoints in $\cup d_j$
on $\partial B_1 = \widehat P$).
\item
$\beta_2^1$ is isotopic to $e_1 \cup e_3$ and $\beta_2^3$ is
isotopic to $e_2 \cup e_4$ in $B_2$ (keeping the endpoints in $\cup d_j$
on $\partial B_1 = \widehat P$).
\end{itemize}

The boundary of $F_1^\prime$ frames the 4-punctured sphere $P \subset 
\widehat P$. This allows us to assign a rational number slope on 
properly embedded arcs in $P$ (up to proper isotopy in $P$).
We take $\beta_1^1, \beta_1^3$ to have slope $\frac{0}{1}$ and 
$\beta_1^2, \beta_1^4$ to have slope $\frac{1}{0}$. It follows that
$\beta_2^1, \beta_2^3$ must have slope $\frac{1}{n}$: otherwise the
edges $e_1, e_2, e_3, e_4$ of $\widehat{G}$  will form a non-trivial 2-bridge knot, but
we are assuming that $\widehat{G} - \epsilon$ is planar. 

By re-choosing $F_1$ we may assume that $\beta_1^2$ is 
$\beta_2^3$ and that $\beta_1^4$ is $\beta_2^1$. 
Then re-choosing $F_2$ (by twisting along the disk disjoint from $\beta_2^1 \cup 
\beta_2^3$ in $B_2$), we may further take $\beta_2^2$ to be $\beta_1^1$
and $\beta_2^4$ to be $\beta_1^3$. Then $F_1^\prime \cup F_2^\prime$
is a 2-sphere in $S^3$ containing $\widehat{G}$  (Alternatively, instead of rechoosing $F_2$,  
after rechoosing $F_1$ so that $\beta_1^2 = 
\beta_2^3$ and $\beta_1^4 = \beta_2^1$, 
we can extend $F_1$ to a disk containing $\widehat{G} - \epsilon^\prime$.  Thus $\widehat{G} - 
\epsilon^\prime$ is planar and $\widehat{G}$  is planar by Theorem \ref{Thompson}).
\qed (Theorem \ref{thm:oppositeedges})
\renewcommand{\qed}{}
\end{proof}

We will now use Theorem \ref{thm:oppositeedges} to prove Corollary \ref{tp}.  

\begin{cor}\label{tp}
Suppose $\widehat{G}$  is a tetrahedral graph with the following properties:
\begin{enumerate}
\item  There exists three edges $\epsilon_1$, $\epsilon_2$, $\epsilon_3$ such that $\widehat{G} - \epsilon_i$ is planar.
\item  The three edges $\epsilon_1$, $\epsilon_2$, $\epsilon_3$ share a common vertex.
\item  There exists two additional edges, $\epsilon_{4}$ and $\epsilon_{5}$ such that $X(\widehat{G} - \epsilon_{4})$ and $X(\widehat{G} 
- 
\epsilon_{5})$ 
have compressible boundary.
\item  $X(\widehat{G})$  has compressible boundary.
\end{enumerate}
Then $\widehat{G}$  is planar.
\end{cor}

\begin{proof} By Lemma \ref{G-bij}, $\widehat{G} - \epsilon_{4}$ and $\widehat{G} 
- \epsilon_{5}$ are planar.   Thus by Theorem \ref{thm:oppositeedges}, $\widehat{G}$  is 
planar.
\end{proof}

We finish this section with examples that show that if any single 
hypothesis of Corollary \ref{tp} is dropped then its conclusion, that 
$\widehat{G}$ is planar, may not hold. The tetrahedral graph in Figure 
\fouredges ~is not planar as $\widehat{G}-e_2-b_{31}$ is knotted.  
However,
 $\widehat{G} - e_1$,  $\widehat{G} - e_3$, 
$\widehat{G} - b_{23}$, and $\widehat{G} - 
b_{12}$ are all planar. 
Furthermore, $X(\widehat{G})$ is a handlebody as the edge $b_{23}$
can be isotoped onto a neighborhood of $\widehat{G}-b_{23}$. (Alternatively, we can determine that 
$X(\widehat{G})$ 
is a handlebody by removing the vertex adjacent to the $e_n$'s, creating a wagon wheel graph, 
${G}$, 
which carries a rational tangle $\T$. Since $X(\widehat{G}) = X(G) = X(\T)$ and $\T$ is rational,  
$X(\widehat{G})$ is a handlebody.)
 Thus 
 $X(\widehat{G})$ has compressible boundary.
$X(\widehat{G} - b_{31})$, and $X(\widehat{G} - e_2)$ also have compressible 
boundary.  
Thus properties 1, 3, 4 in Corollary \ref{tp} hold, but property 2 does not hold.

%%\doublefig{1.25}{fouredges}{\fouredges}{1.25}{counter}{\counter} 

 In Figure \counter, $\widehat{G} - e_1$, $\widehat{G} - b_{12}$, $\widehat{G} - 
b_{31}$, $\widehat{G} - b_{23}$ are all planar graphs.  The edges $e_1$, $b_{12}$, 
and $b_{31}$ share a common vertex.  Since $\widehat{G} - b_{23}$ is planar, 
$X(\widehat{G} - b_{23})$ has compressible boundary. Removing the vertex adjacent to 
the $e_n$'s in $\widehat{G}$ results in a wagon wheel graph, $G$, which carries a 
rational tangle.  By Lemma \ref{lem:TcarriedbyG} $X({G})$ is a handlebody, and thus 
$X(\widehat{G})$ has compressible boundary.
 On the other hand, ${G}-e_2$, ${G}-e_3$ carry tangles which are non-trivial disk sums that 
do not have compressible boundary (by Lemma 3.3 of \cite{Wu3}). By Lemma \ref{lem:TcarriedbyG}, 
neither $X({G}-e_2)$ nor $X({G}-e_3)$ have compressible boundary. Hence neither 
$X(\widehat{G}-e_2)$ nor
                $X(\widehat{G}-e_3)$ have compressible boundary. In particular, 
$\widehat{G}$ is not planar. Thus the tetrahedral graph in Figure \counter 
~satisfies properties 1, 2, 4 in Corollary \ref{tp}, but only satisfies half of 
property 3 ($X(\widehat{G} - b_{23})$ has compressible boundary, but there does not 
exist a fifth such edge).

Note that the tangle, $\T$, in Figure \aaa ~is carried by a nonplanar wagon wheel 
graph, $G$.  If $G$ were planar,
 then $\T$ would be the PJH tangle by Lemmas \ref{lem:standardT} and 
\ref{lem:standard}. To see that $\T$ is not the PJH tangle, cap off these tangles by 
arcs along the tangle circle that are complementary to the $c_i$ ({\it i.e.} the 
$x_{ij}$ of Figure \tangcircA) to form 3 component links.  These two 3-component 
links are not isotopic.  For example, they have different hyperbolic volumes as 
computed by SnapPea.
 Let $\widehat{G}$ be the tetrahedral graph corresponding to $G$.  Every subgraph of $\widehat{G}$ is 
planar.  Hence conditions 1, 2, 3 of Corollary \ref{tp} hold, but $X(\widehat{G})$ does not have 
compressible boundary (by Theorem \ref{Thompson}).

\myfig{1.0}{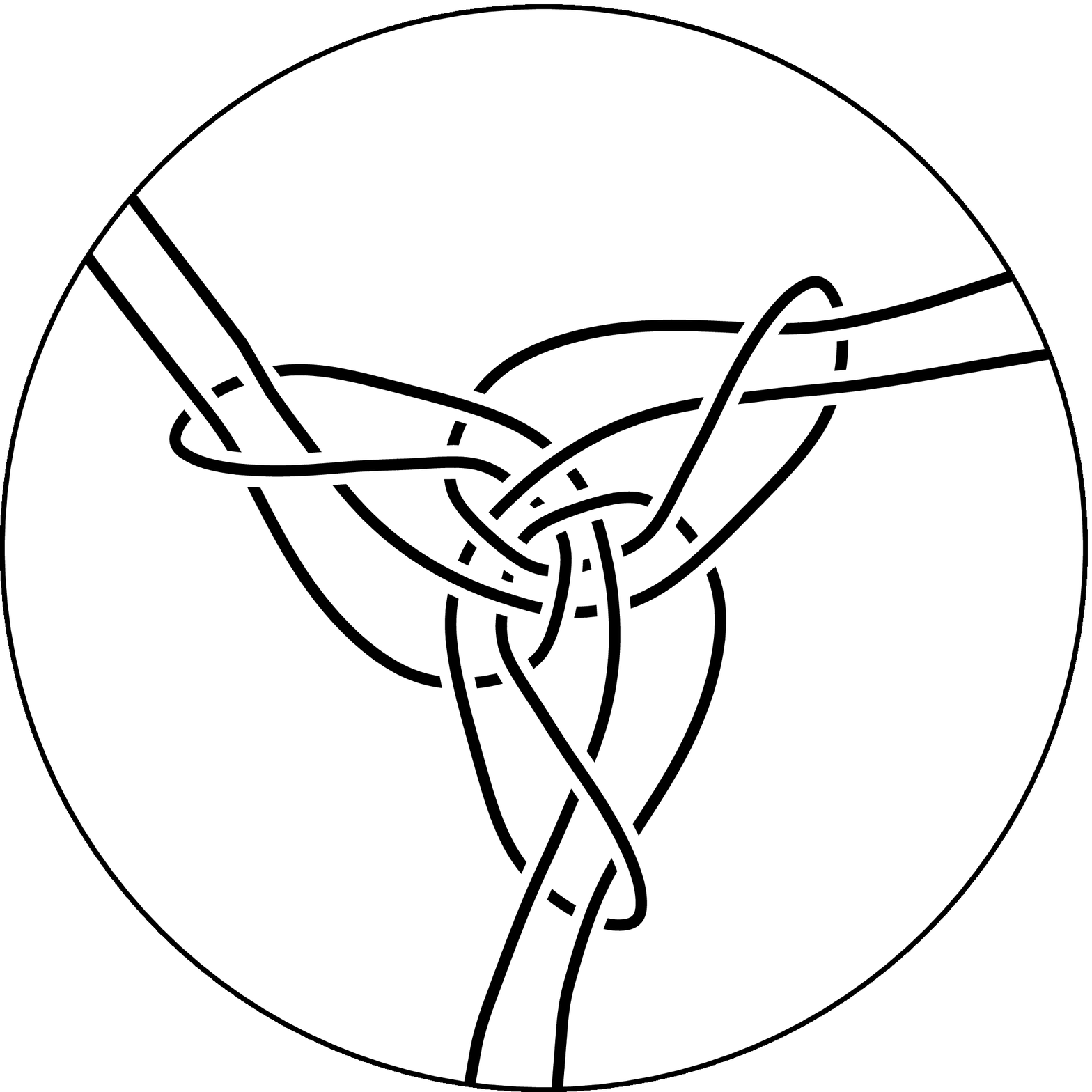}{\aaa}      %% 3.15

%% file: section4.tex
%%%%%%%%%%%%%%%%%%%%%%%%%%%%
\section{If $\T$ is not the PJH tangle, then $Cr(\T) \ge 8$ up to free isotopy}	%% 4

\begin{defn}
Two 3-tangles $\T_1,\T_2$ are {\em freely isotopic} if there is an 
isotopy of the 3-ball, which is not necessarily fixed on its boundary, 
taking $\T_1$ to $\T_2$.
\end{defn}

In the last section we showed 
that if a solution tangle is split or has two parallel strands, then it must be the PJH tangle. Here we 
show that if a tangle can be freely 
isotoped to have fewer than eight crossings, then it is either split 
or has two parallel strands.  Hence a solution tangle which can be freely isotoped to have fewer than eight crossings 
must be the PJH tangle (Corollary \ref{cor:PJH7}). As a rational tangle is one which
            can be freely isotoped to have no crossings, one can think of 
Corollary \ref{cor:PJH7} as a generalization of the result (Corollary \ref{corPJH})
            that the only rational solution tangle is the PJH tangle. Corollary
            \ref{cor:PJH7} also gives a lower bound on the crossing number of an exotic
            solution tangle. Using the {\em in trans} deletion result, we will improve this lower bound in the next
            section after imposing a normal framing on the tangle (working
            with tangle equivalence rather than free isotopy).

\begin{lemma}\label{lem:4.1}
If one string of a 3-tangle $\T$ crosses the union of the other two 
strings at most once, then $\T$ is split.
\end{lemma}

\begin{proof}
If the string passes over the union, isotop the strand to the front  hemisphere
of the tangle sphere, otherwise isotop to the back.
\end{proof}

\begin{defn}
To a projection of a 3-string tangle $\T$ we associate the 4-valent graph 
$\Gamma (\T)$, that is {obtained} 
by placing a vertex at each crossing. 
If $\T$ is not split, 
we label in sequence $e_1,\dots,e_6$ the distinct edges
which are incident to the tangle circle.
Let $v_1,\ldots,v_6$ be the vertices of $\Gamma (\T)$ which are endpoints 
of $e_1,\ldots,e_6$.
\nonumfig{1.0}{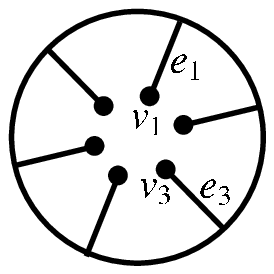}
\end{defn}

\begin{lemma}\label{lem:4.2}
Assume $\T$ is not split. 
If $v_i=v_j$ for some $i\ne j$, then the crossing
number of its projection can be reduced by free isotopy.
\end{lemma}

\begin{proof}
Assume $v_i = v_j$ with $i\ne j$. 
If $e_i,e_j$ are opposing at $v_i$, then $e_i\cup e_j$  is a string of 
$\T$ intersecting the other strings exactly once. 
By Lemma~\ref{lem:4.1}, $\T$ is split. 
So we assume $e_i,e_j$ are not opposing at $v_i$.
\nonumfig{1.0}{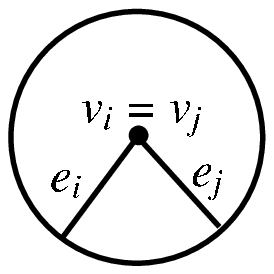}

\noindent
But then we can untwist $e_i,e_j$ to reduce the crossing number.  
\end{proof}

\begin{defn}
If $\T$ is not split, 
let $f_i$ be the face of $\Gamma (\T)$ containing $e_i,e_{i+1}$.
\end{defn}

\begin{lemma}\label{lem:4.3}
Assume $\T$ is not split. 
No two edges of $f_i$ correspond to the same edge of $\Gamma (\T)$. 
If two vertices of $f_i$ correspond to the same vertex of $\Gamma (\T)$, 
then the crossing number of the projection can be reduced by an isotopy 
fixed on the boundary. 
Finally, if $v_j$ is incident to $f_i$, then $j \in \{i,i+1\}$ or a 
crossing can be reduced by a free isotopy.
\end{lemma}

\begin{proof}
Assume that {two} edges of $f_i$ correspond to the same edge of $\Gamma (\T)$ {(Figure~\LemfigA)}.
Then there would be a circle in the interior of the tangle circle 
intersecting the edge once. 
{Then a string of $\T$ intersects this circle exactly once contradicting the Jordan Curve Theorem.}
\doublefig{1.0}{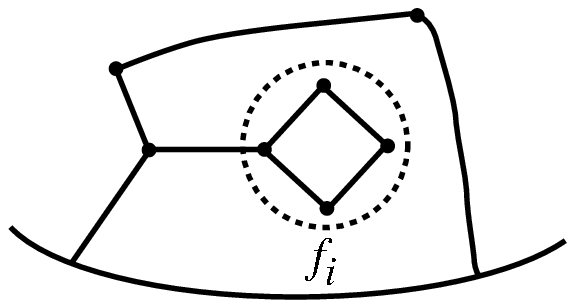}{\LemfigA}{1.0}{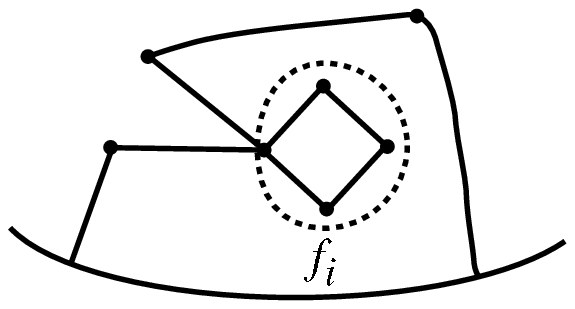}{\LemfigB}

Similarly, if two vertices of $f_i$ correspond to the same vertex, $v$, 
of $\Gamma (\T)$, there would be a circle intersecting $\Gamma (\T)$ 
only in $v$ {(Figure~\LemfigB)}. 
This would give rise to a crossing that could be reduced by an isotopy 
rel $\partial B^3$.

Now assume $v_j$ is incident to $f_i$ with $j\ne i$, $i+1$ { (Figure~\LemfigC)}. 
Then $e_j\notin f_i$ for otherwise $\T$ would be split. 
Thus $e_j$ lies in the exterior of $f_i$. 
\sdoublefig{1.4}{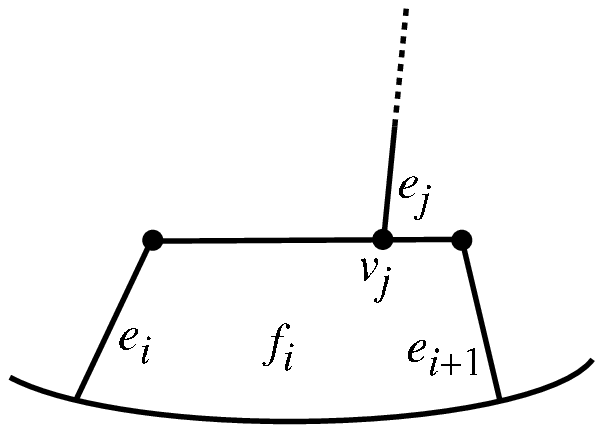}{\LemfigC}{.9}{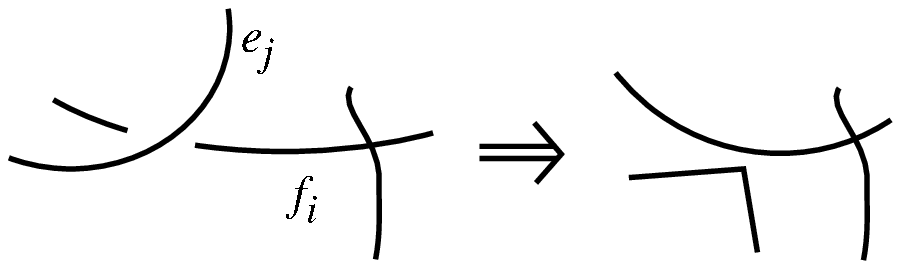}{\LemfigD}
But now we can untwist the crossing at $v_j$ by bringing $e_j$ into 
$f_i$ { (Figure~\LemfigD)}.
\end{proof}

\begin{thm}\label{lem:4.4}
If $\T$ is a 3-string tangle which can be freely isotoped
to a projection with at most seven crossings, then 
either $\T$ is split or has two parallel strands.
\end{thm}

\begin{proof}
Note that $\T$ being split or having parallel strands 
is invariant under free isotopy. 
Freely isotope $\T$ to a minimal projection with $\le 7$ crossings and 
let $\Gamma (\T)$ be the corresponding graph. { Assume $\T$ is not split.}
By Lemma~\ref{lem:4.2}, $\Gamma (\T)$ has at least 6 vertices (thus 
$\T$ has at least 6 crossings). 
If the six $v_i$ are the only vertices incident to $\bigcup_i f_i$, then 
(Lemma~\ref{lem:4.3}) we must have Figure~\lifesaverA.
\doublefig{1.35}{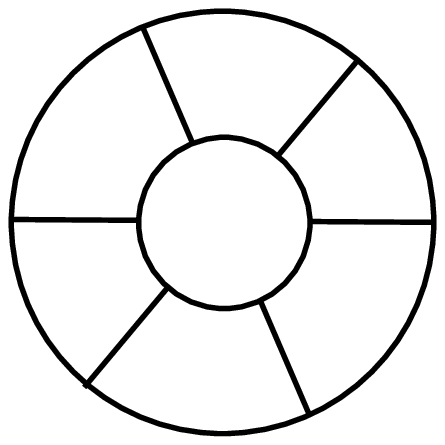}{\lifesaverA}{1.35}{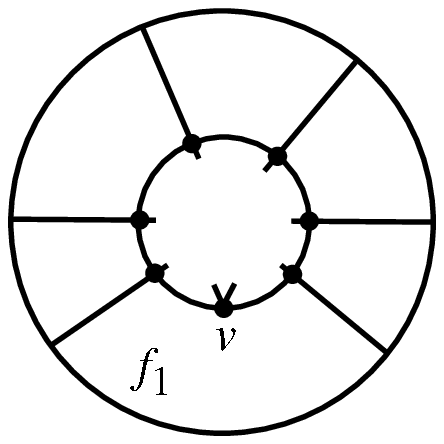}{\lifesaverB}

\noindent 
But then $\T$ has a closed curve, a contradiction.  

Thus $\Gamma (\T)$ must have another vertex $v$ lying on $f_1$, say, 
which is not a $v_i$. 
Since $\T$ has at most 7 crossings, there is exactly one such $v$ and we 
have Figure~\lifesaverB.
Enumerating the possibilities, {we see there are nine cases, six of which are shown in Figure~\hockeyballs, while the 
remaining three can be obtained from the top three via reflection}: 
\myfig{2.0}{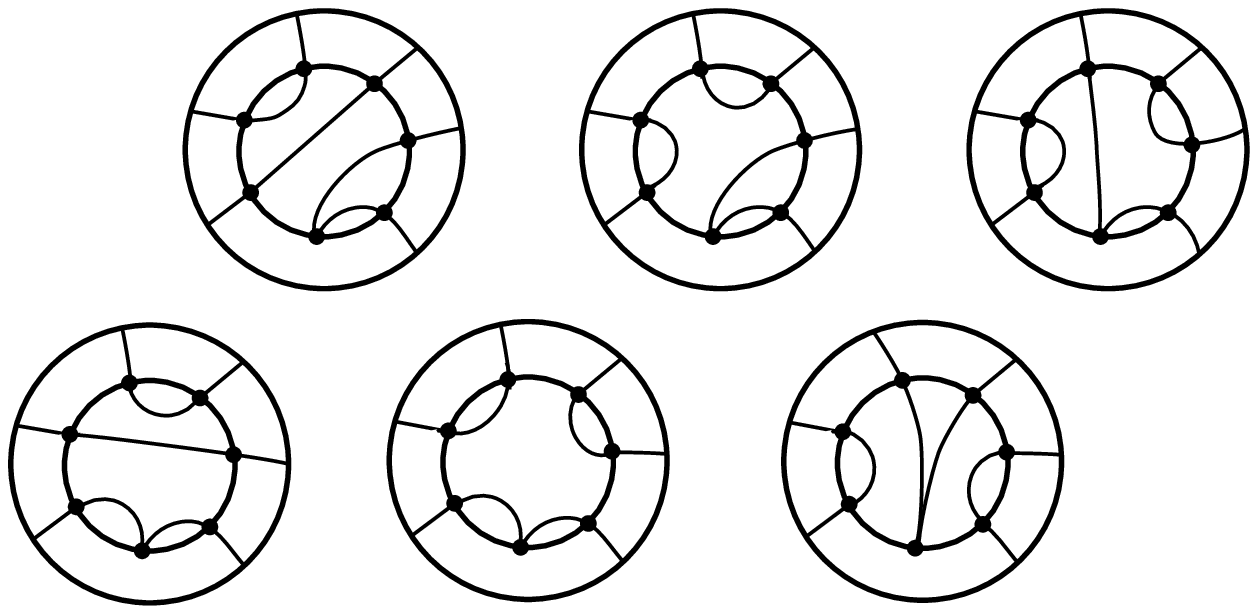}{\hockeyballs}

{Note  in each case we see two parallel strands (since we assumed $\T$ is not split).  For example, see 
Figure~\parallelstrandsC.}
\myfig{1.75}{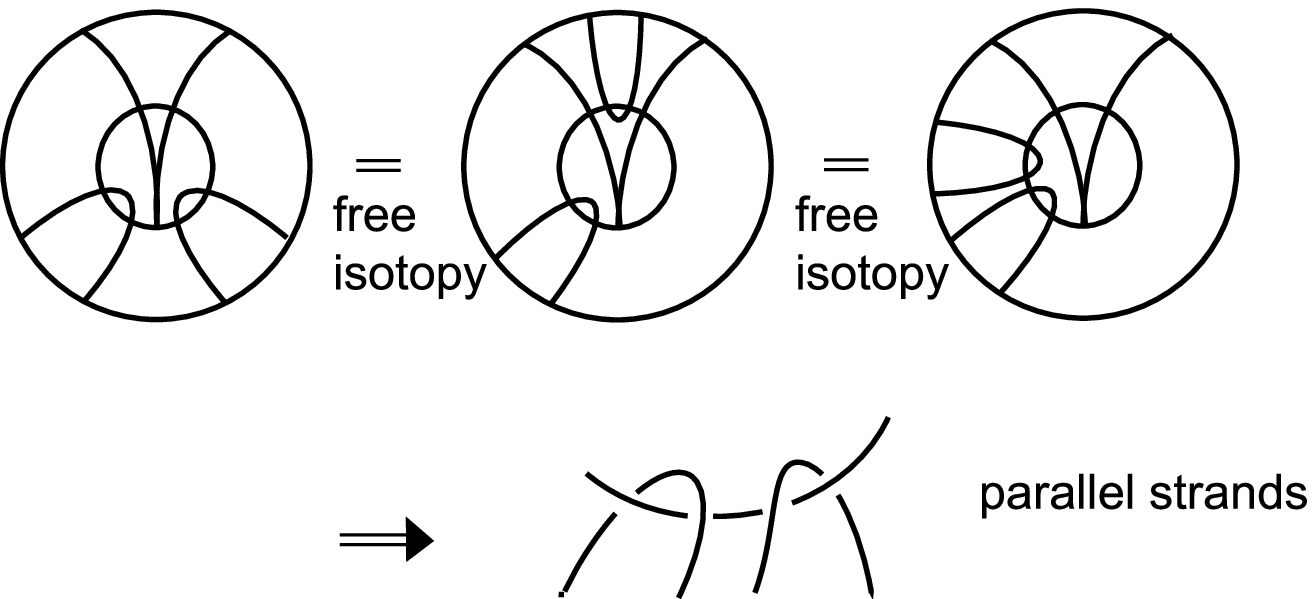}{\parallelstrandsC}

\end{proof}

Theorem \ref{lem:4.4} and Corollary \ref{corPJH} imply the following:

\begin{cor}\label{cor:PJH7}
 Assume $\T$ is a solution tangle.
If 
$\T$ can be freely isotoped to a projection with at most seven crossings, then $\T$ is 
the PJH tangle. \end{cor}

%% file: section5.tex
\section{If $\T$ is not the PJH tangle, then $\T$ has at least 10 crossings}	%% 5

Because the three DNA segments captured by the
 Mu transpososome in \cite{PJH} are relatively short (50, 175 and 190 base pairs), one can argue that 
any projection of the transpososome must
have few crossings.
Let $\T$ be an {\em in trans} solution tangle -- a tangle satisfying both the {\em in cis} and the {\em 
in trans} deletion experiments. 
In this section we show that if $\T$ is not the PJH
tangle (Figure \tangcircA), then $\T$ must be
complicated, as measured by its minimal crossing number:

\begin{prop}\label{prop:5.1}
Let $\T$ be an {in trans} solution tangle. 
If $\T$ has a projection with fewer than 10 crossings, then $\T$ 
is the PJH  tangle.
\end{prop}

In this section, we are interested in crossing number up to tangle equivalence, that is, up to isotopy of 
the tangle fixed 
on the boundary. So we must fix a framing for our {\em in trans} solution tangle. We take that of Proposition 
\ref{prop:synapsetangle}, the normal form.

\begin{remark}
 As pointed out in section 2, one goes from an {\em in trans} solution tangle in normal framing to one in the 
[PJH] framing by adding one left-handed twist at $c_2$. Thus, it formally follows from Proposition 
\ref{prop:5.1} that if an {\em in trans} solution tangle under the PJH framing has fewer than 9 crossings, it 
must be the [PJH] solution of Figure \PJHsolutionUnmarked.
 \end{remark}

The proof of Proposition \ref{prop:5.1} breaks down into checking cases of possible tangle 
diagrams. We focus on the fact that $s_{12} \cup s_{13} = \T - s_{23}$ is the $-1/2$ tangle. 
 Lemma \ref{lem:5.2} shows that each pair of strings must cross at least twice.  
Hence if an {\em in trans} solution tangle has less than 10 crossings, $s_{12} \cup s_{13}$ 
contains at most 5 crossings.  Theorems 5.5, 5.7, 5.12, 5.14 handle the cases when $s_{12} \cup 
s_{13}$ = 2, 3, 4, 5, respectively.  In each case it is shown that an {\em in trans} solution
 tangle with 
less than 10 crossings  can be isotoped to a tangle with  7 crossings.  By Corollary 
\ref{cor:PJH7}, such a tangle is the PJH tangle.
This 
section holds more generally under normal framing when the {\em in trans} deletion product is the (2,2) torus 
link while the 
remaining deletion products are $(2,L_i)$ torus links where $L_i \geq 4$.

\begin{lemma}\label{lem:5.2}
Let $\T$ be a solution tangle. 
In any projection of $\T$, $|s_{ij} \cap s_{ik}| \ge 2$.
\end{lemma}

\begin{proof}
Assume not. 
Then $|s_{ij}\cap s_{ik}| =0$. 
But then $\ell k (s_{ij} \cup x_{ij},s_{ik}\cup x_{ik}) = 0$,
contradicting Lemma~\ref{lem:synapsetangle}.
\end{proof}

Let $\T$ be an {\em in trans} solution tangle with crossing number less than 10. 
We assume that $s_{23}$ is the enhancer strand. 
To simplify notation we set $e= s_{23}$, $\alpha_2 = s_{12}$, 
$\alpha_3 = s_{13}$. 
Recall that $\T-e = \alpha_2 \cup \alpha_3$ is the $-1/2$-tangle.

\begin{assump}\label{assump:5.3}
$\T$ has a projection, in normal form (see the beginning of Section~2), 
with at most nine crossings. 
Furthermore, among all such projections, take $\alpha_2\cup \alpha_3$ 
to have the fewest crossings. 
\end{assump}

\begin{remark}
We will often blur the distinction between $e$, $\alpha_2$, $\alpha_3$ and 
their projections. 
\end{remark}

We will first prove in Theorem \ref{thm:5.4} that Proposition~\ref{prop:5.1} holds when $\alpha_2\cup \alpha_3$ has 
exactly two crossings.
In this case let $R_1,R_2,R_3$ be the closures of the complementary regions of 
$\alpha_2\cup \alpha_3$ in Figure~\Setup2.  {
Note all figures in this section have been rotated by 90 degrees, and thus  the tangle $\T-e = -1/2$ is displayed as an 
integral 2 tangle.} 

\myfig{1.5}{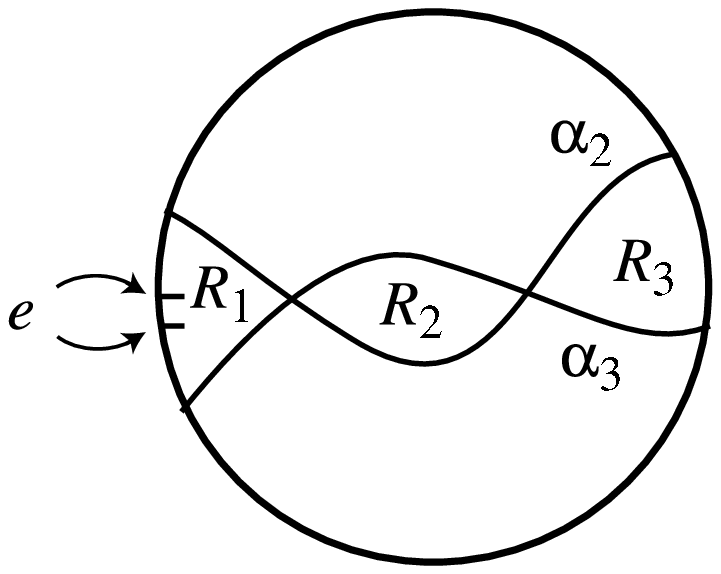}{\Setup2}

\begin{lemma}\label{lem:5.5}
If $e$ crosses $R_1\cap (\alpha_2\cup \alpha_3)$ fewer than four times then 
$\T$ can be freely isotoped to eliminate two crossings. 
Furthermore if $e\cap R_1$ contains only four crossings of $\T$ then we 
can reduce $\T$ by two crossings unless $e\cap R_1$ is as pictured in 
Figure~\IIsetup~
(up to symmetry).

\myfig{1.1}{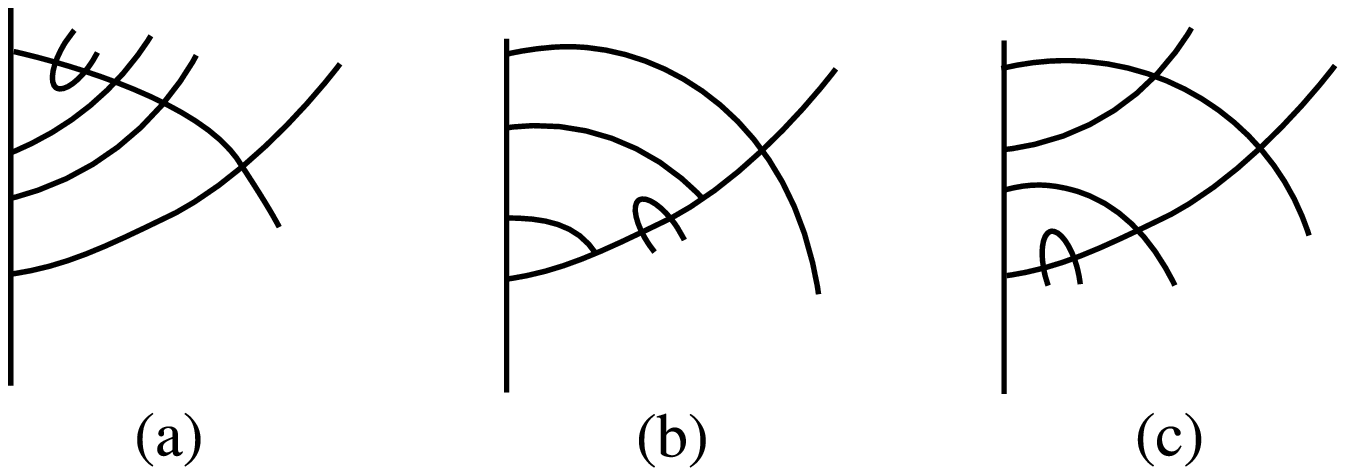}{\IIsetup}
\end{lemma}

\begin{proof}

Assume $e$ intersects $R_1\cap (\alpha_2\cup \alpha_3)$ in fewer than 
four crossings. 
By Lemma~\ref{lem:5.2}, $e$ must cross $R_1\cap (\alpha_2\cup \alpha_3)$ 
exactly twice. 
Then $R_1$ writes $\T\cup c_1$ as a disk sum. 
As $\T\cup c_1$ is a rational 2-string tangle, $e\cap R_1$ must be an integral 
summand \cite{L, ES1}. 
But then we can reduce the crossing number by at least two under a free 
isotopy by removing the crossings in the integral tangle $e\cap R_1$ as well as the two crossings where $e$ intersects 
$R_1\cap (\alpha_2\cup \alpha_3)$.
See Figure~\Reduceone. 
\myfig{1.5}{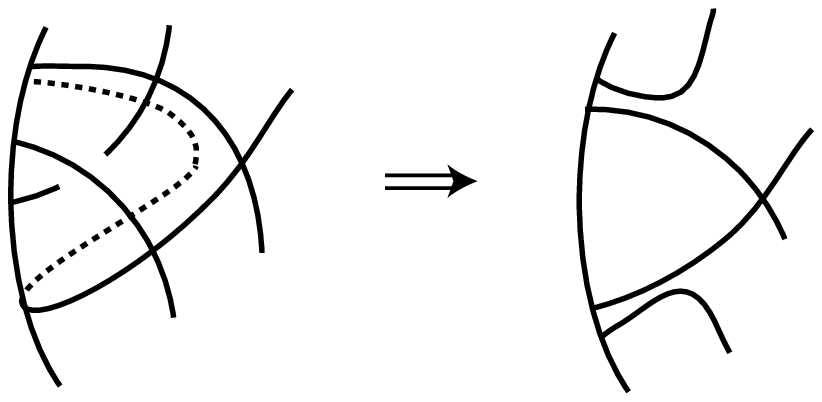}{\Reduceone}
Thus we assume $e\cap R_1$ must have at least 4 crossings in $R_1$ 
with $\alpha_2\cup \alpha_3$. 
Assume these are the only crossings of $\T$ in $R_1$. 
Then there are no crossings in $\text{int }R_1$. 
This allows for two crossing reductions except in the cases pictured 
in Figure~\IIsetup. 
See Figure~\Reducetwo.~\qed (Lemma~\ref{lem:5.5})
\myfig{1.0}{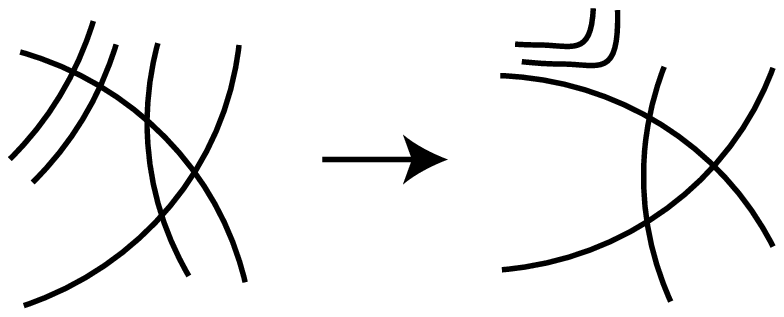}{\Reducetwo}
\renewcommand{\qed}{}
\end{proof}

\begin{thm}\label{thm:5.4} 
Let $\T$ be an {in trans} solution tangle satisfying Assumption~\ref{assump:5.3}. 
Under this projection, assume $\alpha_2\cup \alpha_3$ has exactly two 
crossings. 
Then $\T$ can be freely isotoped to have at most seven crossings. 
\end{thm}

\begin{proof}
Let $R_1,R_2,R_3$ be the closures of the complementary regions of 
$\alpha_2\cup \alpha_3$ in Figure~\Setup2.

We divide the proof of Theorem~\ref{thm:5.4} into two cases. 

\noindent {\bf Case I:}
$e$ has 4 intersections in $R_1$.

\noindent {\bf Case II:}
$e$ has at least $5$ intersections in $R_1$.

~

\noindent {\it Proof of Case I:}
$R_1$ contains exactly 4 crossings  (involving $e$). 
In this case, there are 3 possible configurations up to symmetry (shown in Figure~\IIsetup)
that do not immediately allow a crossing reduction by two in $R_1$. 

Since there are at most $9$ crossings, if $e$ intersects $R_2$ then 
$e\cap R_2$ must contribute exactly 2 crossings and $e\cap R_3$ must be empty. 
But this allows a reduction of two crossings:  In case (a) and (b) in Figure~\IIsetup, the two crossing involving 
$\alpha_2$ and $\alpha_3$ can be removed as shown in Figure~\IIredone.
\myfig{.65}{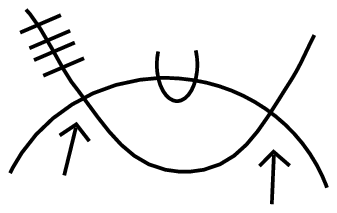}{\IIredone} 
In case (c) one reduction comes from $\alpha_1\cap \alpha_2$, the other from $e\cap (\alpha_2 \cup \alpha_3)$ 
(Figure \figure105).

\myfig{1.1}{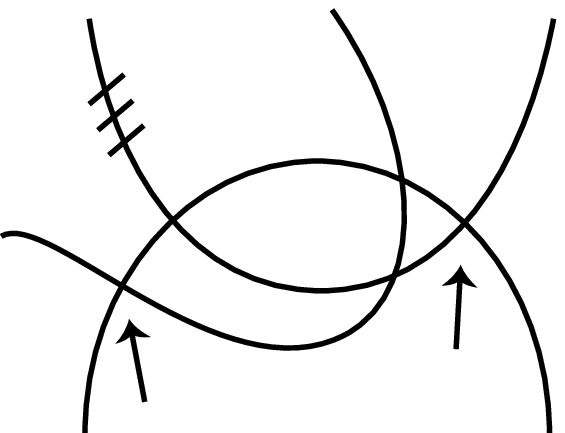}{\figure105}

Thus we take $e$ disjoint from $R_2$, and we must have (a) or (b) in 
Figure~\IIredtwo.
\myfig{1}{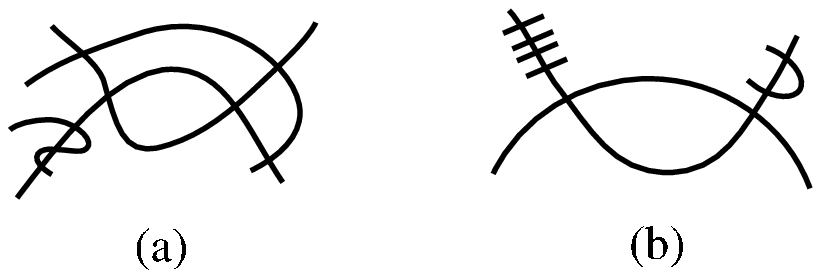}{\IIredtwo}
But (a) gives two crossing reductions and (b) contradicts Lemma~\ref{lem:5.2}.
\QED{Case I}

~

\noindent {\it Proof of Case II:}
Because $e$ must pick up at least two crossings in $R_2\cup R_3$, 
(else we can reduce by two crossings), $e$ must have exactly 4 crossings with 
$(\alpha_2 \cup \alpha_3)\cap R_1$ and exactly 1 self-crossing inside $R_1$.
This accounts for nine crossings. 
For subcases 1 and 2, we assume the two crossings are with $R_3$  --- if not, a similar argument 
works if the two crossings are with $R_2$ by using the fact  that the crossing at $R_3$ would be reducible. 

Let $\delta_2,\delta_3$ (in $\alpha_2,\alpha_3$, resp.) be the arc 
components of $[\alpha_2\cup \alpha_3 - (\alpha_2\cap \alpha_3)]\cap R_1$. 

\noindent {\bf Subcase 1).}
$|\delta_2 \cap e| = 2 = |\delta_3 \cap e|$ (intersection here refers 
to the projection) 

Then we may assume we have 
\myfig{1.0}{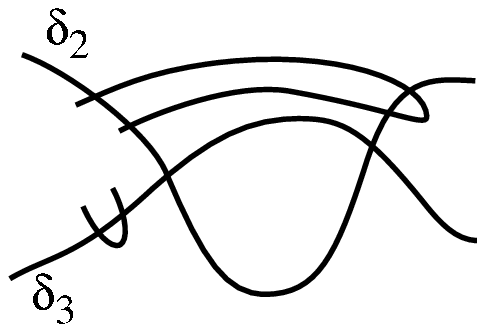}{\Ioneset}
Since $\text{int }R_1$ has exactly one crossing, the possible cases are 
\myfig{1.0}{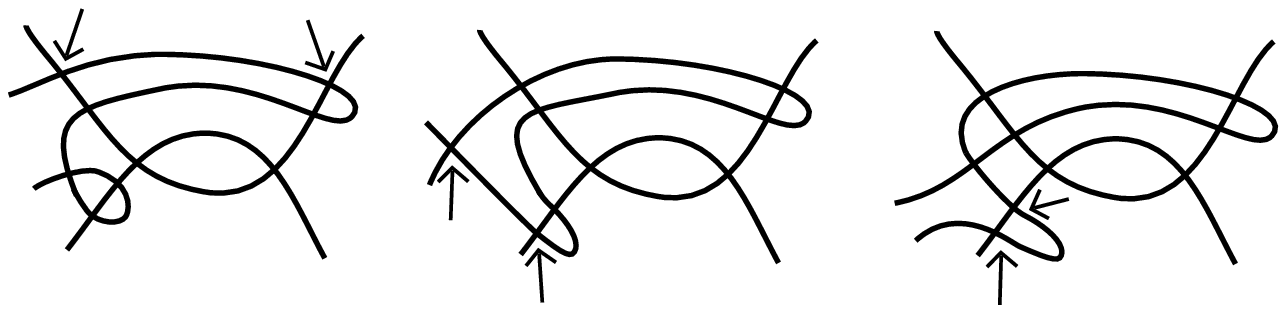}{\Ionered}
In each case we can reduce two crossings.

\noindent {\bf Subcase 2).}
$|\delta_2\cap e| = 3$, $|\delta_3 \cap e| =1$ (case when $|\delta_2\cap e| = 1$, $|\delta_3 \cap e| =3$ is similar).

{We have two possibilities
\myfig{1.2}{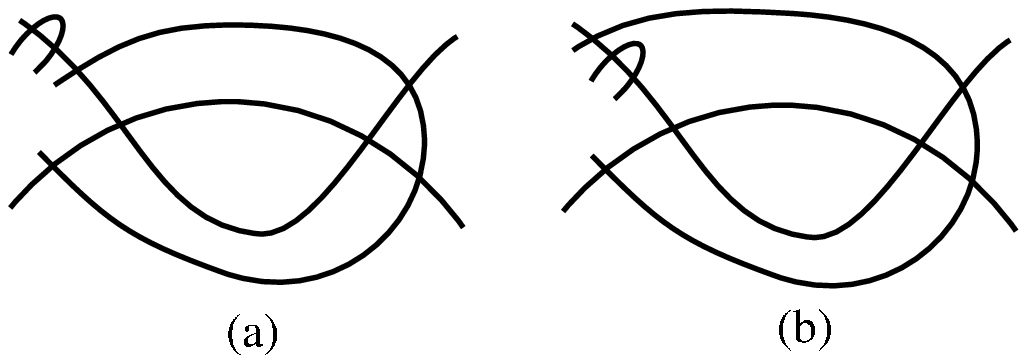}{\Itwosetup}

Recall that $R_1$ can have only one more crossing in its interior.  If this crossing can be reduced, we can reduce two 
more crossings by case 1.  Hence, we assume that this crossing cannot be reduced.    
Thus 
we see that both of the cases in Figure \Itwosetup~ allow a reduction of two crossings:
\myfig{1.6}{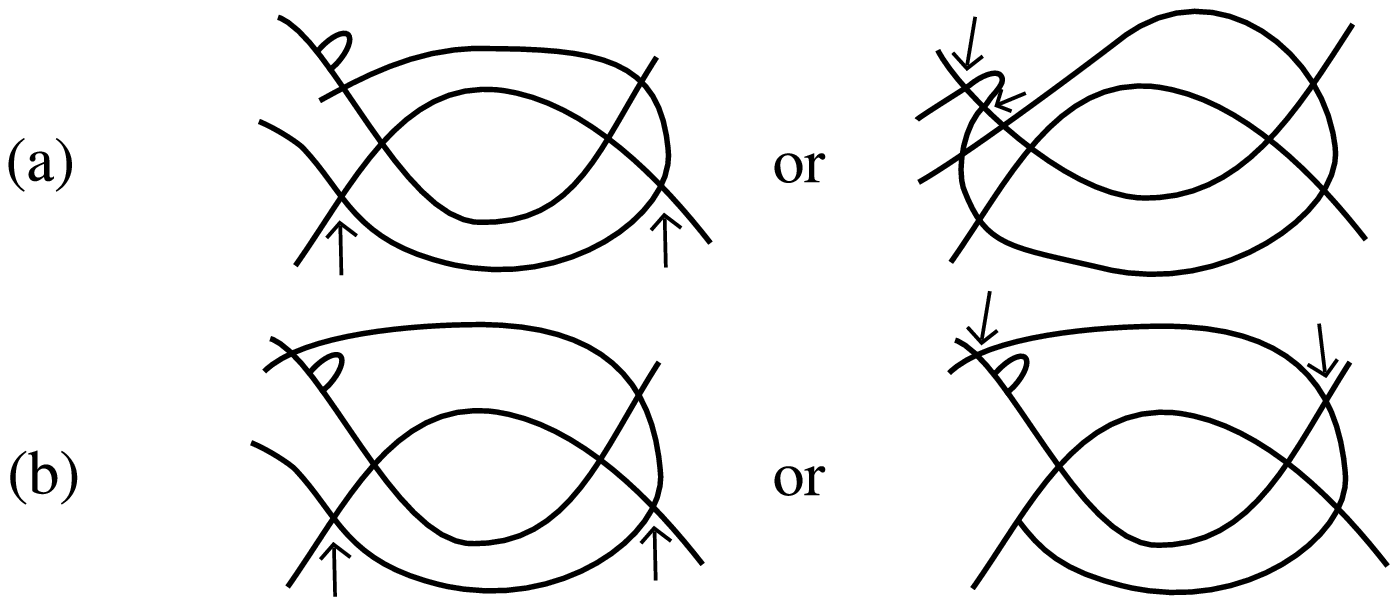}{\Itworeduce}}

\noindent {\bf Subcase 3).} 
$|\delta_2 \cap e| =4$, $|\delta_3\cap e| =0$ (case when $|\delta_2\cap e| = 0$, $|\delta_3 \cap e| =4$ is similar).

Recall $|e \cap (R_2 \cup R_3)| = 2$.  If $|e \cap R_3| = 2$, then $|e\cap \alpha_3|<2$, contradicting 
Lemma~\ref{lem:5.2}. If $|e \cap R_2| = 2$, then
we can remove the two crossings of $\alpha_2 \cap \alpha_3$.

Subcases (1)--(3) exhaust all possibilities, showing a reduction 
in crossing number by two in Case~II.
\QED{Case II}

This finishes the proof to Theorem~\ref{thm:5.4}.\qed (Theorem~\ref{thm:5.4})
\renewcommand{\qed}{}
\end{proof}

\begin{defn}
$\alpha_i$ has a trivial self-intersection w.r.t. $\alpha_j$ if $\alpha_i$ 
has a self-intersection such that the subarc of $\alpha_i$ connecting 
the double points is disjoint from $\alpha_j$, $\{i,j\} = \{2,3\}$.
\end{defn}

\begin{lemma}\label{lem:5.6}
If $\T$ is an {in trans} solution tangle satisfying Assumption~\ref{assump:5.3}
and $\alpha_i$ has a trivial self-intersection w.r.t. $\alpha_j$ for 
$\{i,j\} = \{2,3\}$, then we can reduce $\T$ by two crossings.
\end{lemma}

\begin{proof} 
Assume $\alpha_2$ has a trivial self-intersection w.r.t. $\alpha_3$, 
and let $\delta \subset \alpha_2$ be the subarc connecting its double points.
As $|\alpha_2\cap \alpha_3|$ is minimal, Figure~\Untwist\  shows that 
we may assume $|\delta \cap e| \ge4$.
\myfig{.65}{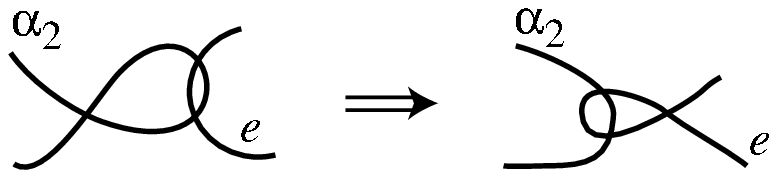}{\Untwist}

{By Lemma~\ref{lem:5.2}, $|\alpha_2\cap \alpha_3| \ge2$ and 
$|e\cap \alpha_3| \ge2$. 
Thus $|\delta \cap e| =4$, $|\alpha_2\cap \alpha_3| =2$, 
$|e\cap \alpha_3| =2$ and $\alpha_2$ has a self crossing, accounting for all nine crossings. }
Let $\alpha'_2 =\alpha_2 - \text{int}(\delta)$.
Let $R_1,R_2,R_3$ be the closures of the complementary regions of 
$\alpha'_2 \cup \alpha_3$ as in Figure~\Setup2. 
Then $|e\cap \partial R_1| \ge2$ and $|e\cap (\partial R_2\cup \partial R_3)|
\ge 2$, otherwise we can reduce by two the number of crossings. 
Thus $|e\cap (\alpha_2\cup \alpha_3)| \ge 8$. 
Because $|\alpha_2 \cap \alpha_3|\ge2$, we have too many crossings.
\end{proof}

\begin{thm}\label{thm:5.7}
Let $\T$ be an {in trans} solution tangle satisfying Assumption~\ref{assump:5.3}. 
If $\alpha_2\cup \alpha_3$ contributes three crossings to $\T$
(including self-crossings), then 
$\T$ can be freely isotoped to have at most seven crossings.
\end{thm}

\hskip -13pt {\it Proof.}
Let $\T$ be such an {\em in trans} solution tangle. 
Then $\alpha_3$, say, must have a self-intersection which we may assume 
is not trivial w.r.t. $\alpha_2$. 
Thus $\alpha_2\cup \alpha_3$ must be as in Figure~\Setupthree, 
with $R_1,R_2$ complementary components of $\alpha_2\cup \alpha_3$, 
$R_1$ containing the endpoints of $e$, 
and $\delta_1,\delta_2,\ldots,\delta_8$ the arc components of 
$\alpha_{{2}}\cup \alpha_3 - (\alpha_2\cap \alpha_3)$.
\myfig{1.5}{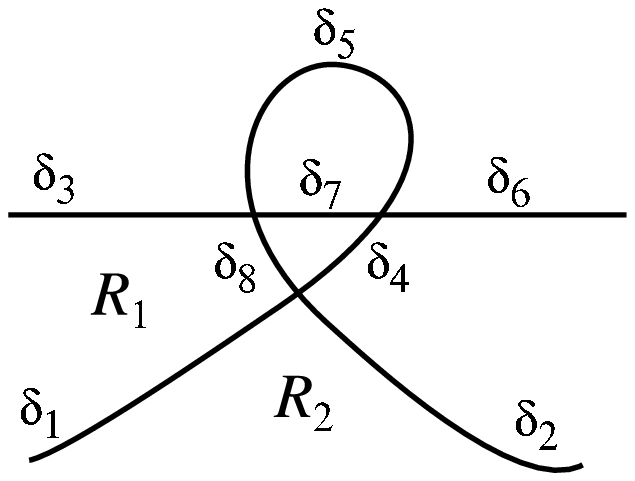}{\Setupthree}

\begin{claim}\label{clm:5.8}
If $|e\cap (R_1\cap (\alpha_2\cup \alpha_3))|<4$, then $\T$ can be 
reduced by two crossings. 
\end{claim}

\begin{proof} 
We may assume that $|e\cap (R_1\cap (\alpha_2\cup \alpha_3))|=2$. 
Then a copy of $\partial R_1$ writes $\T\cup c_1$ as a disk sum. 
Since $\T\cup c_1$ is {a rational tangle}, one of the summands must be 
integral w.r.t.\ to the disk. 
This must be the left-hand side, $R_1\cap e$. 
After a free isotopy we may take this to be $\frac01$. 
{Hence we can take $e$ to have no self-crossings in $R_1$.}

Let $\delta_1,\delta_2$ be as pictured in Figure~\Setupthree.
If $|e\cap \delta_1|=2$, then we can freely isotop away two crossings. 
If $|e\cap \delta_1| =1$, then $|e\cap \delta_2|\ge 3$ (else we can reduce 
two crossings in $R_2)$. 
Since $|e\cap \alpha_2| \ge2$ by Lemma~\ref{lem:5.2}, $|e\cap \alpha_2|$ 
must be two, accounting for all intersections. 
That is, if $|e\cap \delta_1| =1$ we are as in Figure~\threeredone, where 
we can reduce by two crossings.
\myfig{1.25}{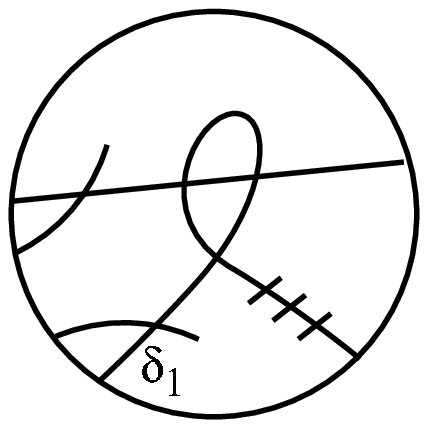}{\threeredone}
Thus we assume {\em $e$ is disjoint from $\delta_1$.}

Similarly we show that $e$ is disjoint from $\delta_3$ 
(see Figure~\Setupthree). 
If $|e\cap \delta_3|=2$, then we can reduce by two crossings in $R_1$. 
So assume $|e\cap \delta_3| =1$. 
If $ {e} \cap \delta_2$ is empty we can reduce two crossings 
(the self-intersection and $e\cap {\delta_3}$ ), so $|e\cap \delta_2|\ge2$. 
The only possibility is shown now in Figure~\threeredtwo, which we can 
reduce by two crossings. 
\wdoublefig{1.2}{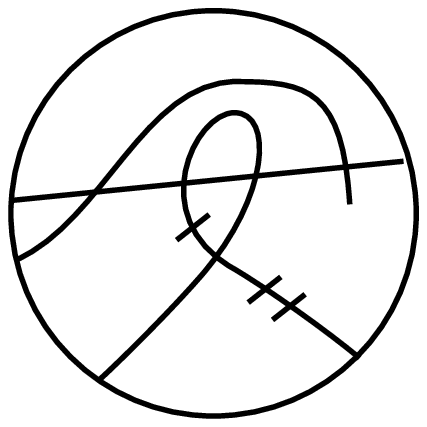}{\threeredtwo}{1.2}{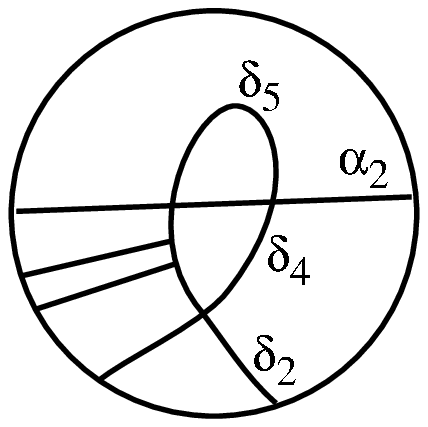}{\threeredthree}

So we assume $e$ is disjoint from arcs $\delta_1$ and {$\delta_3$}  in 
$\partial R_1$. 
See Figure~\threeredthree.

Now $|e\cap \alpha_2| \ge2$ by Lemma~\ref{lem:5.2}, and $|e\cap\delta_2|\ge2$
(else we can reduce two crossings). 
This accounts for all nine crossings. 
But we must also have a crossing between $e$ and $\delta_4\cup \delta_5$.
\end{proof}

\begin{claim}\label{subclm:5.10}
$e$ must intersect $\delta_1$.
\end{claim}

\begin{proof}
By Claim~\ref{clm:5.8}, we may assume  $|e\cap (R_1\cap (\alpha_1\cup 
\alpha_3))|\geq 4$. 
Label the arc components of $\alpha_2 \cup \alpha_3 - (\alpha_2\cap 
\alpha_3)$ as in Figure~\Setupthree.
Assume $e$ is disjoint from $\delta_1$. 
It cannot also be disjoint from $\delta_2$, otherwise we can eliminate the 
crossing at $\delta_1\cap \delta_2$ and argue as in Case~I of 
Theorem~\ref{thm:5.4}. 
Thus $e$ must intersect $\delta_2$ exactly twice, thereby accounting for all 
crossings. 
But the crossings of $e$ at $\delta_2$ lead to additional crossings.
\end{proof}

\begin{claim}\label{clm:5.9}
$\T$ can be reduced by two crossings if 
$|e\cap (R_1 \cap (\alpha_2 \cup \alpha_3))|<6$.
\end{claim}

\begin{proof} 
By Claim~\ref{clm:5.8}, we may assume  $|e\cap (R_1\cap (\alpha_1\cup 
\alpha_3))|=4$. 

\begin{subclaim}\label{subclm:5.11} 
$e\cap (\delta_6\cup \delta_7)$ is non-empty.
\end{subclaim}

\begin{proof} 
{Suppose $e\cap (\delta_6\cup \delta_7)$ is empty. Then}   $|e\cap \delta_3| =2$ by Lemma~\ref{lem:5.2}. 
If $e$ crosses $\delta_5$ then it must cross twice, accounting 
for all crossings of $\T$. 
Then there are no crossings in $\text{int }R_1$. 
Since $|e\cap \delta_1| >0$ and $|e\cap \delta_3| =2$, we have two 
crossing reductions in $R_1$.
Thus $e$ is disjoint from $\delta_5$. 

Similarly $e$ must not cross $\delta_2$. 
If $|e\cap \delta_2|\ne0$, then there can be no crossings in 
$\text{int }R_1$. 
Since $|e\cap \delta_1| >0$ and $|e\cap \delta_3|=2$, 
we see the crossing reductions for $\T$. {Since $e$ does not cross $\delta_2$, but does cross $\delta_1$, $e$ must 
cross $\delta_1$ twice.  Hence $e$ also does not cross $\delta_8$.}

Thus we assume $e$ does not cross $\delta_2$, $\delta_5$, 
$\delta_6$, $\delta_7$, {or $\delta_8$}. 
As the 2-string tangle $\alpha_2,\alpha_3$ forms a two crossing tangle, 
we see that the crossings between $\alpha_2$ and $\alpha_3$ may be 
reduced in $\T$.\qed (\ref{subclm:5.11}) 
\renewcommand{\qed}{}
\end{proof}

Assume first that $e$ crosses $\delta_7$. 
Then $e$ must cross $\delta_5\cup \delta_7$ twice, accounting for 
all crossings. 
Then $e$ does not cross $\delta_2$ and crosses $\delta_1$ exactly 
two times. 
These crossings can be reduced. 

Thus $e$ crosses $\delta_6$. 
If $e$ does not cross $\delta_2$ then we again see two crossings at 
$\delta_1$ that can be reduced. 
If $e$ crosses $\delta_2$ then it crosses it once, accounting for all crossings.
Thus $e$ crosses each of $\delta_1$ and $\delta_3$ an odd number of times 
and has no self-crossings in $\text{int }R_1$. 
Enumerating the possibilities one sees that we can reduce the crossings 
of $e$ at $\delta_3, \delta_6$ or at $\delta_1,\delta_2$.

This finishes the proof of Claim~\ref{clm:5.9}.\qed (\ref{clm:5.9})

By Claim~\ref{clm:5.9} we may assume 
$|e\cap R_1\cap (\alpha_2\cup\alpha_3)|=6$.
This accounts for all crossings of $\T$. 
By Lemma~\ref{lem:5.2}, $|e\cap \delta_3|\ge2$. 
As all crossings are accounted for, two of the crossings of $e$ 
with ${\delta_1 \cup}\delta_3$ can be eliminated by a free isotopy.
\QED{Theorem~\ref{thm:5.7}}\end{proof}

\begin{thm}\label{thm:5.12} 
Let $\T$ be an {in trans} solution tangle satisfying assumption~\ref{assump:5.3}. 
If $\alpha_2\cup\alpha_3$ has four crossings (including self-crossings), 
then there is a free isotopy of $\T$ reducing it to at most seven crossings.
\end{thm}

\begin{proof} 
Assume $\T$ is as hypothesized. 
By Lemma~\ref{lem:5.6}, $\alpha_2\cup \alpha_3$ has no trivial loops and 
we enumerate the possibilities for $\alpha_2\cup \alpha_3$ in 
Figure~\Setupfour\ where $R_1$ is the closure of the (planar) complementary 
region of $\alpha_1\cup \alpha_2$ containing the ends of $e$ and where 
$\delta_1,\delta_2$ are the extremal arc components of $\alpha_2\cup \alpha_3
- (\alpha_2\cap \alpha_3)$ in $R_1$.
\myfig{2.25}{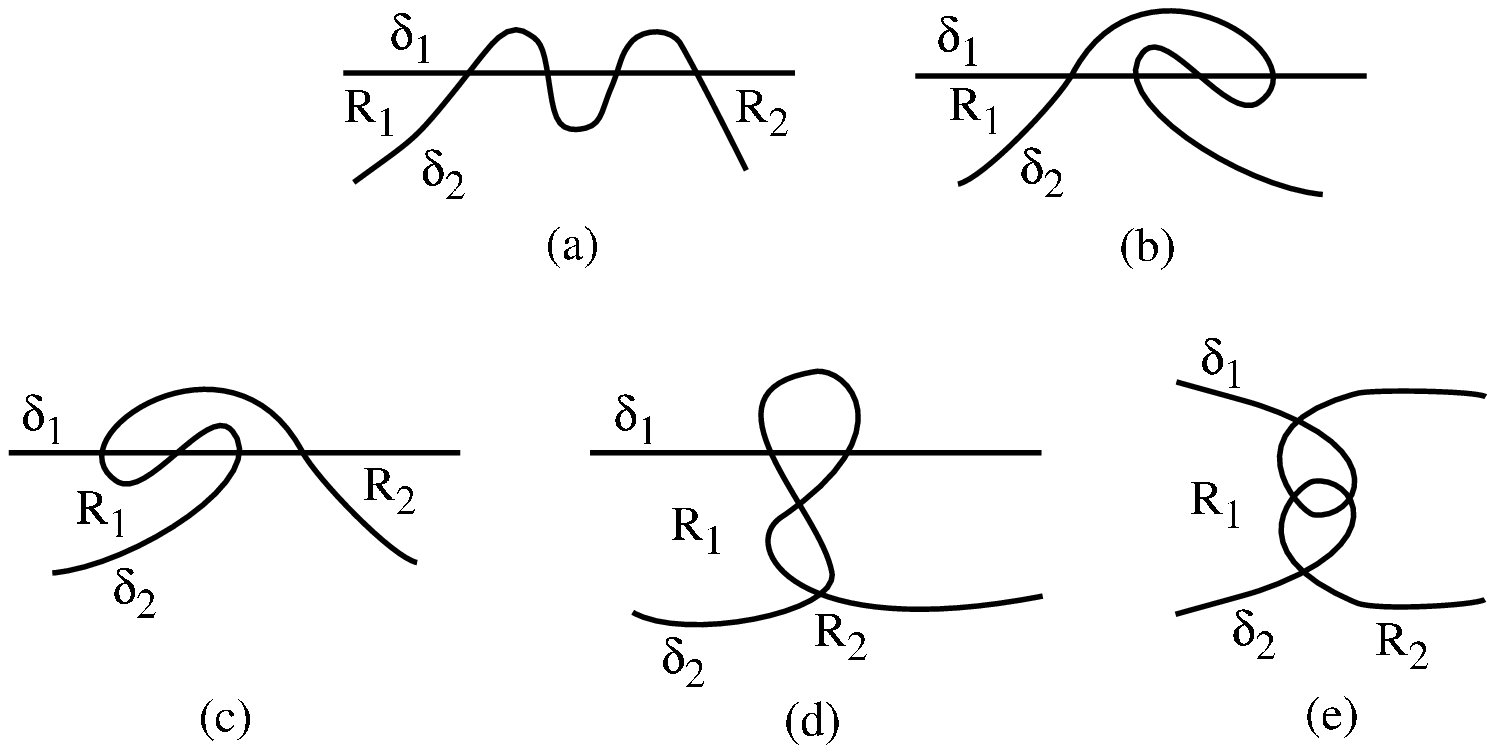}{\Setupfour}
\renewcommand{\qed}{}
\end{proof}

\begin{lemma}\label{lem:5.13}
With $\T$ as in the hypothesis of Theorem~\ref{thm:5.12}, if $e$ 
contains a self-crossing then that self-crossing appears in $R_1$.
\end{lemma}

\begin{proof} 
Assume a self-crossing appears in a complementary region $R\ne R_1$ of 
$\alpha_2\cup \alpha_3$. 
{Since $|\alpha_2\cap \alpha_3| = 4$, $|e \cap \alpha_2| = 2$ $|e \cap \alpha_3| = 2$, $e$ can have at most one 
self-crossing. Hence}
$R$ has exactly 5 crossings involving $e$: four from $e$ crossing 
$\alpha_2\cup \alpha_3$ and one self-crossing. 
By Lemma~\ref{lem:5.2}, $e$ must cross $\alpha_2$ twice and $\alpha_3$ 
twice in $R$.
Exactly one pair of these crossings must be in a component of 
$\alpha_2\cup\alpha_3 - (\alpha_2\cap\alpha_3)$ shared with $R_1$ 
(all crossings are accounted for). 
If this component is $\delta_1$ or $\delta_2$ we may reduce by two crossings.
So we assume this pair is in a different component.
Looking at the possibilities of Figure~\Setupfour, we immediately rule 
out (a) and (b). 
{ In cases (c), (d), and (e), we can reduce two crossings of $\alpha_2\cup\alpha_3$ since we have accounted for all 
crossings   (for case (d), note  $e$ must cross each of $\alpha_2$ and $\alpha_3$ twice in $R$).}

\end{proof}

We continue the proof of Theorem~\ref{thm:5.12}. 

\subsection*{Case I: $e$ crosses $R_1\cap (\alpha_2\cup\alpha_3)$ 
exactly twice}
If both crossings are in $\delta_1\cup\delta_2$, then we can reduce by 
two crossings. 
This rules out configurations (a),(b) of Figure~\Setupfour. 
If $\T$ has nine crossings, then one must be a self-crossing. 
By Lemma~\ref{lem:5.13} it must occur in $R_1$ and hence can be untwisted 
to reduce the crossing number by one. 
Thus we assume $\T$ has only 8 crossings. 
Hence $e$ must be disjoint from $\delta_1 \cup \delta_2$. 
In cases (c), (d), and (e), $e$ must then be disjoint from the 
region labelled $R_2$. 
But then we can reduce the crossing at $R_2$.\qed (Case I)

\subsection*{Case II: $e$ crosses $R_1\cap (\alpha_2\cup\alpha_3)$ four times}
If $\T$ has nine crossings then one must be a self-crossing of $e$. 
By Lemma~\ref{lem:5.13}, the self-crossing occurs in $R_1$. 
Thus, whether $\T$ has eight or nine crossings, all crossings of $\T$ 
involving $e$ lie in $R_1$. 
By Lemma~\ref{lem:5.2}, two of the crossings of $e$ are with $\alpha_2$ 
and two with $\alpha_3$. 
Looking at Figure~\Setupfour, we see that $e$ must have exactly two crossings 
with $\delta_i$ for $i=1$ or 2 (in cases (c) and (e)
we could otherwise reduce by two crossings). 
Either we can reduce by two crossings or near $\delta_i$ we have 
Figure~\Setupdl. 
\myfig{1.0}{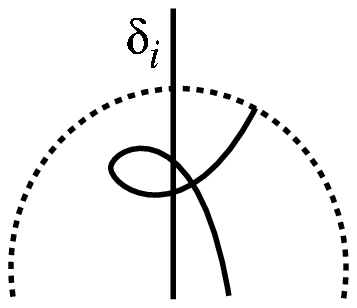}{\Setupdl}
This writes $\T\cup c_1$ as a disk sum. 
As $\T\cup c_1$ is integral w.r.t.\ the disk slope, each summand is integral 
w.r.t.\ to the disk slope. 
But this means we can eliminate the crossings at $\delta_i$ by a free 
isotopy of $\T$.\qed (Case II)
\QED{Theorem~\ref{thm:5.12}}

\begin{thm}\label{thm:5.14} 
Let $\T$ be an {in trans} solution tangle satisfying assumption~\ref{assump:5.3}. 
If $\alpha_2\cup\alpha_3$ has five crossings (including self-crossings), 
then there is a free isotopy of $\T$ reducing it to at most seven crossings.
\end{thm}

\begin{proof}
Let $\T$ be such a tangle. 
WLOG assume $\alpha_2$ has at most one self-intersection. 
By Lemma~\ref{lem:5.2} $e$ must cross each of $\alpha_2,\alpha_3$ exactly 
twice. 
Since $e$ can contribute at most 4 crossings this accounts for all crossings 
coming from $e$. 
In particular, $e$ has no self-crossings. 
\renewcommand{\qed}{}
\end{proof}

\subsection*{Case I} 
$\alpha_2$ crosses $\alpha_3$ twice.

We have three possibilities as shown in Figure~\SetupFT\ (Lemma~\ref{lem:5.6}),
where $R_1$ is the complementary region  of $\alpha_2\cup \alpha_3$ containing
the ends of $e$ {($\ell$ is discussed below)}. 
\myfig{1.5}{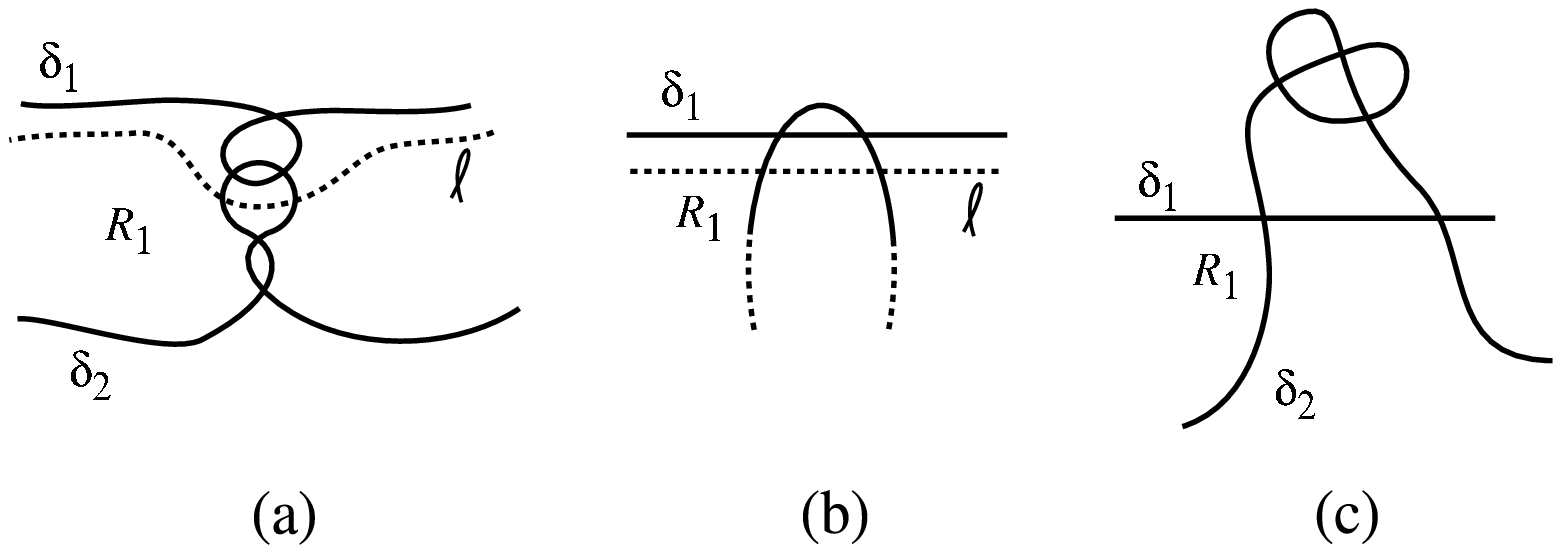}{\SetupFT}
In case (c) we may reduce the crossings of $\T$ by two ($\T$ has no local 
knots, hence $|e\cap (\delta_1\cup\delta_2)| =2$). 
So we restrict our attention to (a) and (b) and consider the arc $\ell$ 
in Figure~\SetupFT. {Since $|e \cap \alpha_2|= 2$},
$\ell$ intersects $e$ twice{.  Hence}  $e\cup\ell$ is as in Figure~\SetupFell, 
where $e',e''$ are components of $e-\ell$.
\myfig{1.0}{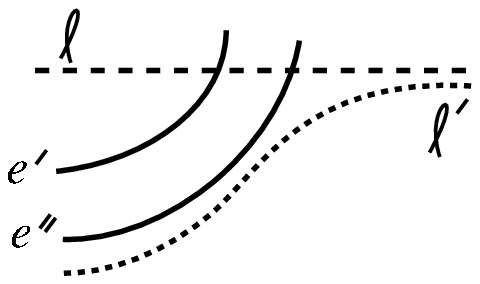}{\SetupFell}
Let $\ell'$ be the arc pictured  {in Figure~\SetupFell, consisting of $e''$ and part of $\ell$}. 

\begin{claim}\label{clm:5.15} 
$|\ell' \cap \alpha_3|=2$. 
\end{claim}

\begin{proof} 
$|\ell'\cap\alpha_3|$ is even and at least 2, since the endpoints of 
$\alpha_3$ are below $\ell'$. 
If $|\ell'\cap\alpha_3| \ge4$, then $|e'' \cap \alpha_3|=2=
|(\ell'-e'')\cap\alpha_3|$. 
Going back to Figure~\SetupFT\ we see we can eliminate the crossings of 
$e$ and  {$\alpha_2$}.\qed (\ref{clm:5.15})

By Claim~\ref{clm:5.15}, $\ell'$ writes $\T\cup c_2$ as a disk sum. 
By Lemma~2.1,  if $\T$ is a normal form solution tangle, $\T\cup c_2$ is the {-$1/4$} tangle.
As the crossing number of the summand below $\ell'$ is at most 3, 
the tangle below $\ell'$ is integral w.r.t.\ the disk slope --- 
allowing us to eliminate the crossings there.
Thus there is at most one crossing below $\ell'$. 
Now we can eliminate two crossings from $\T$, 
{where $\alpha_3$ crosses $e''$ or $\alpha_2$ near where $\alpha_3$ crosses $\ell'$}.
\qed (Case I)
\renewcommand{\qed}{}
\end{proof}

\subsection*{Case II}
$\alpha_2$ crosses $\alpha_3$ four times.

WLOG assume $\alpha_2$ has no self-intersections and $\alpha_3$ only one. 
Then $\alpha_2\cup\alpha_3$ must be one of the cases in 
Figure~\SetupFF.
\myfig{4.5}{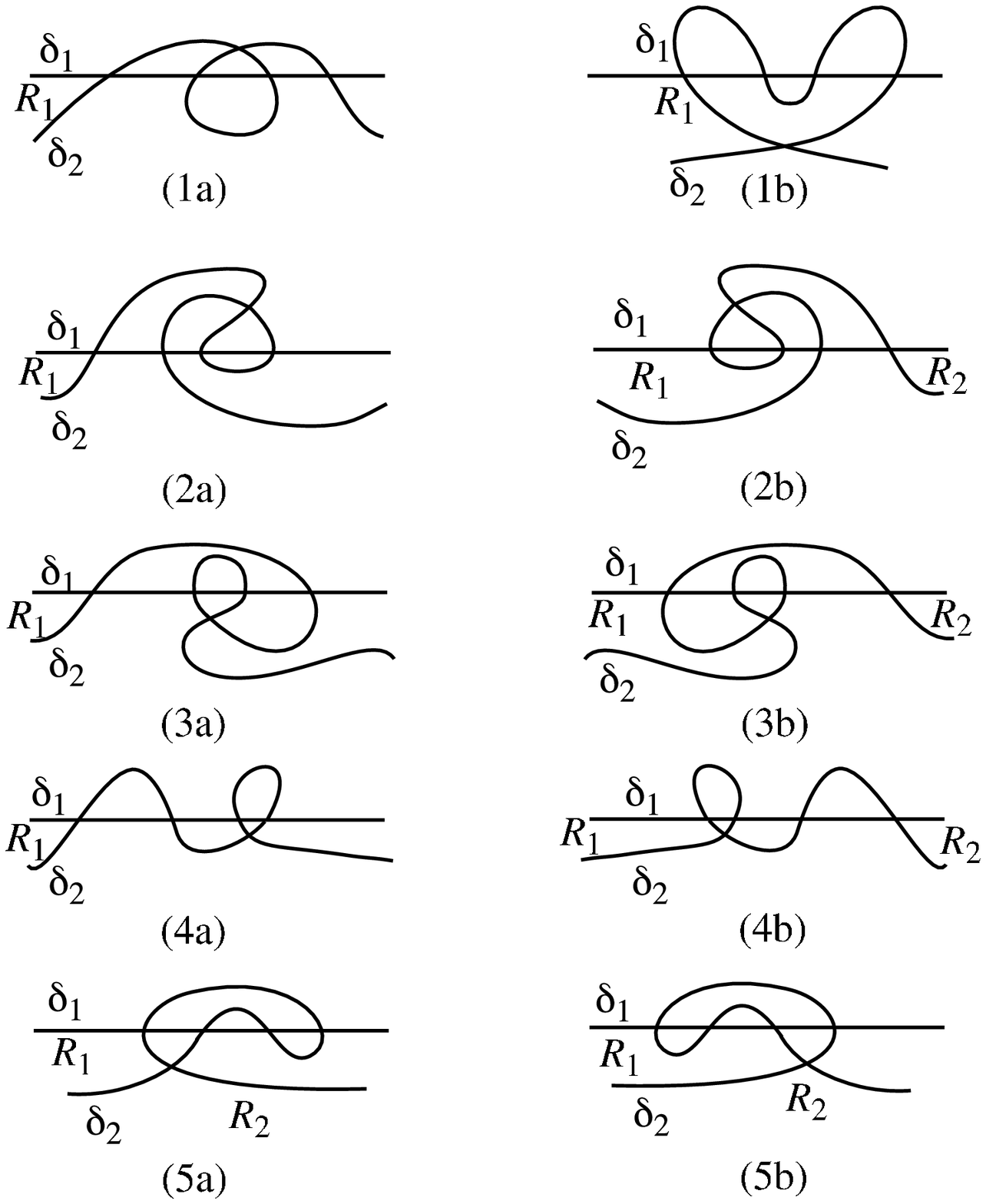}{\SetupFF}

\begin{claim}\label{clm:5.16}
$|e\cap (\delta_1\cup\delta_2)| \le1$.
\end{claim}

\begin{proof}
Assume $|e\cap (\delta_1\cup \delta_2)|\ge2$. 
If $|e\cap R_1\cap (\alpha_2\cup \alpha_3)|\le 2$, then the crossings at 
$\delta_1\cup\delta_2$ can be eliminated. 
So assume $|e\cap R_1\cap (\alpha_2\cup\alpha_3)|=4$, accounting for all 
crossings of $e$. 
Then $|e\cap \delta_i| =2$ for some $i=1,2$, and we can reduce 
these crossings.\qed (\ref{clm:5.16})
\renewcommand{\qed}{}
\end{proof}

Claim~\ref{clm:5.16} eliminates cases (1{a}), (2a), (3a), (4a). 
If $|e\cap\delta_i|=1$ then it must give rise to a reducible crossing in $R_1$.
In cases (2b), (3b), (4b) we could then reduce also the crossing at $R_2$. 
In these cases then we assume $e$ disjoint from $\delta_1\cup\delta_2$. 
By inspection we now see we can reduce by two crossings 
$(|e\cap (\alpha_2\cup\alpha_3)|=4)$.

We are left with {(1b),} (5a) and (5b).
If $|e\cap\delta_2|=1$, then we can eliminate two crossings at $e\cap R_2$. 
If $|e\cap\delta_1|=1$, then we can eliminate this crossing and the 
self-crossing of $\alpha_3$. 
With $|e\cap \alpha_2| =2=|e\cap\alpha_3|$ and $e\cap (\delta_1\cup\delta_2)$
empty, we see that in all possibilities we can eliminate two crossings 
from $\T$.\QED{Theorem \ref{thm:5.14}}

%% file: mu.bbl
\begin{thebibliography}{MMM}

\bibitem[BM]{BM} T Baker, K. Mizuuchi, 
{\em DNA-promoted assembly of the active tetramer of the Mu transposase},
 Genes Dev. {\bf 6(11)} 
(1992) 2221--32.

\bibitem[BSC]{BSC}
J. Bath, D.J. Sherratt and S.D. Colloms,
{\em Topology of Xer recombination on catenanes produced by lamda 
integrase}, J. Mol. Biol. {\bf 289(4)} 
(1999), 873--883.

\bibitem[BL]{BL}
S. Bleiler and R. A. Litherland, {\em Lens spaces and Dehn surgery}, 
Proceedings of the American Mathematical Society {\bf 107(4)} (1989), 
1127-1131.

\bibitem[BV]{BV}
D. {Buck}, C. {Verjovsky Marcotte}, {\em
Tangle solutions for a family of DNA-rearranging proteins}\/, 
Math. Proc. Camb. Phil. Soc. {\bf  139(1)}  (2005),  59--80.

\bibitem[C1]{C1}
H. {Cabrera-Ibarra}, {\em
On the classification of rational 3-tangles}\/,
Journal of Knot Theory and its Ramifications {\bf 12} (2003), 921--946.

\bibitem[C2]{C2}
H. {Cabrera-Ibarra}, {\em
Results on the classification of rational 3-tangles}\/,
Journal of Knot Theory and its Ramifications {\bf 13} (2004), 175--192.

\bibitem[CG]{CG}
A. J. Casson,  C. McA. Gordon {\em Reducing Heegaard splittings.}\/ Topology Appl. {\bf 27} (1987), no. 3, 275--283. 

\bibitem[CH]{CH} G. Chaconas and R. M. Harshey, {\em Transposition of phage Mu DNA}. In Mobile DNA II. N. L. Craig, R. 
Craigie, M. Gellert, and
A. M. Lambowitz (ed), (2002) pp 384--402, American Society for Microbiology.

\bibitem[CLW]{CLW}
G. Chaconas, B.D. Lavoie, and M.A. Watson,
{\em DNA transposition: jumping gene machine, some assembly required},
Curr. Biol. {\bf 6} (1996), 817--820.

\bibitem[CBS]{CBS}
S.D. Colloms, J. Bath and D.J. Sherratt, 
{\em Topological selectivity in Xer site-specific recombination}, 
Cell {\bf88} (1997), 855--864. 

\bibitem[CGLS]{CGLS}
M. Culler, C. Gordon,  J. Luecke, P.B. Shalen {\em Dehn surgery on knots.}\/,  Bull. Amer. Math. Soc. (N.S.) {\bf 13}  
(1985), no. 1, 43--45.

\bibitem[D]{D}
I. {Darcy}, {\em
Biological distances on DNA knots and links: Applications to Xer recombination}\/,
Journal of Knot Theory and its Ramifications {\bf 10} (2001), 269--294.

\bibitem[D1]{D1} 
I. K. Darcy, A. Bhutra, 
 J. Chang, 
 N. Druivenga, 
 C. McKinney,  R. K. Medikonduri,  S. Mills, 
 J. Navarra Madsen, 
 A. Ponnusamy, 
 J. Sweet, 
 T. Thompson, 
{\em Coloring the Mu Transpososome}\/,
 {BMC Bioinformatics},
 {\bf 7},
(2006),
{{Art. No. 435}}.


\bibitem[DMBK]{DMBK} P. L. Deininger, J. V. Moran, M. A. Batzer, H. H. Kazazian Jr., 
{\em Mobile elements and mammalian genome evolution},
Curr Opin Genet Dev. 2003 {\bf 13(6)} (2003), 651--8.

\bibitem[EE]{EE}
J. {Emert}, C. {Ernst}, {\em
$N$-string tangles}\/,
Journal of Knot Theory and its Ramifications {\bf 9} (2000), 987--1004.

\bibitem[ES1]{ES1}
C. {Ernst}, D. W. {Sumners}, {\em
A calculus for rational tangles: applications to DNA recombination}\/, 
Math. Proc. Camb. Phil. Soc. {\bf 108} (1990), 489--515.

\bibitem[ES2]{ES2}
C. {Ernst}, D. W. {Sumners}, {\em
Solving tangles equations arising in a DNA recombination model}\/, Math.
Proc. Camb.Phil. Soc. {\bf 126} (1999), 23-36.


\bibitem[GEB]{GEB}
I. Goldhaber-Gordon, M.H. Early, and T.A. Baker,
{\em MuA transposase separates DNA sequence recognition from catalysis},
Biochemistry {\bf 42} (2003), 14633-14642.

\bibitem[G]{G} 
C. Gordon {\em On the primitive sets of loops in the boundary of a 
		handlebody}, Topology and Appl. {\bf 27} (1987), 285--299.


\bibitem[GBJ]{GBJ}
I. Grange, D. Buck, and M. Jayaram,
{\em Geometry of site alignment during int family recombination: antiparallel 
synapsis by the $FLp$ recombinase},
J. Mol. Biol. {\bf298} (2000), 749--764.

\bibitem[GGD]{GGD}
F. Guo, D.N. Gopaul, and G.D. van Duyne,
{\em Structure of Cre recombinase complexed with DNA in a site-specific 
recombination synapse}, 
Nature {\bf 389} (1997), 40--46.

\bibitem[HJ]{HJ}
R. Harshey and M. Jayaram, 
{\em The mu transpososome through a topological lens},
Crit Rev Biochem Mol Biol. {\bf 41(6)} (2006), 387--405. 


\bibitem[HS]{HS}
M. Hirasawa and K. Shimakawa, 
{\em Dehn surgeries on strongly invertible knots which yield lens spaces}, 
PAMS {\bf 128} (2000), no. 11, 3445--3451.


\bibitem[J]{J} 
W. Jaco, 
{\em Adding a 2-handle to 3-manifolds: an application to Property $R$},
PAMS {\bf92} (1984), 288--292. 



\bibitem[KBS]{KBS}
E. Kilbride, M.R. Boocock, and W.M. Stark, 
{\em Topological selectivity of a hybrid  site-specific recombination system 
with elements from Tn3  res/resolvase and bacteriophase PL lox $P$/Cre}, 
J. Mol. Biol. {\bf 289} (1999), 1219--1230.


\bibitem[KMOS]{KMOS}
P. Kronheimer, T. Mrowka, P. Ozsvath, Z. Szabo {\em Monopoles and lens space surgeries}\/, 
 Ann. of Math. (2) {\bf 165(2)}  (2007),   457--546. 


\bibitem[L]{L}
W. B. R. Lickorish, {\em  Prime knots and tangles},
Trans. Amer. Math. Soc.
{\bf 267(1)} (1981),
321--332.

\bibitem[MBM]{MBM} 
M. Mizuuchi, T. A. Baker, K. Mizuuchi,
{\em Assembly of the active form of the transposase-Mu DNA complex: a critical control point
in Mu transposition}, Cell {\bf 70(2)} (1992)  303--11.



\bibitem[PJH]{PJH}
S. Pathania, M. Jayaram, and R. Harshey, 
{\em Path of DNA within the Mu Transpososome: Transposase interaction 
bridging two Mu ends and the enhancer trap five DNA supercoils}, 
Cell {\bf109} (2002), 425--436. 

\bibitem[PJH2]{PJH2}
S. Pathania, M. Jayaram, and R. Harshey, 
{\em A unique right end-enhancer complex precedes synapsis of Mu ends: the enhancer is sequestered 
within the transpososome throughout transposition},
EMBO J. {\bf 22(14)} (2003), 3725--36. 



\bibitem[S]{S} D. Sankoff, 
{\em Rearrangements and chromosomal evolution},
 Curr Opin Genet Dev.  {\bf 13(6)} (2003) 583--7.

\bibitem[Sch]{Sch} 
M. Scharlemann, 
{\em Outermost forks and a theorem of Jaco} 
Proc. Rochester Conf., AMS Contemporary Math. Series {\bf 44} (1985),
189--193.

\bibitem[ST]{ST}
M. Scharlemann and A. Thompson, 
{\em Detecting unknotted graphs in $S^3$},
J. Differential Geom. {\bf 34} (1991) 539--560.



\bibitem[SECS]{SECS}
D. W. {Sumners}, C. {Ernst}, N.R. {Cozzarelli}, S.J. {Spengler} {\em
Mathematical analysis of the mechanisms of DNA recombination using tangles}, 
Quarterly Reviews of Biophysics {\bf 28} (1995).

\bibitem[T]{T}
A. Thompson, 
{\em A polynomial invariant of graphs in 3-manifolds},
Topology {\bf 31} (1992) 657--665.


\bibitem[VCS]{VCS}
M. {Vazquez}, S.D. {Colloms}, D.W. Sumners,
{\em Tangle analysis of Xer recombination reveals only three solutions, 
all consistent with a single three-dimensional topological pathway},
J. Mol. Biol. {\bf346} (2005), 493--504.

\bibitem[VS]{VS}
M. {Vazquez}, D. W. {Sumners}, {\em
Tangle analysis of Gin site-specific recombination}\/, 
Math. Proc. Camb. Phil. Soc. {\bf 136} (2004), 565--582.

\bibitem[VDL]{VDL}
A.A. Vetcher, A. Y. Lushnikov, J. Navarra-Madsen, R. G. Scharein, Y. L. Lyubchenko, 
I. K. Darcy, 
S. D. Levene,	
{\em DNA topology and geometry in Flp and Cre recombination}\/,
J Mol Biol. {\bf 357(4)} (2006), 1089--104.

\bibitem[Wu1]{Wu1} Y. Q. Wu, 
{\em A generalization of the handle addition theorem}, 
PAMS {\bf114} (1992), 237--242. 


\bibitem[Wu2]{Wu2} 
Y. Q. Wu, 
{\em On planarity of graphs in 3-manifolds}, Comment. Math. 
Helv. {\bf 67} (1992) 635--64.

\bibitem[Wu3]{Wu3} 
Y. Q. Wu, {\em The classification of nonsimple algebraic tangles},
Math. Ann.  {\bf 304(3)}  
(1996), 457--480.


\bibitem[YJPH]{YJPH}
Z. Yin, M. Jayaram, S. Pathania, and R. Harshey, 
{\em The Mu transposase interwraps distant DNA sites within a functional transpososome in the 
absence of DNA supercoiling},
J Biol Chem. {\bf 280(7)} (2005), 6149--56. 


\bibitem[YJH]{YJH}
Z. Yin, A. Suzuki, Z. Lou,
 M. Jayaram, and R. Harshey, 
{\em Interactions of phage Mu enhancer and termini that specify the assembly of a topologically 
unique interwrapped transpososome},
J Mol Biol. {\bf 372(2)} (2007), 382--96. 



\end{thebibliography}
